\documentclass[11pt,a4paper]{article}
\usepackage[margin=27mm,headheight=14pt]{geometry}
\usepackage{amsmath,amssymb,amsthm}
\usepackage[unicode=false,hyperfootnotes=false,colorlinks=true,linkcolor=black,citecolor=black,urlcolor=blue]{hyperref}
\hypersetup{pdftitle={Counting Rational Points on Cayley's Cubic Surface}}
\numberwithin{equation}{section}
\newtheorem{theorem}{Theorem}[section]
\newtheorem{lemma}[theorem]{Lemma}
\newtheorem{proposition}[theorem]{Proposition}
\newtheorem{corollary}[theorem]{Corollary}
\theoremstyle{definition}
\newtheorem{definition}[theorem]{Definition}
\theoremstyle{remark}
\newtheorem{remark}[theorem]{Remark}
\allowdisplaybreaks[2]
\title{Counting Rational Points on Cayley's Cubic Surface}
\author{Xiaodong Zhao}
\date{}
\begin{document}
\maketitle
\begingroup
\renewcommand{\thefootnote}{}
\footnotetext{%
}
\endgroup
\begin{abstract}
We study rational points of bounded maximum-coordinate height on the
split four-nodal Cayley cubic surface over \(\mathbb Q\), with its nine
lines removed. Writing \(U\) for this open subset and \(N_U(B)\) for
its counting function, we establish
\[
N_U(B)\sim c_{\mathrm{Peyre}}B(\log B)^6.
\]
The argument starts from the universal-torsor parametrization and
combines two counting methods over complementary ranges of the
auxiliary variables. The two-dimensional lattice count is controlled
by averages of reciprocal shortest-vector lengths and character-moment
bounds, while Poisson summation treats the complementary
parameter range.
The main contribution is the uniform control needed to combine
these methods when the size parameters of the auxiliary variables
vary independently: the estimates retain every coprimality
condition, the original height weight, and the full set of
auxiliary tuples contributing to the geometric error. Divisibility estimates without an additional power of the
height parameter control the sieve tail and justify removing the
sieve after taking the height limit. We track the derivative losses
and boundary contributions throughout the summation, and identify
the common continuous main term with Peyre's constant by computing
the effective-cone volume and the real and finite local densities.
\end{abstract}

\section{Introduction}\label{in:section}

Consider the split Cayley cubic surface
\begin{equation}\label{in:surface}
 S:\quad F(x_1,x_2,x_3,x_4)
 =x_2x_3x_4+x_1x_3x_4+x_1x_2x_4+x_1x_2x_3=0.
\end{equation}
It has four singularities, at the points
$[1:0:0:0]$, $[0:1:0:0]$, $[0:0:1:0]$ and
$[0:0:0:1]$. Moreover there are exactly $9$
lines in the surface. Six lines are
given by the vanishing of two coordinates; the other three are
\(x_i+x_j=x_k+x_l=0\), for the three partitions of four indices into pairs.
Write \(U\) for their complement.
For \(x=[x_1:\cdots:x_4]\in U(\mathbb Q)\), choose
\((x_1,\ldots,x_4)\in\mathbb Z^4\) with
\(\gcd(x_1,x_2,x_3,x_4)=1\), and put
\[
H(x)=\max_i|x_i|,\qquad
N_U(B)=\#\{x\in U(\mathbb Q):H(x)\le B\}.
\]
The height is unchanged under \(\mathbf x\mapsto-\mathbf x\). Thus
\[
2N_U(B)=\#\bigl\{\mathbf x\in\mathbb Z^4:
[\mathbf x]\in U(\mathbb Q),\ \gcd(x_1,x_2,x_3,x_4)=1,
\ \max_i|x_i|\le B\bigr\}.
\]

Let \(Y\) be the minimal resolution of \(S\). By
\cite[Example~14]{DT2007}, \(Y\) is the blowup of
\(\mathbb P^2_{\mathbb Q}\) at the six pairwise intersections
of four rational lines in general position. Hence
\(\operatorname{rank}\operatorname{Pic}(Y)=7\), and the prediction
of Manin and Peyre for the height above
\cite{BM1990,FMT1989,Peyre1995} is
\[
N_U(B)\sim c_{\mathrm{Peyre}}B(\log B)^6.
\]
Each of the nine rational lines contributes \(\asymp B^2\)
points of height at most \(B\), exceeding the predicted order
\(B(\log B)^6\); we therefore remove all nine lines.

For this surface and this open subset, Heath-Brown
\cite{HB2003} proved
\begin{equation}\label{in:HB-order}
 B(\log B)^6\ll N_U(B)\ll B(\log B)^6.
\end{equation}
Heath-Brown identifies uniform primitive lattice-point asymptotics,
rather than upper bounds, as an obstruction in his argument
\cite[Introduction and Section~6]{HB2003}.

The Cayley cubic also provides a concrete setting for explicit descent
on singular del Pezzo surfaces.
Derenthal and Tschinkel \cite[Example~14]{DT2007} describe a universal
torsor over the minimal resolution \(Y\) using thirteen affine
coordinates and nine defining relations. This torsor is not a complete
intersection. Its presentation motivates counting methods that preserve
the relations among the auxiliary variables throughout the summation.
In another arithmetic direction, Wang \cite[Theorem~2]{Wang2014}
combines the torsor parametrization with the circle method to obtain
quantitative approximation by rational points whose primitive
coordinates satisfy the condition that \(\lvert x_1x_2x_3x_4\rvert\)
has at most twelve prime factors, counted with multiplicity.
In particular, such points are Zariski dense in \(S\).

The reciprocal form of the Cayley equation also connects this surface
with Egyptian fractions. For positive integers \(n,t_1,t_2,t_3\),
\[
\frac{4}{n}=\frac1{t_1}+\frac1{t_2}+\frac1{t_3}
\quad\Longrightarrow\quad
[-n:4t_1:4t_2:4t_3]\in S(\mathbb Q).
\]
Elsholtz and Tao \cite[Remark~2.10]{ElsholtzTao2013}
relate their parametrizations for the Erd\H{o}s--Straus equation
to Heath-Brown's parametrization of the Cayley cubic.
Bright and Loughran \cite[Sections~1--2]{BrightLoughran2020}
study the associated affine surfaces, whose projective closures
are isomorphic to \(S\), through their Brauer groups and
Brauer--Manin obstructions. For \(n\ge2\), they show that there is
no Brauer--Manin obstruction to the existence of positive integral
solutions of the Erd\H{o}s--Straus equation, while identifying an
obstruction to strong approximation.
The integrality, positivity, and size restrictions in these problems
differ from the maximum-coordinate height condition used here.

These connections place the present problem within a broader study
of rational and integral points on the same underlying surface.
The refinement of \eqref{in:HB-order} to an asymptotic with Peyre's
constant measures the overall distribution of rational points of
bounded height and tests the predicted relationship between
the effective cone, the real density, and the finite local densities
\cite{Peyre1995}.

Related asymptotics concern the \(E_6\) and \(D_4\) cubics
\cite{BBD2007,LeB2014} and the Cayley ruled cubic \cite{BBS2016},
which is distinct from the four-nodal surface considered here.
Gao \cite{Gao2025} studies exceptional sets arising from
quasi-\'etale covers; the Cayley example in Section~5.1 of that work
is a nonsplit form whose minimal resolution has Picard rank~\(4\).
The surface in the present paper is split and has Picard rank~\(7\).

\begin{theorem}\label{in:main}
With the surface $S$, open subset $U$, and height $H(x)$ specified above,
\begin{equation}\label{in:asymptotic}
 N_U(B)=c_{\mathrm{Peyre}}B(\log B)^6
       +o\!\left(B(\log B)^6\right)\qquad(B\to\infty),
\end{equation}
where
\begin{equation}\label{in:constant}
 c_{\mathrm{Peyre}}=
 \frac1{80}\left(\int_0^2\frac{\log(1+s)}s\,ds\right)
 \prod_p\left(1-\frac1p\right)^7
          \left(1+\frac7p+\frac1{p^2}\right).
\end{equation}
The Euler product converges absolutely to a positive number.
\end{theorem}

The passage from Heath-Brown's bounds \cite{HB2003} to
Theorem~\ref{in:main} requires uniform asymptotic estimates
whose error terms remain summable as the size parameters
of the auxiliary variables vary independently. The underlying tools include
primitive plane-point bounds, character-moment estimates
\cite{CS2010,HB1995}, and bilinear and trilinear counting results
\cite{Reuss2015}. Our contribution is to obtain compatible uniform
estimates from these tools while retaining the original height and all
coprimality conditions, and to show that their complementary ranges
cover the count up to negligible boundary contributions.

There are two separate errors in the two-dimensional count.
Propositions~\ref{ge:four-proposition} and~\ref{ge:reuss-proposition}
control the geometric error by averaging reciprocal shortest-vector
lengths over the full set of auxiliary tuples, including tuples that
produce no admissible integer point. Theorem~\ref{ar:main} then transfers
the resulting geometric main term to a continuous measure. Its
character estimates are uniform in independently varying
interval parameters and retain the joint weights arising from
recovery of the eliminated variables.
The dependence of these weights on the modulus is handled before
the modulus average is taken.

The need for complementary estimates is visible even on symmetric
blocks. In the notation of
\eqref{in:scales}--\eqref{in:short-scales}, if
\(Z_e\asymp B^{1/4}\) for all six index pairs, then
\[
 P_0\asymp B^{3/2},\qquad Y_i\asymp B^{1/12},\qquad
 Q\asymp B^{1/3},\qquad J\asymp B^{1/2},\qquad
 P_y\asymp B^{1/3}.
\]
The two alternatives for the nonquadratic-character contribution in
\eqref{ar:relative-error} consequently satisfy
\[
 \frac{Q}{\sqrt{P_y}}\asymp B^{1/6},\qquad
 \frac{Q}{\sqrt{JP_y}}+\sqrt{\frac QJ}\ll B^{-1/12}.
\]
The second alternative supplies a saving where the first does not.
The term \(P_0/B\) in the Poisson error \eqref{po:frequency}
also gives no saving on these blocks. Thus the uniform
two-dimensional analysis must include blocks with \(P_0>B\).

In the complementary parameter range, Poisson summation on the original
torsor lattice yields the same continuous measure and local factors
(Theorem~\ref{po:main}). A separate divisibility argument gives the
sieve-tail bound of Theorem~\ref{ta:weak-tail} without a
\(B^\varepsilon\) loss, allowing the sieve threshold to tend to
infinity after the height limit. The
\hyperref[ta:proof-strategy]{preliminary roadmap} specifies the
assignment of parameter ranges to the two methods and the dependencies between
the estimates; Section~\ref{gl:section} verifies that all errors and
derivative losses remain summable.

Section~\ref{ta:section} establishes the divisibility, sieve-tail,
and boundary estimates. Section~\ref{ge:section} treats the
two-dimensional geometric error, and Section~\ref{ar:section}
establishes its arithmetic main term.
Section~\ref{po:section} treats the range assigned to Poisson summation.
Section~\ref{co:section} computes the common volume and identifies
Peyre's constant. Finally, Section~\ref{gl:section} assembles the local
estimates and removes the truncations in the specified order.

\section{Uniform sieve tails and boundary estimates}\label{ta:section}

We first prove a divisibility estimate for canonical points and a sieve-tail bound for weak representations, uniform in \(B\) and \(R\).
We then estimate logarithmic boundary regions and proximity to a selected excluded line.

\subsection{Notation and preliminary identities}

An index pair is an unordered pair
\(e=\{i,j\}\subset\{1,2,3,4\}\) with \(i\ne j\).
A pairing \(ij|kl\) is a partition of \(\{1,2,3,4\}\) into the
two pairs \(\{i,j\}\) and \(\{k,l\}\). The three pairings are
\(12|34\), \(13|24\), and \(14|23\); the two pairs in each
pairing are called complementary. Let
\[
z_{ij}=z_{ji}\in\mathbb Z_{>0}\quad(i<j),
\qquad y_i\in\mathbb Z\setminus\{0\}.
\]
Set
\begin{equation}\label{in:torsor-monomials}
P(z)=\prod_{i<j}z_{ij},\qquad
A_i(z)=\prod_{\substack{j<k\\j,k\ne i}}z_{jk},\qquad
C_i(z)=\prod_{j\ne i}z_{ij}.
\end{equation}
Then
\[
A_iC_i=P,\qquad \prod_iA_i=\prod_iC_i=P^2.
\]
By \cite[Section~2, Lemma~1]{HB2003}, the parametrization is
\begin{equation}\label{in:torsor}
\sum_iA_i(z)y_i=0,\qquad
x_i=C_i(z)\prod_{j\ne i}y_j,
\end{equation}
with
\begin{equation}\label{in:coprimality}
\gcd(z_e,z_f)=1\;(e\ne f),\qquad
\gcd(y_i,y_j)=1\;(i<j),\qquad
\gcd(y_i y_j,z_{ij})=1\;(i<j),
\end{equation}
and
\[
A_1(z)y_1+A_j(z)y_j\ne0\qquad(j=2,3,4).
\]
The map \((z,y)\mapsto[x_1:x_2:x_3:x_4]\) is two-to-one onto
\(U(\mathbb Q)\), with preimages \((z,y)\) and \((z,-y)\)
corresponding to the primitive representatives \(x\) and \(-x\).
When coprimality conditions are relaxed, use the unreduced height
\begin{equation}\label{in:raw-height}
H_{\mathrm{raw}}(z,y)
=\max_i\left|C_i(z)\prod_{j\ne i}y_j\right|,
\end{equation}
without cancelling a common factor from the coordinates.

By \cite[Section~2, equation~(2.16)]{HB2003}, for complementary
pairs \(ij|kl\),
\begin{align}
z_{ik}z_{il}y_j+z_{jk}z_{jl}y_i
&=z_{ij}v_{ij},
\qquad v_{ij}+v_{kl}=0,
\notag\\
A_i(z)y_i+A_j(z)y_j
&=z_{ij}z_{kl}v_{ij},
\qquad v_{ij}\in\mathbb Z\setminus\{0\}.
\label{in:q-identity}
\end{align}

For \(B\ge2\), consider a dyadic block
\(Z_e\le z_e<2Z_e\), where \(Z_e=2^m\) for an integer \(m\ge0\).
Write \(\mathbf Z=(Z_e)_e\), and set
\[
\begin{aligned}
(a,b,c,d,u,v)
&=(z_{12},z_{34},z_{14},z_{23},z_{13},z_{24}),\\
(Z_a,Z_b,C,D,U,V)
&=(Z_{12},Z_{34},Z_{14},Z_{23},Z_{13},Z_{24}).
\end{aligned}
\]
After permuting the indices, put
\begin{equation}\label{in:scales}
\begin{gathered}
A=Z_aZ_b\ge H=CD\ge J=UV,\qquad P_0=AHJ,\\
T=(BP_0)^{1/3},\qquad Q=T/A,\qquad Q'=T/H.
\end{gathered}
\end{equation}
For a pairing \(ij|kl\), the corresponding product of size
parameters is \(Z_{ij}Z_{kl}\); the three such products are
\(A,H,J\).
Here \(H=CD\); the height function is written \(H(x)\). Set
\begin{equation}\label{in:short-scales}
Y_i=\frac{T}{A_i(\mathbf Z)},\qquad
P_y=\prod_iY_i=\frac{T^4}{P_0^2},\qquad
F_0=\frac{T^2}{P_0}=\sqrt{P_y}.
\end{equation}
For fixed \(M\ge1\), the formulas for \(x_i\)
in \eqref{in:torsor} give
\[
\begin{gathered}
B/M\le|x_i|\le B\quad(1\le i\le4)
\quad\Longrightarrow\\
|A_i(z)y_i|
=\frac{\bigl(P(z)\prod_j|x_j|\bigr)^{1/3}}{|x_i|}
\asymp_M T,\qquad |y_i|\asymp_M Y_i.
\end{gathered}
\]
We shall use the following identities:
\begin{equation*}
B=\frac{HJP_y}{Q},\quad
U^2Y_2Y_4=V^2Y_1Y_3=QQ',\quad
Z_{ij}\sqrt{Y_kY_l}=\sqrt{F_0Z_{ij}Z_{kl}}
\quad(ij|kl).
\end{equation*}

Divisor and logarithmic factors whose arguments have polynomial
size in \(B\) may be absorbed into \(B^\varepsilon\), using smaller
auxiliary exponents when several factors occur. Unbounded frequency
sums are estimated before this absorption. The final sieve-tail
bounds with fixed threshold contain no \(B^\varepsilon\) loss.
All smooth-weight norms refer to the normalized coordinates
specified in the relevant estimate. Unless otherwise stated,
implied constants may depend on \(\varepsilon\), fixed supports,
derivative orders, and fixed height ratios, but are uniform in
\(B\) and the individual dyadic scales; any dependence on sieve
parameters is stated explicitly.

\phantomsection
\label{ta:proof-strategy}
We record the division of the parameter ranges and the dependencies used in the
proof. The statements in this roadmap summarize the estimates proved
below; they are not additional inputs to those proofs.
Fix the coordinate-height ratio \(M\) and the height-smoothing
parameter \(h\), as in Section~\ref{gl:section}. For
\(0<\delta<1/100\), put
\[
 (Q_1,Q_2,Q_3)=\left(\frac TA,\frac TH,\frac TJ\right).
\]
After truncating the coordinate heights, retain the blocks satisfying
\[
 \min\{\min_e Z_e,\min_iY_i,\min_jQ_j\}\ge B^\delta,
 \qquad
 P_0\le B^{1-\delta}\ \text{or}\ P_0\ge B^{1+\delta},
\]
and remove the transition band
\[
 B^{-\delta}<\frac{HJ}{T}<B^\delta.
\]
Propositions~\ref{ta:height-tail} and~\ref{ta:boundary} control the
height tail and the discarded logarithmic ranges, respectively.
On the retained blocks the assignment is exhaustive:
\begin{center}
\small
\renewcommand{\arraystretch}{1.25}
\begin{tabular}{@{}p{0.48\textwidth}p{0.48\textwidth}@{}}
\hline
Range & Counting method and geometric input\\
\hline
\(P_0\le B^{1-\delta}\), \(A\le F_0B^{3\delta/4}\)
& Poisson summation, Theorem~\ref{po:main}\\
\(P_0\le B^{1-\delta}\), \(A>F_0B^{3\delta/4}\)
& Two-dimensional count, Proposition~\ref{ge:four-proposition}\\
\(P_0\ge B^{1+\delta}\), \(A/J\ge B^\delta\)
& Two-dimensional count, Proposition~\ref{ge:four-proposition}\\
\(P_0\ge B^{1+\delta}\), \(A/J<B^\delta\)
& Two-dimensional count, Proposition~\ref{ge:reuss-proposition}\\
\hline
\end{tabular}
\end{center}
All three two-dimensional ranges use Theorem~\ref{ar:main} to
evaluate the geometric main terms. Outside the deleted \(HJ/T\)
band, one of its two nonquadratic-character bounds gives a saving.
In the Poisson range, the scale identity
\[
 \frac{Z_{ij}}{F_0}
 =\frac{\sqrt{Z_{ij}Z_{kl}/F_0}}{\sqrt{Y_kY_l}}
 \le B^{-5\delta/8}\qquad(ij|kl)
\]
verifies the additional bound on \(Z_{ij}/F_0\) required
for Poisson summation.
Section~\ref{gl:geometric-core} proves the geometric estimates
for these four ranges, and the following subsection treats their
arithmetic errors.

For the two-dimensional count, select the pairing \(ij|kl\)
for which \(Z_{ij}Z_{kl}\) is largest, and relabel the indices
so that this pairing is \(12|34\). Put \(q=|v_{12}|\).
Theorem~\ref{ta:near-line} removes \(q\ll B^{-\kappa}Q\);
the remaining range is divided into smooth bands \(q\asymp sQ\),
with \(B^{-\kappa}\le s\ll1\).
For \(0<\kappa,\varepsilon<\delta/100\), the local errors obtained
in Section~\ref{gl:section} are
\[
\begin{array}{c|c}
\text{contribution}&\text{bound}\\ \hline
\text{two-dimensional geometric error}
 &B^{1+\varepsilon+2\kappa-\delta/6}\\
\text{two-dimensional arithmetic error}
 &B^{1+\varepsilon+7\kappa-\delta/2}\\
\text{Poisson error}
 &B^{1-5\delta/16+\varepsilon/2}.
\end{array}
\]
The first two bounds are per block and band; the third is per block
with the finite sieve threshold \(R\) fixed. There are
\(O(\log^6(2B))\) blocks in the variables \(z_e\) and \(O(\log(2B))\) bands per
two-dimensional block, so all three totals are
\(o(B(\log B)^6)\).

The proof dependencies are as follows. The reciprocal-shortest-vector
estimates of Section~\ref{ge:section} control the geometric error in
Proposition~\ref{ar:long}; Theorem~\ref{ar:main} evaluates its main
term, yielding Corollary~\ref{ar:count}.
The frequency estimate~\eqref{po:frequency} is proved directly from
the original torsor lattice, independently of
Theorem~\ref{ta:near-line}, and is therefore available in the proof
of that theorem. Separately, the divisibility estimates imply
Theorem~\ref{ta:weak-tail}; this tail bound, through
Corollary~\ref{po:free-tail}, permits removal of the finite sieve in
Corollary~\ref{po:core-main}. The forward use of the frequency
estimate and the subsequent removal of the sieve thus involve
distinct arguments.

Finally, the two methods give the same continuous measure, whose
volume is evaluated in Section~\ref{co:section}.
Section~\ref{gl:section} assembles the counts and takes limits in the
order
\[
 B\longrightarrow\infty,\quad R\longrightarrow\infty,\quad
 \delta\longrightarrow0,\quad h\longrightarrow0,\quad
 M\longrightarrow\infty.
\]
The parameters \(\kappa\) and \(\varepsilon\) are chosen after
\(\delta\) and before the first limit. In particular, the constants
in the finite-sieve Poisson estimate need not be uniform in \(R\).
The power savings on the retained blocks are thus fixed-parameter
estimates; the final asymptotic follows by removing the truncations
through these ordered limits.

\subsection{Canonical and weak parameters}

A \emph{canonical representation} is a tuple \((z,y)\) satisfying
\eqref{in:torsor}, \eqref{in:coprimality}, and
\(A_1(z)y_1+A_j(z)y_j\ne0\) for \(j=2,3,4\).
Throughout this section we choose \(y_1>0\), removing the sign multiplicity.

\begin{definition}\label{ta:weak-definition}
A \emph{weak representation} is a tuple \((z,y)\) satisfying all
the canonical conditions above, including \(y_1>0\), except that
\(\gcd(y_i,y_j)=1\) for \(i\ne j\) is not required. With \(H_{\mathrm{raw}}\) as in \eqref{in:raw-height}, define
\[
\mathcal N_{\mathrm{weak}}(B)
=
\#\bigl\{
(z,y)\text{ weak}: H_{\mathrm{raw}}(z,y)\le B
\bigr\},
\]
and
\[
\mathcal T_{>R}(B)
=
\#\left\{
(z,y)\text{ weak}:
\begin{array}{l}
H_{\mathrm{raw}}(z,y)\le B,\\
\exists\,p>R\text{ prime},\quad
\exists\,1\le i<j\le4:\\
p\mid y_i,\quad p\mid y_j
\end{array}
\right\}.
\]
\end{definition}

For a fixed canonical tuple \((z,y)\), the tuples
\[
 (z,gy),\qquad g\in\mathbb Z_{\ge1},\quad \gcd(g,P(z))=1,
\]
are weak and represent the same rational point \(x\), whereas
\(H_{\mathrm{raw}}(z,gy)=g^3H(x)\).
Thus the unreduced height is necessary for a finite count.

For a weak representation \((z,y)\), write
\[
x_{\mathrm{raw}}
=(x_{\mathrm{raw},1},\ldots,x_{\mathrm{raw},4}),
\qquad
x_{\mathrm{raw},i}
=C_i(z)\prod_{j\ne i}y_j.
\]
These coordinates are taken without cancelling a common factor.

\begin{lemma}\label{ta:normalization}
Every weak representation determines a canonical representation
\((z^*,y^*)\) and positive integers \(g,h\), with \(h\mid P(z^*)\), such that
\begin{equation}\label{ta:normalization-height}
x_{\mathrm{raw}}=g^3h x^*,
\qquad
H_{\mathrm{raw}}=g^3h H(x^*),
\end{equation}
where
\[
x_i^*=C_i(z^*)\prod_{j\ne i}y_j^*.
\]
For fixed \((z^*,y^*),g,h\), at most one weak representation has these
normalization data. Moreover, for every prime \(p\),
\[
p\mid y_i,\quad p\mid y_j,\quad i\ne j
\qquad\Longrightarrow\qquad p\mid gh.
\]
\end{lemma}

\begin{proof}
Set $g=\gcd(|y_1|,|y_2|,|y_3|,|y_4|).$
No prime divides three of the integers \(y_i/g\).
Indeed, if \(p\) divides \(y_i/g,y_j/g,y_k/g\), then every
pair \(\{a,b\}\), with \(1\le a<b\le4\), contains at least one
of \(i,j,k\). Hence the coprimality conditions
\(\gcd(y_a y_b,z_{ab})=1\) imply \(p\nmid z_{ab}\) for every
\(a<b\), and therefore \(p\nmid P(z)\). The torsor equation then yields
\[
A_l(z)\frac{y_l}{g}
\equiv
-\sum_{r\ne l}A_r(z)\frac{y_r}{g}
\equiv0\pmod p,
\qquad
p\nmid A_l(z).
\]
Thus \(p\mid y_l/g\), contradicting
\(\gcd(|y_1|/g,\ldots,|y_4|/g)=1\).

For complementary pairs \(\{i,j\}\) and \(\{k,l\}\), define
\[
\begin{gathered}
h_{ij}=\gcd\!\left(\frac{|y_k|}{g},\frac{|y_l|}{g}\right),
\qquad
\mathbf h=(h_{ij})_{i<j},
\qquad
h=P(\mathbf h),\qquad z_{ij}^*=z_{ij}h_{ij},
\qquad
y_i^*=\frac{y_i}{gA_i(\mathbf h)}.
\end{gathered}
\]
The preceding observation and the conditions
\(\gcd(y_r y_s,z_{rs})=1\) for \(r<s\) imply
\[
\begin{gathered}
\gcd(h_e,h_f)=\gcd(z_e,h_f)=1
\;(e\ne f)\quad\text{and}\quad
\gcd\!\left(z_{ij}h_{ij},\frac{|y_i|}{g}\right)=1.
\end{gathered}
\]
Indeed, a prime dividing \(h_{ij}\) divides exactly the two quotients
\(y_k/g,y_l/g\); it can divide no original variable \(z_e\) other than \(z_{ij}\).

Since the factors below are pairwise coprime,
\[
A_i(\mathbf h)
=
\prod_{j\ne i}
\gcd\!\left(\frac{|y_i|}{g},\frac{|y_j|}{g}\right)
\ \Bigm|\ \frac{|y_i|}{g}.
\]
Consequently, \(y_i^*\in\mathbb Z\setminus\{0\}\), \(y_1^*>0\), and
\[
\begin{gathered}
\gcd(z_e^*,z_f^*)=1\quad(e\ne f),
\qquad
\gcd(|y_i^*|,z_{ij}^*)=1,\\
\gcd(|y_i^*|,|y_j^*|)
\ \Bigm|\
\gcd\!\left(
\frac{|y_i|}{g h_{kl}},
\frac{|y_j|}{g h_{kl}}
\right)
=1
\qquad(i\ne j).
\end{gathered}
\]
Multiplicativity of \(A_i\) gives
\[
\begin{gathered}
A_i(z^*)y_i^*=\frac{A_i(z)y_i}{g},
\qquad
\sum_i A_i(z^*)y_i^*=0,\\
A_1(z^*)y_1^*+A_j(z^*)y_j^*
=
\frac{A_1(z)y_1+A_j(z)y_j}{g}
\ne0
\qquad(j=2,3,4).
\end{gathered}
\]
Thus \((z^*,y^*)\) is canonical.
By \cite[Section~2, Lemma~1]{HB2003}, the corresponding vector
\(x^*\) is primitive and represents a point of \(U(\mathbb Q)\).

Using \(A_i(\mathbf h)C_i(\mathbf h)=h\) and
\(\prod_i A_i(\mathbf h)=h^2\), we obtain
$\prod_{j\ne i}A_j(\mathbf h)
=\frac{h^2}{A_i(\mathbf h)}
=hC_i(\mathbf h).$
Hence
\[
\begin{aligned}
x_{\mathrm{raw},i}
=C_i(z)\prod_{j\ne i}
   \bigl(gA_j(\mathbf h)y_j^*\bigr)
=g^3h\,C_i(z)C_i(\mathbf h)\prod_{j\ne i}y_j^*
=g^3h\,C_i(z^*)\prod_{j\ne i}y_j^*
 =g^3h x_i^*.
\end{aligned}
\]
Since \(g^3h>0\), this proves the height identity in
\eqref{ta:normalization-height}. Moreover,
\[
P(z^*)=P(z)P(\mathbf h)=P(z)h,
\qquad h\mid P(z^*).
\]

Finally, \(h=\prod_e h_e\), \(h_e\mid z_e^*\), and the \(z_e^*\) are
pairwise coprime. Therefore
\begin{equation*}
h_e=\gcd(h,z_e^*),
\qquad
z_e=\frac{z_e^*}{h_e},
\qquad
y_i=gy_i^*A_i(\mathbf h).
\end{equation*}
These formulas determine at most one weak representation from
\((z^*,y^*),g,h\). For complementary pairs \(\{i,j\}\) and \(\{k,l\}\),
\[
p\mid y_i,\quad p\mid y_j,\quad p\nmid g
\quad\Longrightarrow\quad
p\mid h_{kl}
\quad\Longrightarrow\quad
p\mid h.
\]
Thus every prime dividing two distinct \(y\)-coordinates divides \(gh\).
\end{proof}

\subsection{A divisibility-weighted product estimate}

We first record the primitive plane estimate used throughout this
section. If \(\mathbf a\in\mathbb Z^3\) is primitive and \(H_i>0\),
Heath-Brown \cite[Lemma~6]{HB2003} gives
\begin{equation}\label{ta:plane}
 \#\{\mathbf n\in\mathbb Z^3_{\mathrm{prim}}:
       \mathbf a\cdot\mathbf n=0, |n_i|\le H_i\}
 \ll 1+\frac{H_1H_2H_3}{\max_i H_i|a_i|}.
\end{equation}
In every application below the normal vector and the counted vector
are primitive. In particular, the constant term in this estimate is
never suppressed.

\begin{lemma}\label{ta:product-lemma}
Let \(K_1,\ldots,K_7\ge \tfrac12\) and \(D\in\mathbb Z_{>0}\). Put
\[
\begin{gathered}
V_1=K_1K_2K_3,\quad
V_2=K_4K_5K_6,\quad
V=V_1V_2,\quad
J=\max(V_1,V_2),\quad
\mathcal R=\frac{\max(V_1,V_2)}{\min(V_1,V_2)}.
\end{gathered}
\]
For either fixed choice of sign, let \(N\) count the integer solutions of
\[
 n_1n_2n_3\pm n_4n_5n_6=n_7n_8,
 \qquad K_i<n_i\le 2K_i\quad(1\le i\le7),
\]
with $n_8\ne0$,
\(\gcd(n_1n_2n_3,n_4n_5n_6)=1\) and \(D\mid n_1\).
Then
\begin{equation}\label{ta:weighted-product}
N\ll VD^{-5/6}\mathcal R^{1/3},
\end{equation}
with an absolute implied constant. The same bound holds if any
of the intervals \((K_i,2K_i]\) is replaced by \([K_i,2K_i)\).
\end{lemma}

\begin{proof}
The proof uses only the containing closed intervals \([K_i,2K_i]\),
so it covers both endpoint conventions.
Assume \(N>0\). Write \(n_1=Dm_1\), and set
\[
L_1=\max(K_1/D,K_2,K_3),
\qquad
L_2=\max(K_4,K_5,K_6).
\]
Nonemptiness gives
\[
\frac{K_1}{D}\ge\frac12,
\qquad
L_1^3\ge\frac{V_1}{D},
\qquad
L_2^3\ge V_2,
\qquad
L_1L_2\ge(V/D)^{1/3}.
\]

Choose a variable of scale \(L_1\) from \((m_1,n_2,n_3)\),
and one of scale \(L_2\) from \((n_4,n_5,n_6)\), denoting them by
\(x,z\), respectively. The remaining four variables have
$O\!\left(\frac{V}{DL_1L_2}\right)$
possible values. For each such choice, write
\[
\begin{gathered}
n_1n_2n_3=ax,\qquad n_4n_5n_6=bz,\qquad
ax\pm bz=n_7n_8,\\
L_1\le x\le2L_1,\qquad L_2\le z\le2L_2,\qquad
aL_1\asymp V_1,\qquad bL_2\asymp V_2.
\end{gathered}
\]
Discard choices admitting no solution. The coprimality condition implies
\[
\gcd(ax,bz)=1
\quad\Longrightarrow\quad
\gcd(a,b)=\gcd(x,z)=1.
\]
Moreover,
\[
|n_7n_8|
\le 8(V_1+V_2)\le16J,
\qquad
0<|n_8|\le\frac{16J}{K_7}.
\]

Apply \eqref{ta:plane}, first fixing \(n_7\), then fixing
\(n_8\ne0\). The corresponding normal and point vectors are
\[
\begin{array}{c|c|c}
\text{fixed variable}&\text{normal vector}&\text{point vector}\\ \hline
n_7&(a,\pm b,-n_7)&(x,z,n_8)\\
n_8&(a,\pm b,-n_8)&(x,z,n_7).
\end{array}
\]
Both vectors are primitive in each application, and the denominator
in the cited estimate is
\[
\gg\max(aL_1,bL_2)\asymp J.
\]
Using box bounds \((2L_1,2L_2,16J/K_7)\) and
\((2L_1,2L_2,2K_7)\), respectively, gives
\[
\begin{aligned}
N
\ll
\frac{V}{DL_1L_2}
\sum_{K_7\le n_7\le2K_7}
\left(1+\frac{L_1L_2}{K_7}\right)
\ll
\frac{V}{DL_1L_2}\left(K_7+L_1L_2\right),\\[1ex]
N
\ll
\frac{V}{DL_1L_2}
\sum_{0<|n_8|\le16J/K_7}
\left(1+\frac{L_1L_2K_7}{J}\right)
\ll
\frac{V}{DL_1L_2}\left(\frac{J}{K_7}+L_1L_2\right).
\end{aligned}
\]
The second sum includes both signs of \(n_8\), and
\[
\#\left\{n_8\in\mathbb Z:
0<|n_8|\le\frac{16J}{K_7}\right\}
\le\frac{32J}{K_7}.
\]
Taking the smaller of the two bounds yields
\begin{equation}\label{ta:weighted-product-intermediate}
\begin{aligned}
N
\ll
\frac{V}{DL_1L_2}
\left(L_1L_2+\min\{K_7,J/K_7\}\right)
\le
\frac VD+\frac{V\sqrt J}{DL_1L_2}
\le
\frac VD+(V/D)^{2/3}J^{1/2}.
\end{aligned}
\end{equation}

Finally,
\[
D\le2K_1\le8V_1\le8J,
\qquad
J^2=V\mathcal R,
\qquad D,\mathcal R\ge1.
\]
Therefore
\[
\begin{aligned}
(V/D)^{2/3}J^{1/2}
=
VD^{-5/6}\mathcal R^{1/3}
\left(\frac DJ\right)^{1/6}
\le
8^{1/6}VD^{-5/6}\mathcal R^{1/3}\quad\text{and}\quad
VD^{-1}
\le VD^{-5/6}\mathcal R^{1/3}.
\end{aligned}
\]
Substitution into \eqref{ta:weighted-product-intermediate}
proves \eqref{ta:weighted-product}.
\end{proof}

\subsection{Uniform divisibility in the canonical parametrization}

For \(x\in U(\mathbb Q)\), let \(z_e(x)\) be the positive parameter indexed by \(e\) in its canonical representation; it is unchanged by \(y\mapsto-y\).
For a fixed index pair \(e\), an integer \(D\ge1\), and a real number \(X\ge1\), define
\[
 N_{D,e}(X)=\#\{x\in U(\mathbb Q):H(x)\le X,\ D\mid z_e(x)\}.
\]

\begin{proposition}\label{ta:single-edge}
Uniformly for \(D\in\mathbb Z_{\ge1}\), \(X\ge1\), and every index
pair \(e\),
\begin{equation}\label{ta:single-edge-bound}
N_{D,e}(X)\ll X\log^6(2X)D^{-5/48}.
\end{equation}
The implied constant is absolute.
\end{proposition}

\begin{proof}
Decompose into dyadic blocks
$X_i<|x_i|\le2X_i$ and $Z_{ij}\le z_{ij}<2Z_{ij},$
with all scales powers of two, \(X_i\ge1/2\), and \(Z_{ij}\ge1\).
Signs and index permutations contribute only absolute factors.
By \cite[Section~6, equations~(6.3)--(6.5)]{HB2003}, we may assume
$X_1\ge X_2\ge X_3\ge X_4$ and $X_3\le6X_4.$
Put
\begin{equation}\label{ta:dyadic-notation}
\begin{gathered}
\mathcal X=X_1X_2X_3X_4,\quad
P_0=\prod_{i<j}Z_{ij},\quad
W=X_1,\quad M=\frac{W}{X_4},\quad
\sigma=1+\log(2+W)
             \left(\frac{P_0^2}{\mathcal X}\right)^{1/9}.
\end{gathered}
\end{equation}
The same reference gives
\begin{equation}\label{ta:height-scales}
|y_i|\asymp Y_i,\qquad
Y_i=\frac{(\mathcal XP_0)^{1/3}}
          {X_iA_i(\mathbf Z)}.
\end{equation}
All comparison constants are absolute. Each range
\(|y_i|\asymp Y_i\) is covered by boundedly many dyadic intervals.
With \(X,e\) fixed, write \(\mathbf X=(X_i)_{i=1}^4\), and let
\(N_D=N_D(\mathbf X;\mathbf Z)\) be the contribution of the dyadic block
\((\mathbf X,\mathbf Z)\) to \(N_{D,e}(X)\).

Choose a pairing from $\{13,24\}, \{14,23\}$
which does not contain \(e\). Write it as \(\{ij,kl\}\), and put
\(R=Z_{ij}Z_{kl}\). Then $e\in\{ik,il,jk,jl\}.$
We use the identities
\cite[Section~2, equations~(2.16)--(2.17)]{HB2003}:
\begin{align}\label{ta:aux-linear}
z_{ik}z_{il}y_j+z_{jk}z_{jl}y_i
 =z_{ij}v_{ij}
\quad\text{and}\quad
v_{ij}v_{ik}
 =z_{il}^2y_jy_k-z_{jk}^2y_iy_l.
\end{align}
Here \(v_{ij}\ne0\), and canonical coprimality gives $\gcd(z_{ik}z_{il}|y_j|,\,
     z_{jk}z_{jl}|y_i|)=1.$
Moreover,
\[
\frac{Z_{ik}Z_{il}Y_j}{Z_{jk}Z_{jl}Y_i}
=\frac{X_i}{X_j},
\qquad
\max\!\left(\frac{X_i}{X_j},\frac{X_j}{X_i}\right)\le M.
\]
Applying Lemma~\ref{ta:product-lemma}, with \(z_e\) as the factor
divisible by \(D\), bounds the number of admissible tuples
\((z_{ik},z_{il},z_{jk},z_{jl},y_i,y_j,z_{ij},v_{ij})\) by
\begin{equation}\label{ta:first-eight}
\ll D^{-5/6}M^{1/3}
     Z_{ik}Z_{il}Z_{jk}Z_{jl}Y_iY_j.
\end{equation}

For each such tuple, apply the primitive-plane estimate
\eqref{ta:plane}, as in
\cite[Section~6, calculation preceding~(6.7)]{HB2003}, to
$(v_{ij},-z_{il}^2y_j,z_{jk}^2y_i)
\cdot(v_{ik},y_k,y_l)=0.$
Its normal and point are primitive, since $\gcd(z_{il}^2|y_j|,z_{jk}^2|y_i|)=1$ and
$\gcd(|y_k|,|y_l|)=1.$
The auxiliary identity \eqref{ta:aux-linear} also gives
\[
|v_{ik}|\ll
\frac{(\mathcal XP_0)^{1/3}}{Z_{ik}Z_{jl}}
\max(X_i^{-1},X_k^{-1}).
\]
Use this box bound together with \(|y_k|\ll Y_k\) and
\(|y_l|\ll Y_l\).
The remaining variable is uniquely determined by
$
z_{kl}
=-\frac{z_{il}z_{jl}y_k+z_{ik}z_{jk}y_l}{v_{ij}}.$
Multiplying the resulting uniform bound by \eqref{ta:first-eight}
yields
\begin{equation*}
\begin{aligned}
N_D\ll D^{-5/6}M^{1/3}
\bigg\{
\frac{(\mathcal XP_0)^{2/3}}{Z_{kl}^2X_iX_j}
+\max(X_i^{-1},X_k^{-1})
       \min(X_iX_l,X_jX_k)
\bigg\}.
\end{aligned}
\end{equation*}
Both the constant and volume terms of the plane estimate are retained.

Use the following pairs of index choices:
\[
\begin{array}{c|c|c}
R&(i,j,k,l)&
\max(X_i^{-1},X_k^{-1})\min(X_iX_l,X_jX_k)\\ \hline
Z_{14}Z_{23}&(1,4,2,3),\ (2,3,1,4)&X_4\\
Z_{13}Z_{24}&(1,3,2,4),\ (2,4,1,3)&\le X_3.
\end{array}
\]
For the corresponding first terms,
\[
\begin{aligned}
\min\!\left\{
\frac{(\mathcal XP_0)^{2/3}}{Z_{kl}^2X_iX_j},
\frac{(\mathcal XP_0)^{2/3}}{Z_{ij}^2X_kX_l}
\right\}
\le
\frac{(\mathcal XP_0)^{2/3}}
     {Z_{ij}Z_{kl}\mathcal X^{1/2}}
=
\frac{\mathcal X^{1/6}P_0^{2/3}}{R}.
\end{aligned}
\]
Consequently,
\begin{equation*}
N_D\ll D^{-5/6}M^{1/3}
\left(\frac{\mathcal X^{1/6}P_0^{2/3}}{R}+X_3\right).
\end{equation*}

Relabel the indices so that $e=ik$ and $R=Z_{il}Z_{jk}.$
This is possible because the chosen pairing does not contain \(e\).
The eight-tuple estimate in
\cite[Section~6, immediately after~(6.7)]{HB2003}
bounds the number of admissible tuples $(v_{ij},v_{ik},z_{il},z_{jk},y_i,y_j,y_k,y_l)$
by
\begin{equation}\label{ta:HB-eight}
\begin{aligned}
\ll
\sigma R\mathcal X^{1/3}P_0^{-2/3}
\max\!\left\{
\frac{X_iX_l}{X_jX_k},
\frac{X_jX_k}{X_iX_l}
\right\}^{1/4}\le
\sigma M^{1/2}R\mathcal X^{1/3}P_0^{-2/3}.
\end{aligned}
\end{equation}
The substitutions in that estimate are
\[
\prod_rY_r=\mathcal X^{1/3}P_0^{-2/3},
\qquad
\frac{Z_{il}^2Y_jY_k}{Z_{jk}^2Y_iY_l}
=\frac{X_iX_l}{X_jX_k},
\qquad
\left(\prod_rY_r\right)^{-1/3}
=\left(\frac{P_0^2}{\mathcal X}\right)^{1/9}.
\]

Retain only eight-tuples admitting a completion with \(D\mid z_{ik}\).
Writing \(z_{ik}=Dw\), equation \eqref{ta:aux-linear} becomes
\[
Dz_{il}y_jw+z_{jk}y_i z_{jl}-v_{ij}z_{ij}=0.
\]
For every admissible completion,
$\gcd(Dz_{il}|y_j|,z_{jk}|y_i|)=1$ and
$\gcd(w,z_{jl})=1.$
Thus the normal and point
\[
(Dz_{il}y_j,z_{jk}y_i,-v_{ij}),
\qquad
(w,z_{jl},z_{ij})
\]
are primitive. Applying \eqref{ta:plane} with box scales
\(Z_{ik}/D,Z_{jl},Z_{ij}\), up to absolute factors, gives
\[
\begin{aligned}
\#\{(w,z_{jl},z_{ij})\}
\ll
1+\frac{Z_{ik}Z_{jl}Z_{ij}}
 {D\max(Z_{ik}Z_{il}Y_j,Z_{jl}Z_{jk}Y_i)}
=
1+\frac{Z_{ik}Z_{jl}Z_{ij}Z_{kl}}
 {D(\mathcal XP_0)^{1/3}}\min(X_i,X_j).
\end{aligned}
\]
The final variable is recovered by the formula above. Since
$
Z_{ik}Z_{jl}Z_{ij}Z_{kl}=\frac{P_0}{R}$ and
$\min(X_i,X_j)\le W,$
equation \eqref{ta:HB-eight} implies
\begin{equation*}
N_D\ll \sigma M^{1/2}
\left(R\mathcal X^{1/3}P_0^{-2/3}+\frac WD\right).
\end{equation*}

Using \(\mathcal X\le W^4\), \(M\ge1\), and \(\sigma\ge1\), the two
counting orders give
\[
\begin{aligned}
\frac{N_D}{\sigma M^{1/2}}
&\ll
\min\!\left\{
\frac{D^{-5/6}W^{2/3}P_0^{2/3}}R,\,
RW^{4/3}P_0^{-2/3}
\right\}
+WD^{-5/6}+WD^{-1}\\
&\le
WD^{-5/12}+WD^{-5/6}+WD^{-1}
\le3WD^{-5/12}.
\end{aligned}
\]
Hence
\begin{equation}\label{ta:divisibility-block-rough}
N_D\ll \sigma WM^{1/2}D^{-5/12}.
\end{equation}

The unweighted block estimate in
\cite[Section~6, calculation preceding~(6.8)]{HB2003} is
\begin{equation}\label{ta:HB-block}
N_D\le N_1\ll \sigma K(\mathbf X)\quad\text{and}\quad
K(\mathbf X)
=\mathcal X^{1/4}(X_1/X_2)^{1/8}
 +(X_1/X_2)^{1/4}X_3.
\end{equation}
Since \(X_3\le6X_4\) and \(X_4\le X_2\le W\), we have
\[
\begin{aligned}
\mathcal X^{1/4}(W/X_2)^{1/8}
&=W^{3/8}X_2^{1/8}(X_3X_4)^{1/4}
 \ll WM^{-1/2},\\
(W/X_2)^{1/4}X_3
&\ll W^{1/4}X_4^{3/4}
 =WM^{-3/4}.
\end{aligned}
\]
Interpolating this with \eqref{ta:divisibility-block-rough},
\begin{equation}\label{ta:divisibility-block}
\begin{aligned}
N_D
\ll
\bigl(\sigma WM^{1/2}D^{-5/12}\bigr)^{1/4}
\bigl(\sigma WM^{-1/2}\bigr)^{3/4}
=\sigma WM^{-1/4}D^{-5/48}.
\end{aligned}
\end{equation}

For fixed \(\mathbf X\), every nonempty block satisfies
$\prod_{f\ni i}Z_f\le2X_i.$
There are \(O(\log^6(2+W))\) possible blocks in the variables \(z_e\).
To retain the stronger estimate needed below, put
\[
k_f=\log_2(2Z_f),\qquad
d_i=\log_2(16X_i)-\sum_{f\ni i}k_f.
\]
Then \(k_f,d_i\) are nonnegative integers and
$\frac{P_0^2}{\mathcal X}
=2^{\,4-\sum_i d_i}.$ For fixed \((d_i)\), the choices of \(k_{12},k_{13}\) determine
\((k_{14},k_{23},k_{24},k_{34})\), whose coefficient matrix is
\[
\begin{pmatrix}
1&0&0&0\\
0&1&1&0\\
0&1&0&1\\
1&0&1&1
\end{pmatrix},
\qquad
\det=-2.
\]
Each chosen exponent has \(O(\log(2+W))\) possible values. Therefore
\[
\begin{aligned}
\sum_{\mathbf Z}(\sigma-1)
\ll
\log^3(2+W)
\sum_{d_1,\ldots,d_4\ge0}
2^{-(d_1+\cdots+d_4)/9}
=
\frac{\log^3(2+W)}{(1-2^{-1/9})^4}.
\end{aligned}
\]
Together with the count of blocks in the variables \(z_e\), this gives
\begin{equation}\label{ta:sigma-sum}
\sum_{\mathbf Z}\sigma\ll\log^6(2+W),
\qquad
\sum_{\mathbf Z}(\sigma-1)\ll\log^3(2+W).
\end{equation}

For fixed \(W=X_1\) and \(M=2^m\), we have $X_4=W2^{-m}$ and $X_3/X_4\in\{1,2,4\},$
and there are at most \(m+1\) choices for \(X_2\).
Using \eqref{ta:divisibility-block} and \eqref{ta:sigma-sum},
\[
\begin{aligned}
\sum_{\mathbf X:\,X_1=W}\sum_{\mathbf Z}
N_D(\mathbf X;\mathbf Z)
&\ll
WD^{-5/48}\log^6(2+W)
\sum_{m\ge0}(m+1)2^{-m/4}\\
&\ll
WD^{-5/48}\log^6(2+W).
\end{aligned}
\]
Finally,
\[
\sum_{\substack{W\le X\\W\ \mathrm{dyadic}}}W\ll X,
\qquad
\log(2+W)\ll\log(2X),
\]
so
\[
\begin{aligned}
N_{D,e}(X)
\ll
D^{-5/48}
\sum_{\substack{W\le X\\W\ \mathrm{dyadic}}}
W\log^6(2+W)
\ll
X\log^6(2X)D^{-5/48}.
\end{aligned}
\]
All implied constants are independent of \(D\), \(X\), and \(e\).
\end{proof}

\subsection{An absolutely convergent global sieve tail}

For a positive integer \(h\) put
\[
 N_h^{P}(X)=\#\{x\in U(\mathbb Q):H(x)\le X, h\mid P(z(x))\},
\]
where \(P\) is defined in \eqref{in:torsor-monomials}.

\begin{corollary}\label{ta:product-divisibility}
For \(h\in\mathbb Z_{\ge1}\) and \(X\ge1\),
\[
N_h^P(X)\ll d_6(h)\,X\log^6(2X)\,h^{-5/288},
\]
where \(d_6\) is the sixfold divisor function.
For \(X<1\), one has \(N_h^P(X)=0\).
\end{corollary}

\begin{proof}
The six parameters \(z_e\) in the canonical representation are pairwise coprime.
For each point counted by \(N_h^P(X)\), set
$h_e=\gcd(h,z_e(x)).$
Since \(h\mid\prod_e z_e(x)\), we have $\prod_e h_e=h$ and $h_e\mid z_e(x).$
For every factorization \(h=\prod_e h_e\), we have $\max_e h_e\ge h^{1/6}.$
Grouping the points by these factors, enlarging the sum to all positive
integer factorizations, and applying Proposition~\ref{ta:single-edge},
we obtain
\[
\begin{aligned}
N_h^P(X)
&\le
\sum_{\substack{h_e\ge1\\ \prod_e h_e=h}}
\min_e N_{h_e,e}(X)
\ll
X\log^6(2X)
\sum_{\substack{h_e\ge1\\ \prod_e h_e=h}}
\left(\max_e h_e\right)^{-5/48}\\
&\le
X\log^6(2X)\,h^{-5/288}
\sum_{\substack{h_e\ge1\\ \prod_e h_e=h}}1
=
d_6(h)\,X\log^6(2X)\,h^{-5/288}.
\end{aligned}
\]
Finally, every rational point has \(H(x)\ge1\), so the count vanishes
when \(X<1\).
\end{proof}

\begin{theorem}\label{ta:weak-tail}
Let \(\mathcal T_{>R}(B)\) be as in
Definition~\ref{ta:weak-definition}. Uniformly for \(B,R\ge2\),
\begin{equation}\label{ta:weak-tail-bound}
\mathcal T_{>R}(B)
\ll B\log^6(2B)R^{-5/288}.
\end{equation}
Consequently,
\[
\lim_{R\to\infty}\limsup_{B\to\infty}
\frac{\mathcal T_{>R}(B)}{B(\log B)^6}=0.
\]
The same bound holds after imposing additional restrictions
on the weak representations.
\end{theorem}

\begin{proof}
By Lemma~\ref{ta:normalization} and the fixed sign convention,
each weak representation is uniquely determined by the
primitive vector \(x^*\) and the integers \(g,h\)
obtained in that lemma, with
\[
h\mid P(z(x^*)),\qquad
H_{\mathrm{raw}}=g^3hH(x^*),
\]
and
\[
p\mid y_i,\quad p\mid y_j,\quad i\ne j
\quad\Longrightarrow\quad p\mid gh.
\]
Since \(H(x^*)\ge1\), Corollary~\ref{ta:product-divisibility}
therefore gives
\begin{align}\label{ta:tail-positive-series}
\mathcal T_{>R}(B)
\le
\sum_{\substack{p>R\\p\ \mathrm{prime}}}
\sum_{\substack{g,h\ge1\\g^3h\le B\\p\mid gh}}
N_h^P\!\left(\frac{B}{g^3h}\right)
\ll
B\log^6(2B)
\sum_{\substack{p>R\\p\ \mathrm{prime}}}
\sum_{\substack{g,h\ge1\\p\mid gh}}
\frac{d_6(h)}{g^3h^{293/288}}.
\end{align}
Here \(\log(2B/(g^3h))\le\log(2B)\) on the original range;
the extension and rearrangement of sums are valid by
nonnegativity.

The Euler product for \(d_6\) and
\(1-(1-u)^6\le6u\) for \(0\le u\le1\) give, for every prime \(p\),
\[
\sum_{\substack{g,h\ge1\\p\mid gh}}
\frac{d_6(h)}{g^3h^{293/288}}
\le\zeta(3)\zeta(293/288)^6
\bigl(p^{-3}+6p^{-293/288}\bigr).
\]
Hence, for \(R\ge2\), comparison with the convergent integrals gives
\[
\sum_{\substack{p>R\\p\ \mathrm{prime}}}
\sum_{\substack{g,h\ge1\\p\mid gh}}
\frac{d_6(h)}{g^3h^{293/288}}
\ll\sum_{n>R}\bigl(n^{-3}+n^{-293/288}\bigr)
\ll R^{-5/288}.
\]
Substitution into \eqref{ta:tail-positive-series} proves
\eqref{ta:weak-tail-bound}, uniformly in \(B,R\).
The iterated limit follows, and additional restrictions can only decrease the count.
\end{proof}

\begin{corollary}\label{ta:weak-count}
Let $\mathcal N_{\mathrm{weak}}(B)$ be as in
Definition~\ref{ta:weak-definition}. Uniformly for \(B\ge2\),
\begin{equation}\label{ta:weak-count-bound}
 \mathcal N_{\mathrm{weak}}(B)\ll B\log^6(2B).
\end{equation}
\end{corollary}
\begin{proof}
By Lemma~\ref{ta:normalization} and
Corollary~\ref{ta:product-divisibility}, we have
\[
\begin{aligned}
\mathcal N_{\mathrm{weak}}(B)
&\le
\sum_{\substack{g,h\ge1\\g^3h\le B}}
N_h^P\!\left(\frac{B}{g^3h}\right)
\ll
B\log^6(2B)
\sum_{g,h\ge1}
\frac{d_6(h)}{g^3h^{293/288}}\\
&=
\zeta(3)\zeta(293/288)^6
B\log^6(2B)
\ll B\log^6(2B).
\end{aligned}
\]
\end{proof}

\subsection{Logarithmic boundaries and the coordinate-height tail}

In this subsection parameters are canonical. For a fixed dyadic
height block use \eqref{ta:dyadic-notation}, \eqref{ta:height-scales},
and the integer coordinates \(k_e=\log_2(2Z_e)\) introduced above.
All scale sums below are over nonempty blocks.

\begin{proposition}\label{ta:boundary}
Let \(B\ge2\), \(\delta\ge0\), and \(L=1+\log(2B)\).
Fix \(\boldsymbol\alpha=(\alpha_e)_e\in\mathbb R^6\setminus\{0\}\)
and \(C>0\). For each dyadic height block \(\mathbf X\), let
\(I_{\mathbf X}\subset\mathbb R\) be an interval satisfying
$|I_{\mathbf X}|\le C(\delta L+1).$
The number of rational points \(x\in U(\mathbb Q)\) with
\(H(x)\le B\) whose parameters \(z_e\) in the canonical representation lie in blocks
satisfying
$\sum_e\alpha_e k_e\in I_{\mathbf X}$
is
\begin{equation}\label{ta:boundary-bound}
\ll_{\boldsymbol\alpha,C}\delta BL^6+BL^5.
\end{equation}
The same estimate holds for any fixed finite union of such
conditions, with the implied constant depending on that family.
\end{proposition}

\begin{proof}
Use the ordered dyadic blocks underlying \eqref{ta:HB-block},
with height bound \(B\). The finitely many index permutations
are treated separately, permuting \(\boldsymbol\alpha\) accordingly.
All sums below are over nonempty blocks.

Each \(z_f\) divides a primitive coordinate \(x_i\) with \(i\in f\), so
\[
z_f\le B,\qquad
0\le k_f=\log_2(2Z_f)\le\log_2(2B)\ll L.
\]
Choose an index pair \(e\) with \(\alpha_e\ne0\).
For fixed \((k_f)_{f\ne e}\),
\[
\#\left\{
k_e\in\mathbb Z_{\ge0}:
\alpha_e k_e\in
I_{\mathbf X}-\sum_{f\ne e}\alpha_f k_f
\right\}
\le
1+\frac{|I_{\mathbf X}|}{|\alpha_e|}
\ll_{\boldsymbol\alpha,C}\delta L+1.
\]
There are \(O(L^5)\) choices for the other five coordinates.
Consequently, uniformly in \(\mathbf X\) and the position of
\(I_{\mathbf X}\),
\[
\sum_{\substack{\mathbf Z\\
\sum_e\alpha_e k_e\in I_{\mathbf X}}}1
\ll_{\boldsymbol\alpha,C}
L^5(\delta L+1).
\]
Since \(\sigma-1\ge0\), \eqref{ta:sigma-sum} gives
\[
\begin{aligned}
\sum_{\substack{\mathbf Z\\
\sum_e\alpha_e k_e\in I_{\mathbf X}}}\sigma
\le
\sum_{\substack{\mathbf Z\\
\sum_e\alpha_e k_e\in I_{\mathbf X}}}1
+\sum_{\mathbf Z}(\sigma-1)
\ll_{\boldsymbol\alpha,C}
\delta L^6+L^5+L^3
\ll_{\boldsymbol\alpha,C}
\delta L^6+L^5.
\end{aligned}
\]

For these ordered height scales,
\cite[Section~6, final summation following~(6.8)]{HB2003}
provides $\sum_{\mathbf X}K(\mathbf X)\ll B.$
Applying \eqref{ta:HB-block} with \(D=1\), we therefore obtain
\[
\begin{aligned}
\sum_{\mathbf X}
\sum_{\substack{\mathbf Z\\
\sum_e\alpha_e k_e\in I_{\mathbf X}}}
N_1(\mathbf X;\mathbf Z)
&\ll
\sum_{\mathbf X}K(\mathbf X)
\sum_{\substack{\mathbf Z\\
\sum_e\alpha_e k_e\in I_{\mathbf X}}}\sigma
\ll_{\boldsymbol\alpha,C}
(\delta L^6+L^5)\sum_{\mathbf X}K(\mathbf X)\\
&\ll_{\boldsymbol\alpha,C}
\delta BL^6+BL^5.
\end{aligned}
\]
Summing over the finitely many index orderings proves
\eqref{ta:boundary-bound}. The assertion for a fixed finite union
follows by addition.
\end{proof}

For each fixed height block, \eqref{ta:height-scales} gives
\[
 \log Y_i=\frac13\log\mathcal X-\log X_i
 +\frac13\sum_{e\ni i}\log Z_e
 -\frac23\sum_{e\not\ni i}\log Z_e.
\]
In particular, the coefficient vector is nonzero.
Nonempty blocks satisfy
\[
Z_e\ge\tfrac12,\qquad Y_i\gg1,\qquad |y_i|\asymp Y_i.
\]
Consequently, each of the conditions
\[
z_e<B^\delta,\qquad
|y_i|<B^\delta,\qquad
|\log P_0-\log B|<\delta\log B
\]
implies, on each fixed height block, that a fixed nonzero
linear combination of the logarithmic parameters
\(k_e=\log_2(2Z_e)\) lies in an interval of length
\(O(\delta L+1)\). The coefficients depend only on the
chosen condition, while the interval may depend on the
height block. Proposition~\ref{ta:boundary} therefore applies.

On \(B/M\le |x_i|\le B\), with \(M\ge1\) fixed, retain
\(T=(BP_0)^{1/3}\) and put
$
Q_{ij|kl}=\frac{T}{Z_{ij}Z_{kl}}.$
By \eqref{in:q-identity} and \(|A_i(z)y_i|\asymp_M T\),
\[
 1\le |v_{ij}|
 =\frac{|A_i(z)y_i+A_j(z)y_j|}{z_{ij}z_{kl}}
 \ll_M Q_{ij|kl}.
\]
Hence \(Q_{ij|kl}<B^\delta\) implies
\[
\begin{gathered}
 -O_M(1)\le\log Q_{ij|kl}<\delta\log B,\\
 \log Q_{ij|kl}
 =\frac13\log B+\frac13\sum_e\log Z_e-\log Z_{ij}-\log Z_{kl}.
\end{gathered}
\]
Proposition~\ref{ta:boundary} therefore bounds each of these exceptional sets,
and their finite union, by \(O_M(\delta BL^6+BL^5)\).

\begin{proposition}\label{ta:height-tail}
Let \((x_1,\ldots,x_4)\) be a primitive integral representative
of \(x\in U(\mathbb Q)\). Uniformly for \(B\ge2\) and \(M\ge1\),
\begin{equation}\label{ta:height-tail-bound}
\#\left\{
x\in U(\mathbb Q):
H(x)\le B,\quad \min_i|x_i|\le B/M
\right\}
\ll BM^{-1/2}\log^6(2B).
\end{equation}
\end{proposition}

\begin{proof}
Since \(\min_i|x_i|\ge1\), the count vanishes when \(M>B\).
Assume \(1\le M\le B\).
Use the dyadic decomposition underlying \eqref{ta:HB-block}.
By \cite[Section~6, (6.3)--(6.4)]{HB2003}, it suffices to
consider the finitely many index orderings with
\[
\tfrac12\le X_4\le X_3\le X_2\le X_1\le B,
\qquad
X_3\le6X_4,\qquad X_4\le B/M.
\]
All sums below are over dyadic scales at least \(1/2\).
For each \(X_4\), there are at most three choices of \(X_3\).
The definition in \eqref{ta:HB-block} therefore gives
\[
K(\mathbf X)
\ll
X_1^{3/8}X_2^{1/8}X_4^{1/2}
+
X_1^{1/4}X_2^{-1/4}X_4.
\]

Enlarging the summation ranges and summing geometric series,
we obtain
\[
\begin{aligned}
\sum_{X_4\le B/M}X_4^{1/2}
 \sum_{X_1\le B}X_1^{3/8}
 \sum_{X_2\le X_1}X_2^{1/8}
&\ll
 \sum_{X_4\le B/M}X_4^{1/2}
 \sum_{X_1\le B}X_1^{1/2}\\
&\ll
 B^{1/2}(B/M)^{1/2}
 =BM^{-1/2},
\end{aligned}
\]
and
\[
\begin{aligned}
\sum_{X_1\le B}X_1^{1/4}
 \sum_{X_4\le B/M}X_4
 \sum_{X_2\ge X_4}X_2^{-1/4}
&\ll
 \sum_{X_1\le B}X_1^{1/4}
 \sum_{X_4\le B/M}X_4^{3/4}\\
&\ll
 B^{1/4}(B/M)^{3/4}
 =BM^{-3/4}.
\end{aligned}
\]
Hence, over the ordered height blocks,
\[
\sum_{\substack{\mathbf X\\X_4\le B/M}}K(\mathbf X)
\ll B(M^{-1/2}+M^{-3/4})
\ll BM^{-1/2}.
\]

Take height bound \(B\) in \(N_1(\mathbf X;\mathbf Z)\).
By \eqref{ta:HB-block} and \eqref{ta:sigma-sum}, uniformly
in \(\mathbf X\),
\[
\sum_{\mathbf Z}N_1(\mathbf X;\mathbf Z)
\ll
K(\mathbf X)\sum_{\mathbf Z}\sigma
\ll K(\mathbf X)\log^6(2B).
\]
The required points are contained in the blocks with
\(X_4\le B/M\), so
\[
\begin{aligned}
\sum_{\substack{\mathbf X\\X_4\le B/M}}
\sum_{\mathbf Z}N_1(\mathbf X;\mathbf Z)
\ll
\log^6(2B)
\sum_{\substack{\mathbf X\\X_4\le B/M}}K(\mathbf X)
\ll BM^{-1/2}\log^6(2B).
\end{aligned}
\]
Summing over the finitely many index orderings proves
\eqref{ta:height-tail-bound}.
\end{proof}

\subsection{Proximity to a line for the pairing with the largest product}

We require a primitive lattice estimate which remains useful in very
thin regions.

\begin{lemma}\label{ta:ellipse}
Let \(\Lambda\subset\mathbb R^2\) be a rank-two lattice,
with covolume \(\det\Lambda\).
A nonzero vector \(v\in\Lambda\) is primitive in \(\Lambda\) if
$v\notin m\Lambda\;(m=2,3,\ldots).$
For an ellipse \(E\) centered at the origin, including its interior,
\[
\#\{v\in\Lambda\cap E:
       v\text{ is primitive in }\Lambda\}
\ll
1+\frac{\operatorname{area}(E)}{\det\Lambda}.
\]
The same estimate holds for a parallelogram centered at the origin.
The implied constant is absolute.
\end{lemma}

\begin{proof}
By \cite[Lemma~4]{Reuss2015}, we have
\[
\#\{v\in\Lambda\cap E:
      v\text{ is primitive in }\Lambda\}
\le
4\left(1+\frac{\operatorname{area}(E)}{\det\Lambda}\right).
\]
For a centered parallelogram, write
\(P=A([-1,1]^2)\), where \(A\in\mathrm{GL}_2(\mathbb R)\).
Then
\[
P\subseteq
A\{u\in\mathbb R^2:u_1^2+u_2^2\le2\}\]
and
\[\operatorname{area}
\bigl(A\{u\in\mathbb R^2:u_1^2+u_2^2\le2\}\bigr)
=2\pi|\det A|
=\frac{\pi}{2}\operatorname{area}(P).
\]
Applying the same estimate to this enclosing ellipse gives
\[
\#\{v\in\Lambda\cap P:
      v\text{ is primitive in }\Lambda\}
\le
4+2\pi\frac{\operatorname{area}(P)}{\det\Lambda}
\le
2\pi\left(1+\frac{\operatorname{area}(P)}{\det\Lambda}\right).
\]
\end{proof}

\begin{theorem}\label{ta:near-line}
For \eqref{in:scales}--\eqref{in:short-scales}, fix \(M\ge1\), \(\delta>0\), and a dyadic block in the variables \(z_e\) with
\(A\ge H\ge J\). Suppose that
\begin{equation}\label{ta:near-core}
\begin{gathered}
 \min_e Z_e\ge B^\delta,\qquad
 \min_iY_i\ge B^\delta,\qquad Q\ge B^\delta,\\
 P_0\le B^{1-\delta}\quad\text{or}\quad P_0\ge B^{1+\delta}.
\end{gathered}
\end{equation}
For every fixed \(0<\kappa<\delta/32\),
\begin{equation}\label{ta:near-line-bound}
\#\left\{
\begin{array}{l}
(z,y)\text{ canonical in the block}:\\
B/M\le |x_i|\le B\quad(1\le i\le4),\\
0<|A_1(z)y_1+A_2(z)y_2|\le B^{-\kappa}T
\end{array}
\right\}
\ll_{\delta,M,\kappa}B^{1-\kappa/12}.
\end{equation}
Summing over all such blocks in the variables \(z_e\) gives
\(o(B(\log B)^6)\).
The same conclusions hold for any fixed positive multiple of the
displayed cutoff, with the implied constant depending additionally
on that multiple.
\end{theorem}

\begin{proof}
For \(0<\rho\) sufficiently small in terms of \(M\), define
\[
N_{\mathrm{near}}(\rho)
:=
\#\left\{
(z,y):
\begin{array}{l}
(z,y)\text{ is a canonical representation},\quad y_1>0,\\
Z_e\le z_e<2Z_e\quad\text{for every index pair }e,\\
B/M\le |x_i|\le B\quad(1\le i\le4),\\
0<|A_1(z)y_1+A_2(z)y_2|\le \rho T
\end{array}
\right\},
\]
where \(x_i=C_i(z)\prod_{j\ne i}y_j\) and
\(T=(BP_0)^{1/3}\).
Retain the notation for the six variables
\[
a=z_{12},\quad b=z_{34},\quad c=z_{14},\quad d=z_{23},
\qquad u=z_{13},\quad v=z_{24},
\]
and write \(U=Z_{13}\), \(V=Z_{24}\).
With \(F_0\) as in \eqref{in:short-scales}, put
\begin{equation*}
D_0=JF_0=\frac{T^2}{AH},
\qquad
K=\sqrt{AH}.
\end{equation*}
By \eqref{in:scales}--\eqref{in:short-scales} and the
definitions of \(D_0\) and \(K\), we have
\begin{equation*}
\begin{gathered}
U^2Y_2Y_4=V^2Y_1Y_3=D_0,\quad
\prod_iY_i=F_0^2,\quad
K\sqrt{D_0}=T,\quad
B=TF_0=\frac{P_0}{T}F_0^2.
\end{gathered}
\end{equation*}

Since \(|A_i(z)y_i|\asymp_M T\), near cancellation forces opposite
signs in both complementary pairs. Thus
\[
(y_1,y_2,y_3,y_4)
=(t_1,-t_2,\epsilon t_3,-\epsilon t_4),
\qquad
t_i>0,\quad \epsilon\in\{-1,1\},
\]
where \(t_i\asymp_M Y_i\).
Set
$
L_1=-ut_2c+vt_1d$ and
$L_2=\epsilon(-vt_3c+ut_4d).$
The equation \(\sum_i A_i(z)y_i=0\) in \eqref{in:torsor}
and the bound
\(0<|A_1(z)y_1+A_2(z)y_2|\le\rho T\)
in the definition of \(N_{\mathrm{near}}(\rho)\) give
\[
bL_1=aL_2,\qquad
0<|bL_1|\le\rho T,\qquad
|L_1|\le\frac{\rho T}{b},\quad
|L_2|\le\frac{\rho T}{a}.
\]
Define
$\Delta=\epsilon(v^2t_1t_3-u^2t_2t_4).$
Canonical coprimality gives
$\gcd(u^2t_2t_4,v^2t_1t_3)=1.$
Hence \(\Delta=0\) would force
\(u=v=t_1=t_2=t_3=t_4=1\), contrary to
\eqref{ta:near-core}.

The adjugate identity
$\Delta c=\epsilon ut_4L_1-vt_1L_2$
and
$UY_4Z_{12}=VY_1Z_{34}=\frac{T}{Z_{23}}$
give
\[
\begin{aligned}
|\Delta|Z_{14}
\ll_M
\rho T\left(
\frac{UY_4}{Z_{34}}+\frac{VY_1}{Z_{12}}
\right)
\ll_M\frac{\rho T^2}{AZ_{23}}.
\end{aligned}
\]
Therefore
\begin{equation}\label{ta:small-det}
0<|\Delta|\ll_M\frac{\rho T^2}{AH}
=\rho D_0.
\end{equation}

For fixed \((u,v,t_1,t_2,t_3,t_4)\), the equation
\(bL_1-aL_2=0\) is the bilinear equation
\[
(a,b)
\begin{pmatrix}
\epsilon vt_3&-\epsilon ut_4\\
-ut_2&vt_1
\end{pmatrix}
\binom{c}{d}=0.
\]
Its coefficient matrix has determinant \(\Delta\ne0\) and
coprime entries; moreover,
$\gcd(a,b)=\gcd(c,d)=1.$
Thus \cite[Theorem~2]{Reuss2015}, followed by the standard divisor
bound, yields
\[
\#\{(a,b,c,d)\text{ in their prescribed blocks}\}
\ll
\tau(|\Delta|)
\left(1+\frac{K}{\sqrt{|\Delta|}}\right)
\ll_\varepsilon
B^{\varepsilon/2}
\left(1+\frac{K}{\sqrt{|\Delta|}}\right).
\]
Here \(z_e,t_i\le B\), so all divisor arguments have polynomial
size in \(B\).

All sums below run over the six short variables in their prescribed
ranges and retain their inherited coprimality conditions.
For fixed \(\Delta=d_0\ne0\) and fixed coprime \(u,v\), the product
pair \((t_2t_4,t_1t_3)\) lies in an arithmetic progression with
step \((v^2,u^2)\). Its ranges satisfy
\[
t_2t_4\asymp_M\frac{D_0}{U^2},
\qquad
t_1t_3\asymp_M\frac{D_0}{V^2}.
\]
Counting the product pairs and their factorizations gives
\begin{equation*}
\begin{aligned}
\#\{\text{short tuples}:\Delta=d_0\}
\ll_{M,\varepsilon}
B^{\varepsilon/2}UV
\left(1+\frac{D_0}{U^2V^2}\right)
=
B^{\varepsilon/2}(J+F_0).
\end{aligned}
\end{equation*}
Consequently,
\[
\begin{aligned}
K\sum_{\substack{u,v,t_1,\ldots,t_4\\
                  0<|\Delta|\ll_M\rho D_0}}
 |\Delta|^{-1/2}
\ll_{M,\varepsilon}
B^{\varepsilon/2}K(J+F_0)\sqrt{\rho D_0}
=
B^{1+\varepsilon/2}\sqrt{\rho}
\left(1+\frac{J}{F_0}\right).
\end{aligned}
\]

For a second estimate, fix \(t_1,t_2,t_3,t_4\).
If \(0<|\Delta|\le E\), for $|v|\ll V$, we have
\[
\begin{aligned}
\left|u-\sqrt{\frac{t_1t_3}{t_2t_4}}\,v\right|
=
\frac{|\Delta|}
{t_2t_4\left(u+\sqrt{t_1t_3/(t_2t_4)}\,v\right)}
\ll_M\frac{EU}{D_0}.
\end{aligned}
\]
These inequalities define a centered parallelogram of area
\(O_M(EUV/D_0)=O_M(E/F_0)\).
Since \(\gcd(u,v)=1\), Lemma~\ref{ta:ellipse} gives
\(O_M(1+E/F_0)\) possibilities for \((u,v)\).
There are \(O_M(\prod_iY_i)=O_M(F_0^2)\) choices for the four
\(t_i\), and hence
$\#\{\text{short tuples}:0<|\Delta|\le E\}
\ll_M F_0^2+EF_0.$
Dyadic summation now gives
\[
\begin{aligned}
K\sum_{\substack{u,v,t_1,\ldots,t_4\\
                  0<|\Delta|\ll_M\rho D_0}}
 |\Delta|^{-1/2}
\ll_M KF_0^2+KF_0\sqrt{\rho D_0}
=
B\left(\sqrt{\frac{F_0}{J}}+\sqrt{\rho}\right).
\end{aligned}
\]
If the determinant range is nonempty, combine the two estimates
and split at \(J/F_0=\rho^{-1/3}\):
\[
\min\left\{
\sqrt{\rho}\left(1+\frac{J}{F_0}\right),
\sqrt{\frac{F_0}{J}}+\sqrt{\rho}
\right\}
\le2\rho^{1/6}.
\]
Thus
\[
K\sum_{\substack{u,v,t_1,\ldots,t_4\\
                  0<|\Delta|\ll_M\rho D_0}}
 |\Delta|^{-1/2}
\ll_{M,\varepsilon}B^{1+\varepsilon/2}\rho^{1/6}.
\]

To bound the constant term in the bilinear estimate, fix five
short variables. As a function of the remaining \(t_i\),
\(\Delta\) has coefficient of size \(\asymp_M D_0/Y_i\).
Therefore
\[
\begin{aligned}
\#\{\text{short tuples}:0<|\Delta|\ll_M\rho D_0\}
\ll_M
\frac{JF_0^2}{Y_i}(1+\rho Y_i)
\ll_M JF_0^2(\rho+B^{-\delta}).
\end{aligned}
\]
If \(J\le P_0/T\), then \(JF_0^2\le B\), so the bilinear
estimate proves
\begin{equation}\label{ta:near-small-R3}
N_{\mathrm{near}}(\rho)
\ll_{M,\varepsilon}
B^{1+\varepsilon}
\left(\rho^{1/6}+B^{-\delta}\right).
\end{equation}

It remains to treat \(J>P_0/T\).
Since \(A\ge H\ge J>P_0/T\), we have
\[
P_0=AHJ>
\left(\frac{P_0}{T}\right)^3
=\frac{P_0^2}{B},
\]
and hence \eqref{ta:near-core} implies
\[
P_0\le B^{1-\delta},
\qquad
A\le
\frac{P_0}{(P_0/T)^2}
=\frac{T^2}{P_0}=F_0.
\]
For complementary index pairs \(ij|kl\), we have
\[
\frac{P_0Z_{ij}}{T^2}
=\frac{Z_{ij}}{F_0}
=\frac{Z_{ij}Z_{kl}}{F_0Z_{kl}}
\le\frac1{Z_{kl}}\le B^{-\delta}.
\]
Together with \(A_i(\mathbf Z)/T=Y_i^{-1}\) and
\(A/T=Q^{-1}\), with \(H,J\le A\), this gives
\begin{equation}\label{ta:near-poisson-core}
\frac{P_0}{B}\le B^{-\delta},\quad
\frac{P_0Z_e}{T^2}\le B^{-\delta},\quad
\frac{A_i(\mathbf Z)}{T}\le B^{-\delta}\quad\text{and}\quad
\frac{A}{T},\ \frac{H}{T},\ \frac{J}{T}\le B^{-\delta}.
\end{equation}

For fixed \(z\), put \(n_i=A_i(z)y_i\), so that
$n_4=-n_1-n_2-n_3,\; |n_i|\asymp_M T.$
In the normalized coordinates \(\xi=(n_1,n_2,n_3)/T\), choose
fixed nonnegative functions \(\psi\in C_c^\infty(\mathbb R^3)\) and
\(\chi\in C_c^\infty(\mathbb R)\), with \(\psi\ge1\) on a fixed compact set containing all
\(\xi\) corresponding to points satisfying
\(B/M\le|x_i|\le B\) for \(1\le i\le4\),
and \(\chi\ge1\) on \([-1,1]\).
Put
$
f_\rho(\xi)=\psi(\xi)\chi\!\left(\frac{\xi_1+\xi_2}{\rho}\right).$
Thus \(f_\rho(\xi)\ge1\) for every point counted by
\(N_{\mathrm{near}}(\rho)\). Moreover,
\[
\int_{\mathbb R^3}f_\rho\ll_M\rho,
\qquad
\|f_\rho\|_{C^{16}}\ll_M\rho^{-16}.
\]
Drop the coprimality conditions on \(y\), retaining pairwise
coprimality of the six variables \(z_e\). By \eqref{po:det}, the lattice in
\((n_1,n_2,n_3)\) has covolume \(P(z)^2\). Its zero-frequency
contribution, summed over the block, is
\[
\sum_z\frac{T^3}{P(z)^2}\int_{\mathbb R^3}f_\rho
\ll_M
\rho P_0\frac{BP_0}{P_0^2}
=\rho B.
\]

Apply \eqref{po:frequency} with \(W_R=1\) and normalized weight
$\mathbf t\longmapsto
 f_\rho\!\left((P(z)/P_0)^{1/3}\mathbf t\right).$
Since \(1\le P(z)/P_0<2^6\), its \(C^{16}\) norm is
\(\ll_M\rho^{-16}\), uniformly in \(z\).
No derivatives in the variables \(z_e\) are required.
The frequency estimate is independent of the present theorem and gives
\[
\begin{aligned}
N_{\mathrm{near}}(\rho)
\ll_{M,\varepsilon}\;&
\rho B+\rho^{-16}B^{1+\varepsilon}
\left\{
\frac{P_0}{B}
+\frac{P_0}{T^2}\sum_eZ_e
+\frac{\sum_iA_i(\mathbf Z)+A+H+J}{T}
\right\}.
\end{aligned}
\]
Using \eqref{ta:near-poisson-core},
\begin{equation}\label{ta:near-poisson-bound}
N_{\mathrm{near}}(\rho)
\ll_{M,\varepsilon}
\rho B+\rho^{-16}B^{1-\delta+\varepsilon}.
\end{equation}

Finally, take
$
\rho=B^{-\kappa},\;
0<\kappa<\frac{\delta}{32}$ and
$0<\varepsilon<
\min\left(\frac{\kappa}{12},\frac{\delta}{4}\right).$
The exponents in \eqref{ta:near-small-R3} satisfy
\[
\varepsilon-\frac{\kappa}{6}<-\frac{\kappa}{12},
\qquad
\varepsilon-\delta<-\frac{\kappa}{12},
\]
whereas \eqref{ta:near-poisson-bound} gives
\[
-\kappa\le-\frac{\kappa}{12},
\qquad
-\delta+\varepsilon+16\kappa
<-\frac{\delta}{4}<-\frac{\kappa}{12}.
\]
This proves \eqref{ta:near-line-bound}; the bounded range of \(B\)
is absorbed into the implied constant.

There are \(O(\log^6(2B))\) possible blocks in the variables \(z_e\), so
\[
\sum_{\text{admissible blocks}}N_{\mathrm{near}}(B^{-\kappa})
\ll_{\delta,M,\kappa}
B^{1-\kappa/12}\log^6(2B)
=o(B(\log B)^6).
\]
The same estimates hold if the upper bound \(B^{-\kappa}T\)
is multiplied by any fixed positive constant.
The implied constants may depend on this multiplier.
\end{proof}

\begin{remark}\label{ta:near-line-scope}
Theorem~\ref{ta:near-line} applies to canonical representations
and the pairing \(ij|kl\) for which \(Z_{ij}Z_{kl}\) is largest.
The factor \(B^\varepsilon\) in \eqref{ta:near-small-R3}
is absorbed by the cutoff \(\rho=B^{-\kappa}\), with
\(\varepsilon>0\) sufficiently small.
For fixed \(\rho>0\), however, that estimate alone does not
bound \(N_{\mathrm{near}}(\rho)/B\) independently of \(B\).
\end{remark}


\section{Averaging reciprocal shortest-vector lengths}
\label{ge:section}

We estimate the sum of reciprocal shortest-vector lengths arising
from the two-dimensional lattice count.

\subsection{Size parameters and the lattice}
\label{ge:scales-subsection}

Let \(B\ge2\), and let the six size parameters
$(Z_{12},Z_{34},Z_{14},Z_{23},Z_{13},Z_{24})
=(Z_a,Z_b,C,D,U,V)$
be at least \(1/2\). Set
$A=Z_aZ_b,\; H=CD,\;J=UV$ and $P_0=AHJ.$
Using the monomials in \eqref{in:torsor-monomials}, define
\begin{equation}\label{ge:scales}
T=(BP_0)^{1/3},\qquad
Q=\frac{T}{A},\qquad
Y_i=\frac{T}{A_i(\mathbf Z)}
\quad(1\le i\le4).
\end{equation}
Assume \(Y_i\ge1/2\), and put \(P_y=\prod_iY_i\).
Since \(\prod_iA_i(\mathbf Z)=P_0^2\),
\begin{equation}\label{ge:identities}
P_y=\frac{T^4}{P_0^2},\qquad
\frac{HJP_y}{Q}
=P_y^{3/4}\sqrt{P_0}
=\frac{T^3}{P_0}
=B.
\end{equation}

Let \(u,v\) be positive integers and \(y_i\) nonzero integers with
$u\asymp U,\; v\asymp V$ and $|y_i|\asymp Y_i,$
satisfying
\begin{equation}\label{ge:short-gcd}
\gcd(u,v)=1,\qquad
\gcd(y_i,y_j)=1\quad(i\ne j),\qquad
\gcd(u,y_1y_3)=\gcd(v,y_2y_4)=1.
\end{equation}
For \((c,d)\in\mathbb Z^2\), define
\begin{equation}\label{ge:linear-forms}
\begin{gathered}
\binom{L_1(c,d)}{L_2(c,d)}
=
\begin{pmatrix}
uy_2&vy_1\\
vy_3&uy_4
\end{pmatrix}
\binom{c}{d},\quad
\Delta=u^2y_2y_4-v^2y_1y_3.
\end{gathered}
\end{equation}
We sum over positive integers \(q\asymp Q_0=sQ\), where
\(0<s\ll1\), subject to
\begin{equation}\label{ge:q-conditions}
\Delta\ne0,\qquad q\mid\Delta,\qquad
\gcd(q,uvy_1y_2y_3y_4)=1.
\end{equation}
The last condition follows from \(q\mid\Delta\) and
\eqref{ge:short-gcd}, since
\[
\gcd(u^2y_2y_4,v^2y_1y_3)=1
\quad\Longrightarrow\quad
\gcd(\Delta,uvy_1y_2y_3y_4)=1.
\]
Replacing a signed auxiliary integer by its absolute value
changes the upper bounds by at most a factor of two.

For each such tuple, let
$\Lambda_q=
\{(c,d)\in\mathbb Z^2:
q\mid L_1(c,d),\ q\mid L_2(c,d)\}.$
The coefficient matrix in \eqref{ge:linear-forms} has determinant
\(\Delta\), and
\[
q\mid\Delta,\qquad
\gcd(q,uy_2,vy_1,vy_3,uy_4)=1.
\]
Thus \cite[Lemma~3]{Reuss2015} applies for \(q\ge2\); for \(q=1\),
one has \(\Lambda_1=\mathbb Z^2\). Hence
\begin{equation}\label{ge:det-lattice}
\det\Lambda_q=q.
\end{equation}

Define the shortest-vector length in the normalized sup norm by
\begin{equation}\label{ge:G}
\lambda_q=
\min_{(c,d)\in\Lambda_q\setminus\{(0,0)\}}
\max\left\{\frac{|c|}{C},\frac{|d|}{D}\right\},
\qquad
\mathcal G=\sum_{u,v,\mathbf y,q}\lambda_q^{-1},
\end{equation}
where the sum runs over the prescribed ranges and
\eqref{ge:short-gcd}--\eqref{ge:q-conditions}.
Empty ranges contribute zero. Additional restrictions may be
imposed on the summation tuples.

Since $Y_i,Z_e\ge1/2$, the identities above give
\[
P_y\ge2^{-4},\qquad P_0=\frac{B^2}{P_y^{3/2}}\ll B^2,
\qquad T=\frac{B}{P_y^{1/2}}\ll B,
\]
and hence $Z_e\ll B^2$ and $Y_i,Q,q\ll B$.
In particular,
\[
 u,v\ll B^2,\qquad |y_i|,q\ll B,\qquad |\Delta|\ll B^6.
\]
Thus the divisor and logarithmic factors in these summation variables
may be absorbed into \(B^\varepsilon\), after choosing the auxiliary exponents sufficiently small.

\subsection{An eight-variable unit estimate}

We use the fourth-moment theorem of Cochrane and Shi
\cite[Theorem~1]{CS2010}.  For an integer interval \(I\) of cardinality
\(L\ge1\) and an integer modulus \(m\ge2\), it gives
\begin{equation}
 \frac1{\varphi(m)}\sum_{\chi\ne\chi_0}
       \left|\sum_{n\in I}\chi(n)\right|^4
 \ll_\varepsilon m^\varepsilon L^2.
 \label{ge:CS}
\end{equation}
This estimate includes imprimitive characters and translated intervals;
the length need not be less than $m$.
Writing $S_I(\chi)=\sum_{n\in I}\chi(n)$, Cauchy--Schwarz gives
\begin{equation}
\sum_{\chi\ne\chi_0}
|S_{I_1}(\chi)S_{I_2}(\overline\chi)|^2
\ll_\varepsilon \varphi(m)m^\varepsilon|I_1||I_2|.
\label{ge:pair-moment}
\end{equation}

\begin{lemma}\label{ge:unit-count}
Let \(m\ge1\). Let \(u,v,c,d,y_1,y_2,y_3,y_4\) range over
integer intervals of respective cardinalities
\[
U,V,C_*,D_*,Y_1,Y_2,Y_3,Y_4\ge1,
\]
with every variable coprime to \(m\). Put
$J=UV,$ $W_*=C_*D_*$ and $P_y=\prod_{i=1}^4Y_i.$
The number \(N_m\) of solutions of
\[
ucy_2+vdy_1\equiv0,\qquad
vcy_3+udy_4\equiv0\pmod m
\]
satisfies, for every \(\varepsilon>0\),
\begin{equation}\label{ge:unit-bound}
N_m\ll_\varepsilon m^\varepsilon
\left\{
\frac{P_yJW_*}{m^2}
+\sqrt{P_yJW_*}
+\frac{P_y\sqrt{JW_*}}{m}
+\frac{JW_*\sqrt{P_y}}{m}
\right\}.
\end{equation}
The estimate is uniform in the interval positions and cardinalities.
\end{lemma}

\begin{proof}
For \(m=1\), one has \(N_1=P_yJW_*\).
Assume \(m\ge2\). All characters below are modulo \(m\), extended
by zero on nonunits; \(\chi_0\) denotes the principal character.

For each variable \(w\), let \(I_w\) be its interval. Using the
preceding notation \(S_I(\chi)=\sum_{n\in I}\chi(n)\), set
\[
F_{r,t}(\chi)
=
|S_{I_r}(\chi)S_{I_t}(\chi)|.
\]
Then
\[
F_{r,t}(\overline\chi)=F_{r,t}(\chi),
\qquad
F_{r,t}(\chi)\le |I_r|\,|I_t|.
\]
By \eqref{ge:pair-moment},
\[
\sum_{\chi\ne\chi_0}F_{r,t}(\chi)^2
\ll_\varepsilon\varphi(m)m^{\varepsilon/4}|I_r|\,|I_t|.
\]

Multiplicative character orthogonality, followed by absolute values,
yields
\begin{equation*}
N_m\le
\frac1{\varphi(m)^2}
\sum_{\chi,\psi}
F_{y_1,y_2}(\chi)
F_{y_3,y_4}(\psi)
F_{c,d}(\chi\psi)
F_{u,v}(\chi\overline\psi).
\end{equation*}

First restrict to $\chi,\ \psi,\ \chi\psi,\ \chi\overline\psi\ne\chi_0,$
and denote this restriction by \(\sum^*\).
Every nonempty fibre of $(\chi,\psi)\longmapsto
(\chi\psi,\chi\overline\psi)$
has cardinality $\#\{\xi:\xi^2=\chi_0\}
\le\tau(m)\ll_\varepsilon m^{\varepsilon/4}.$
Here the inequality follows from the structure of the unit groups
modulo prime powers. Cauchy--Schwarz therefore gives
\[
\begin{aligned}
&\frac1{\varphi(m)^2}
\sum_{\chi,\psi}^{*}
F_{y_1,y_2}(\chi)F_{y_3,y_4}(\psi)
F_{c,d}(\chi\psi)F_{u,v}(\chi\overline\psi)\\
&\quad\le
\frac{\#\{\xi:\xi^2=\chi_0\}^{1/2}}{\varphi(m)^2}
\prod_{(r,t)\in\{(y_1,y_2),(y_3,y_4),(c,d),(u,v)\}}
\left(
\sum_{\xi\ne\chi_0}F_{r,t}(\xi)^2
\right)^{1/2}\\
&\quad\ll_\varepsilon
m^\varepsilon\sqrt{P_yJW_*}.
\end{aligned}
\]

For the remaining cases, the trivial bound for $F_{r,t}$ and
\eqref{ge:pair-moment} give the following contributions, up to a
factor $O_\varepsilon(m^\varepsilon)$:
\[
\begin{array}{c|c}
\text{characters} & \text{contribution}\\ \hline
\text{exactly one of }\chi,\psi\text{ principal}
& P_y\sqrt{JW_*}/\varphi(m)\\[1mm]
\chi,\psi\ne\chi_0,\quad
\chi\psi=\chi_0\text{ or }\chi\overline\psi=\chi_0
& JW_*\sqrt{P_y}/\varphi(m)\\[1mm]
\chi=\psi=\chi_0 & P_yJW_*/\varphi(m)^2.
\end{array}
\]
For example, the first row is bounded by
\[
\frac{2P_y}{\varphi(m)^2}
\left(\sum_{\xi\ne\chi_0}F_{c,d}(\xi)^2\right)^{1/2}
\left(\sum_{\xi\ne\chi_0}F_{u,v}(\xi)^2\right)^{1/2}.
\]
For the second row, $\psi=\overline\chi$ or $\psi=\chi$;
apply the same inequality to the two $y$-pairs.
Counting the intersection twice preserves the upper bound.
Together with the all-nonprincipal case and
$m/\varphi(m)\ll_\varepsilon m^\varepsilon$, this proves
\eqref{ge:unit-bound}.
\end{proof}

\subsection{The content of an actual shortest vector}

\begin{lemma}\label{ge:content}
Let $\mathbf w=e(c_0,d_0)\in\Lambda_q\setminus\{0\},
\; e\in\mathbb Z_{>0},\; \gcd(c_0,d_0)=1,$
attain the minimum in \eqref{ge:G}.
Then \(e\mid q\). If \(m=q/e\), we have
\[
m\mid L_i(c_0,d_0)\quad(i=1,2),
\qquad
\gcd(m,c_0d_0)=1.
\]
\end{lemma}

\begin{proof}
Since $q\mathbb Z^2\subseteq\Lambda_q$, B\'ezout's identity gives
\[
\gcd(e,q)(c_0,d_0)\in\Lambda_q\setminus\{0\},\qquad
\lambda_q\le\frac{\gcd(e,q)}e\lambda_q\le\lambda_q.
\]
Thus $e\mid q$, and linearity gives $m\mid L_i(c_0,d_0)$.
Both $uy_2$ and $vy_1$ are units modulo $m$, so
\[
uy_2c_0\equiv-vy_1d_0\pmod m
\Longrightarrow \gcd(m,c_0)=\gcd(m,d_0)\mid\gcd(c_0,d_0)=1.
\]
Hence $\gcd(m,c_0d_0)=1$.
\end{proof}

\begin{proposition}
\label{ge:four-proposition}
Under the notation and assumptions of
\eqref{ge:scales}--\eqref{ge:G}, for every \(\varepsilon>0\),
\begin{equation}
\frac{\mathcal G}{B}
\ll_{\varepsilon}
B^\varepsilon
\left\{
s^{-1/2}\sqrt{\frac QH}
+s\frac{Q^2}{\sqrt{HJP_y}}
+\frac Q{\sqrt{HJ}}
+\frac Q{\sqrt H\,P_y^{1/4}}
+s^{-1}(C^{-1}+D^{-1})
\right\}.
\label{ge:four}
\end{equation}
The implied constant is independent of \(s\) and the individual
scales, with the interval constants fixed.
\end{proposition}

\begin{proof}
Empty summation ranges contribute zero. Otherwise \(Q_0\gg1\).
For each admissible tuple, choose one shortest vector by a fixed rule
and write it as \(e(c_0,d_0)\), as in Lemma~\ref{ge:content}. Put
$\mathcal A(R)
=
\#\left\{
(u,v,\mathbf y,q):
\lambda_q\le R,\quad c_0d_0\ne0
\right\},$
retaining all restrictions in \eqref{ge:G}.

Since every short-variable scale is at least \(1/2\), we have
$\#\{(u,v,\mathbf y)\}\ll UV\prod_iY_i=JP_y.$
On the original tuples, $\Delta\ne0$ and $q\mid\Delta$; hence
\begin{equation}
\mathcal A(R)
\le\sum_{\substack{u,v,\mathbf y\\\Delta\ne0}}\tau(|\Delta|)
\ll_\varepsilon B^{\varepsilon/2}JP_y.
\label{ge:divisor-cap}
\end{equation}
We retain this bound when enlarging the count below.

By Lemma~\ref{ge:content}, one has
\[
e\mid q,\qquad
m=\frac qe\asymp\frac{Q_0}{e},\qquad
m\mid L_1(c_0,d_0),\ L_2(c_0,d_0),\qquad
\gcd(m,c_0d_0)=1.
\]
For a tuple counted by \(\mathcal A(R)\), we have
$1\le |c_0|\le\frac{CR}{e}$ and
$1\le |d_0|\le\frac{DR}{e}.$
Hence a nonempty range satisfies
\[
e\ll Q_0,\qquad
\#\{m:m\asymp Q_0/e\}\ll Q_0/e.
\]

Split the nonzero coordinates by sign. The resulting intervals have
cardinalities \(O(U),O(V),O(Y_i)\) for the short variables and at most
\(CR/e,DR/e\) for \(c_0,d_0\), so \(W_*\ll HR^2/e^2\).
All eight variables are units modulo \(m\), and the chosen-vector rule
gives an injection into \((e,m,u,v,\mathbf y,c_0,d_0)\).
We may therefore discard the remaining restrictions and apply
Lemma~\ref{ge:unit-count} with the actual interval cardinalities to obtain
\[
\begin{aligned}
\mathcal A(R)
\ll_\varepsilon B^{\varepsilon/4}
\sum_{1\le e\ll Q_0}
\biggl\{
\frac{P_yJH}{Q_0}\frac{R^2}{e}
+Q_0\sqrt{P_yJH}\frac{R}{e^2}
+P_y\sqrt{JH}\frac{R}{e}
+JH\sqrt{P_y}\frac{R^2}{e^2}
\biggr\}.
\end{aligned}
\]
Using
\[
\sum_{1\le e\ll Q_0}\frac1e\ll\log(2+Q_0),
\qquad
\sum_{e\ge1}\frac1{e^2}\ll1,
\]
we obtain
\begin{equation}
\mathcal A(R)
\ll_\varepsilon B^{\varepsilon/2}
\left\{
\frac{P_yJH}{Q_0}R^2
+Q_0\sqrt{P_yJH}\,R
+P_y\sqrt{JH}\,R
+JH\sqrt{P_y}\,R^2
\right\}.
\label{ge:A-four}
\end{equation}

Every nonzero integer vector has normalized length at least
$1/\max(C,D)$. Minkowski's theorem and \eqref{ge:det-lattice} give
$\lambda_q\ll\sqrt{q/H}$; fix
$R_{\max}\asymp\sqrt{Q_0/H}$ large enough that every
$\lambda_q\le R_{\max}$. Since
\(\lambda_q^{-1}=\int_{\lambda_q}^{\infty}R^{-2}\,dR\),
summing over the same tuples gives
\begin{equation}
\sum_{\substack{u,v,\mathbf y,q\\c_0d_0\ne0}}\lambda_q^{-1}
=\int_0^\infty\frac{\mathcal A(R)}{R^2}\,dR,\quad\text{with}\quad
\mathcal A(R)=0\;\left(R<\frac1{\max(C,D)}\right).
\label{ge:layer-cake}
\end{equation}

For \(N,K>0\), we have
\[
\int_0^\infty
\frac{\min\{N,KR^2\}}{R^2}\,dR
=2\sqrt{NK},
\]
and
\[
\int_{1/\max(C,D)}^\infty
\frac{\min\{N,KR\}}{R^2}\,dR
\ll
K\log\left(2+\frac{N\max(C,D)}K\right).
\]
Furthermore, for \(t_j\ge0\),
\[
\min\left\{N,\sum_{j=1}^{4}t_j\right\}
\le\sum_{j=1}^{4}\min\{N,t_j\}.
\]
Apply these formulas to
\eqref{ge:divisor-cap} and \eqref{ge:A-four}, with \(N=JP_y\).
For the two linear terms, the logarithms are
\(O(\log(2B))\), by \(Q_0\gg1\) and the scale bounds in
Section~\ref{ge:scales-subsection}. Absorbing them into \(B^{\varepsilon/2}\), we obtain
\[
\begin{aligned}
\sum_{\substack{u,v,\mathbf y,q\\c_0d_0\ne0}}
\lambda_q^{-1}
\ll_\varepsilon B^\varepsilon
\biggl\{
JP_y\sqrt{\frac H{Q_0}}
+Q_0\sqrt{P_yJH}
+P_y\sqrt{JH}
+J\sqrt H\,P_y^{3/4}
\biggr\}.
\end{aligned}
\]
Since \(Q_0=sQ\) and \(B=HJP_y/Q\), one has
\[
\begin{aligned}
\frac{JP_y}{B}\sqrt{\frac H{Q_0}}
&=s^{-1/2}\sqrt{\frac QH},&
\frac{Q_0\sqrt{P_yJH}}B
&=s\frac{Q^2}{\sqrt{HJP_y}},\\
\frac{P_y\sqrt{JH}}B
&=\frac Q{\sqrt{HJ}},&
\frac{J\sqrt H\,P_y^{3/4}}B
&=\frac Q{\sqrt H\,P_y^{1/4}}.
\end{aligned}
\]

Finally, if \(c_0d_0=0\), Lemma~\ref{ge:content} gives
\[
\gcd(m,0)=1
\quad\Longrightarrow\quad
m=1,\quad e=q.
\]
Since \(\gcd(c_0,d_0)=1\), one has
\[
e(c_0,d_0)\in
\{(\pm q,0),(0,\pm q)\},
\qquad
\lambda_q^{-1}\le\frac{C+D}{q}.
\]
The same divisor bound therefore gives
\[
\begin{aligned}
\frac1B
\sum_{\substack{u,v,\mathbf y,q\\c_0d_0=0}}
\lambda_q^{-1}
\ll_\varepsilon
B^\varepsilon\frac{JP_y(C+D)}{BQ_0}
=B^\varepsilon\frac{C+D}{sH}
=B^\varepsilon s^{-1}(C^{-1}+D^{-1}).
\end{aligned}
\]
Combining the two contributions proves \eqref{ge:four}.
\end{proof}

\subsection{Trilinear recovery and all exceptional sets}


We use Reuss's trilinear estimate
\cite[Theorem~9]{Reuss2015}.
Let \(f\) be an irreducible integral trilinear form with
hyperdeterminant \(\mathcal D\ne0\).
For three centered coordinate boxes, let \(K_1,K_2,K_3\)
be the products of their respective coordinate bounds, each
individual bound being at least \(1\).
The number of solutions of \(f=0\), with each coordinate pair
primitive and all three associated binary determinants nonzero, is
\begin{equation}
\ll_\varepsilon
(\|f\|K_1K_2K_3)^\varepsilon
\left\{
\frac{\sqrt{K_1K_2K_3}}{|\mathcal D|^{1/4}}
+\sqrt{K_1K_2}+K_3
\right\},
\label{ge:Reuss-input}
\end{equation}
where \(\|f\|\) is the largest absolute coefficient.
The coordinate pairs may be permuted.
We write $|\mathcal D|$ in the denominator: the proof of
\cite[Lemma~13]{Reuss2015} treats both real and complex-conjugate
roots of the binary quadratic form.

Retain the chosen shortest vector \(e(c_0,d_0)\).
The axial contribution was bounded in
Proposition~\ref{ge:four-proposition}; hence assume \(c_0d_0\ne0\).
Since \(\Delta\ne0\), we have
\((L_1(c_0,d_0),L_2(c_0,d_0))\ne(0,0)\).
Define
\begin{equation}\label{ge:recover}
 \ell=\gcd(L_1(c_0,d_0),L_2(c_0,d_0))>0\quad\text{and}\quad
 (a',b')=\frac1\ell(L_1(c_0,d_0),L_2(c_0,d_0)).
\end{equation}
By Lemma~\ref{ge:content},
\[
\frac qe\mid\ell,
\qquad
\ell\ge\frac qe,
\qquad
\gcd(u,v)=\gcd(c_0,d_0)=\gcd(a',b')=1.
\]

Suppose \(\lambda_q\le R\).
The scale identities in \eqref{ge:scales} give
\[
UCY_2=VDY_1=Z_aQ,
\qquad
VCY_3=UDY_4=Z_bQ,
\]
and therefore
\[
\begin{aligned}
|L_1(c_0,d_0)|
&\ll (UCY_2+VDY_1)\frac Re
\ll \frac{Z_aQR}{e},\\
|L_2(c_0,d_0)|
&\ll (VCY_3+UDY_4)\frac Re
\ll \frac{Z_bQR}{e}.
\end{aligned}
\]
Since \(q\asymp sQ\), this yields
\begin{equation}
\begin{gathered}
|c_0|\le\frac{CR}{e}\le CR,
\qquad
|d_0|\le\frac{DR}{e}\le DR,\\
|a'|\ll\frac{Z_aQR}{e\ell}
\le\frac{Z_aQR}{q}\ll\frac{Z_aR}{s},
\qquad
|b'|\ll\frac{Z_bQR}{e\ell}
\le\frac{Z_bQR}{q}\ll\frac{Z_bR}{s}.
\end{gathered}
\label{ge:recovery-box}
\end{equation}

The three primitive pairs satisfy
\begin{equation}
f=b'(uc_0y_2+vd_0y_1)
-a'(vc_0y_3+ud_0y_4)=0.
\label{ge:trilinear}
\end{equation}
Fixing \((u,v)\), its bilinear coefficient matrix has determinant
\[
\det
\begin{pmatrix}
-vy_3&uy_2\\
-uy_4&vy_1
\end{pmatrix}
=u^2y_2y_4-v^2y_1y_3=\Delta.
\]
Thus \cite[Lemma~5]{Reuss2015} gives
\begin{equation}
\mathcal D
=\operatorname{disc}_{u,v}(\Delta)
=4y_1y_2y_3y_4\ne0,
\qquad
|\mathcal D|\asymp P_y.
\label{ge:hyperdet}
\end{equation}
Moreover, \(\gcd(y_1,y_2,y_3,y_4)=1\), so \(f\) is primitive.
Its irreducibility follows from
\(\mathcal D\ne0\) and \cite[Lemmas~5 and~7]{Reuss2015}.

Fixing $\mathbf y$ and the recovered three primitive pairs leaves at most
$\#\{q:q\mid\Delta\}\le\tau(|\Delta|)\ll_\varepsilon B^\varepsilon$
original tuples: $u,v,\mathbf y,q$ determine the chosen shortest vector
and its content $e$. Thus the recovery map has this multiplicity.

First consider \(a'b'=0\).
By \eqref{ge:short-gcd}, we have
$\gcd(uy_2,vy_1)=\gcd(uy_4,vy_3)=1.$
Primitivity of \((c_0,d_0)\) therefore gives
\[
\begin{array}{c|c|c}
 & (c_0,d_0)
 & \bigl(L_1(c_0,d_0),L_2(c_0,d_0)\bigr)\\ \hline
a'=0 & \pm(vy_1,-uy_2) & (0,\mp\Delta)\\
b'=0 & \pm(uy_4,-vy_3) & (\pm\Delta,0).
\end{array}
\]
Since \(q\mid\Delta\), these primitive vectors belong to
\(\Lambda_q\). Shortestness forces \(e=1\).
Hence
\[
\lambda_q\asymp
\begin{cases}Z_aQ/H,&a'=0,\\ Z_bQ/H,&b'=0.\end{cases}
\]
There are $O(JP_y)$ choices of $(u,v,\mathbf y)$ and at most
$\tau(|\Delta|)\ll_\varepsilon B^\varepsilon$ moduli for each.
Writing subscripts for the imposed condition and using $B=HJP_y/Q$,
\begin{equation}
\begin{aligned}
\frac{\mathcal G_{a'=0}+\mathcal G_{b'=0}}B
\ll_\varepsilon B^\varepsilon\frac{JP_yH}{BQ}
\left(Z_a^{-1}+Z_b^{-1}\right)
=B^\varepsilon\left(Z_a^{-1}+Z_b^{-1}\right).
\end{aligned}
\label{ge:zero-recovery}
\end{equation}

The three associated binary determinants are
\[
\Delta,\qquad
\eta=c_0^2y_2y_3-d_0^2y_1y_4,\qquad
b'^2y_1y_2-a'^2y_3y_4.
\]
By \eqref{ge:recover},
\begin{equation}
\begin{aligned}
\ell^2(b'^2y_1y_2-a'^2y_3y_4)
=y_1y_2L_2(c_0,d_0)^2
-y_3y_4L_1(c_0,d_0)^2
=-\Delta\eta.
\end{aligned}
\label{ge:det-identity}
\end{equation}
Since \(\Delta\ne0\), all three determinants are nonzero
if and only if \(\eta\ne0\).
On \(a'b'\eta\ne0\), all recovered coordinates are nonzero;
a coordinate bound below \(1\) gives an empty box.
Thus \eqref{ge:Reuss-input} applies to every nonempty such box.

It remains to bound \(\eta=0\).
On this locus,
\[
c_0^2|y_2y_3|=d_0^2|y_1y_4|,
\qquad
\gcd(c_0^2,d_0^2)
=\gcd(|y_1y_4|,|y_2y_3|)=1.
\]
Equality of the reduced positive fractions gives
$c_0^2=|y_1y_4|$ and $d_0^2=|y_2y_3|$. Pairwise coprimality therefore yields
\[
|y_i|\text{ is a square}\quad(1\le i\le4)\quad\text{and}\quad
\#\{\mathbf y:|y_i|\asymp Y_i,\ |y_i|\text{ all squares}\}
\ll\prod_iY_i^{1/2}=\sqrt{P_y}.
\]
Furthermore, we have
\[
\begin{aligned}
\lambda_q
=e\max\left\{\frac{|c_0|}{C},\frac{|d_0|}{D}\right\}
\ge\sqrt{\frac{|c_0d_0|}{H}}
=\frac{|y_1y_2y_3y_4|^{1/4}}{\sqrt H}
\asymp\frac{P_y^{1/4}}{\sqrt H}.
\end{aligned}
\]
Using \(O(J)\) choices of \((u,v)\), the divisor bound for
\(q\mid\Delta\ne0\), and \eqref{ge:identities}, we conclude that
\begin{equation}
\begin{aligned}
\frac{\mathcal G_{\eta=0}}B
\ll_\varepsilon
B^\varepsilon
\frac{J\sqrt{P_y}}B
\frac{\sqrt H}{P_y^{1/4}}=
B^\varepsilon\frac{Q}{\sqrt H\,P_y^{3/4}}
=
B^\varepsilon\frac{\sqrt{J/A}}{\sqrt{P_y}}.
\end{aligned}
\label{ge:square-exception}
\end{equation}
The last equality uses
\[
P_y^{1/4}=\frac{T}{\sqrt{AHJ}},
\qquad Q=\frac TA.
\]
This bounds the exceptional locus arising from square coefficients.

\begin{proposition}
\label{ge:reuss-proposition}
Under the hypotheses of Proposition~\ref{ge:four-proposition},
for every \(\varepsilon>0\),
\begin{equation}
\begin{aligned}
\frac{\mathcal G}{B}\ll_\varepsilon B^\varepsilon
\biggl\{
s^{-1/2}\sqrt{\frac QH}
+\frac Q{\sqrt{HJ}}
+s^{-1}\frac T{H\sqrt{AJ}}
+s^{-1}(C^{-1}+D^{-1})
+Z_a^{-1}+Z_b^{-1}
+\frac{\sqrt{J/A}}{\sqrt{P_y}}
\biggr\}.
\end{aligned}
\label{ge:reuss}
\end{equation}
The implied constant is independent of \(s\) and the individual
scales, with the interval constants fixed.
\end{proposition}

\begin{proof}
Let \(\mathcal A_{\rm gen}(R)\) count the tuples in
\(\mathcal A(R)\) satisfying \(a'b'\eta\ne0\), with the recovered
quantities defined using the fixed shortest-vector rule.
Empty ranges contribute zero.

First suppose \(0<R\le R_{\max}\), where \(R_{\max}\) is as in
\eqref{ge:layer-cake}.
The recovery argument preceding this proposition gives
\[
\gcd(u,v)=\gcd(c_0,d_0)=\gcd(a',b')=1,
\qquad
\Delta\eta\ne0,
\]
and \eqref{ge:det-identity} ensures that the third binary
determinant is also nonzero.
Apply \eqref{ge:Reuss-input} to the ordered pairs
 $(u,v),$ $(c_0,d_0)$ and $(a',b')$.
By \eqref{ge:recovery-box} and \eqref{ge:hyperdet}, we have
\[
 K_1\ll J,\qquad K_2\ll HR^2,\qquad
 K_3\ll AR^2/s^2,\qquad |\mathcal D|\asymp P_y.
\]
Here the actual coordinate bounds include the fixed interval
constants. Since all six coordinates are nonzero, every nonempty
box has each coordinate bound at least \(1\).

Nonemptiness and the scale bounds in Section~\ref{ge:scales-subsection} imply
\[
sQ=Q_0\gg1,\qquad
s^{-1}\ll Q\ll B,\qquad
R_{\max}\asymp\sqrt{\frac{sQ}{H}}\ll B^{1/2}.
\]
Thus \(\|f\|K_1K_2K_3\) has polynomial size in \(B\),
uniformly in the scales and \(s\).

There are \(O(P_y)\) choices of \(\mathbf y\).
Apply \eqref{ge:Reuss-input}, namely
\cite[Theorem~9]{Reuss2015}, together with the previously established
recovery multiplicity, using sufficiently small exponents.
Since \(P_y^{3/4}\sqrt{AHJ}=B\), we obtain
\begin{equation}
\begin{aligned}
\mathcal A_{\rm gen}(R)
&\ll_\varepsilon B^{\varepsilon/2}P_y
\left\{
\frac{\sqrt{J(HR^2)(AR^2/s^2)}}{P_y^{1/4}}
+\sqrt{JH}\,R+\frac{AR^2}{s^2}
\right\}\\
&=B^{\varepsilon/2}
\left\{
\frac BsR^2
+P_y\sqrt{JH}\,R
+\frac{P_yA}{s^2}R^2
\right\}.
\end{aligned}
\label{ge:A-reuss}
\end{equation}
For \(R\ge R_{\max}\), we have
$\mathcal A_{\rm gen}(R)=\mathcal A_{\rm gen}(R_{\max}).$
The right-hand side of \eqref{ge:A-reuss} increases with \(R\);
hence that estimate holds for every \(R>0\).

The divisor bound \eqref{ge:divisor-cap} also gives
$\mathcal A_{\rm gen}(R)
\le\mathcal A(R)
\ll_\varepsilon B^{\varepsilon/2}JP_y.$
Moreover, for $0<R<\frac1{\max(C,D)}$, we have
$\mathcal A_{\rm gen}(R)=0.$
Apply the infinite-integral identity in \eqref{ge:layer-cake}
and the integral estimates proved in
Proposition~\ref{ge:four-proposition}.
Combining the divisor bound with each term of
\eqref{ge:A-reuss} yields
\[
\begin{aligned}
&\sum_{\substack{u,v,\mathbf y,q\\c_0d_0a'b'\eta\ne0}}
\lambda_q^{-1}
=\int_0^\infty\frac{\mathcal A_{\rm gen}(R)}{R^2}\,dR\\
&\quad\ll_\varepsilon B^{\varepsilon/2}
\biggl\{
\sqrt{\frac{BJP_y}{s}}
+P_y\sqrt{JH}
 \log\left(2+\max(C,D)\sqrt{\frac JH}\right)
+\frac{P_y\sqrt{AJ}}s
\biggr\}.
\end{aligned}
\]
Using \(C,D,U,V\ll B^2\) and \(C,D\ge1/2\)
from Section~\ref{ge:scales-subsection}, together with
\(H=CD\) and \(J=UV\), we obtain
\[
\log\left(2+\max(C,D)\sqrt{\frac JH}\right)
\ll\log(2B)\ll_\varepsilon B^{\varepsilon/2}.
\]
Consequently,
\[
\sum_{\substack{u,v,\mathbf y,q\\c_0d_0a'b'\eta\ne0}}
\lambda_q^{-1}
\ll_\varepsilon B^\varepsilon
\left\{
\sqrt{\frac{BJP_y}{s}}
+P_y\sqrt{JH}
+\frac{P_y\sqrt{AJ}}s
\right\}.
\]

Divide by \(B\) and use \eqref{ge:identities}; in particular,
$\frac{P_y\sqrt{AJ}}{sB}=\frac{T}{sH\sqrt{AJ}}.$
This gives the first three terms of \eqref{ge:reuss}.
The remaining tuples lie on
\[
c_0d_0=0,\qquad a'b'=0,\qquad\text{or}\qquad\eta=0.
\]
Adding the axial estimate from
Proposition~\ref{ge:four-proposition}, the bounds
\eqref{ge:zero-recovery}, and \eqref{ge:square-exception}
proves \eqref{ge:reuss}.
\end{proof}

\subsection{A geometric error bound when
\texorpdfstring{$P_0\ge B^{1+\delta}$}{P0 >= B^(1+delta)}}

\begin{corollary}
\label{ge:above-B}
Assume the hypotheses of Proposition~\ref{ge:four-proposition},
with \(A\ge H\ge J\).
Suppose that, for fixed \(\delta>0\),
\[
\begin{gathered}
 P_0\ge B^{1+\delta},\quad
 \min_e Z_e\ge B^\delta,\quad \min_iY_i\ge B^\delta,\quad
 B^{-\kappa}\le s\ll1,
\end{gathered}
\]
where \(0<\kappa\le\delta/12\). Then
\begin{equation}
\mathcal G\ll_\delta B^{1-\delta/24}.
\label{ge:above-B-bound}
\end{equation}
\end{corollary}

\begin{proof}
Using \(P_0=AHJ\), \(Q=(BP_0)^{1/3}/A\), and
\(P_y=\prod_{i=1}^4Y_i\), together with $A\ge H\ge J$, we obtain
\[
\begin{aligned}
s^{-1/2}\sqrt{Q/H}
&=s^{-1/2}(B/P_0)^{1/6}(J^2/(AH))^{1/6}
\le B^{-\delta/6+\kappa/2}\le B^{-\delta/8},\\
Q/\sqrt{HJ}
&=(B/P_0)^{1/3}(HJ/A^2)^{1/6}\le B^{-\delta/3}.
\end{aligned}
\]
\[
\begin{aligned}
s^{-1}(C^{-1}+D^{-1})+Z_a^{-1}+Z_b^{-1}
+\sqrt{J/A}/\sqrt{P_y}
\ll B^{-\delta+\kappa}+B^{-\delta}+B^{-2\delta}
\ll B^{-\delta/8}.
\end{aligned}
\]
The remaining terms of \eqref{ge:four} satisfy
\[
s\frac{Q^2}{\sqrt{HJP_y}}+\frac Q{\sqrt H\,P_y^{1/4}}
=s\frac{\sqrt{HJ}}A+\sqrt{J/A}\ll\sqrt{J/A}.
\]
Taking the better of \eqref{ge:four} and \eqref{ge:reuss},
\[
\frac{\mathcal G}{B}\ll_\varepsilon B^\varepsilon
\left\{B^{-\delta/8}
+\min\left(\sqrt{J/A},\frac{T}{sH\sqrt{AJ}}\right)\right\}.
\]
Finally, $T^3=BP_0$ and $P_0=AHJ$ imply
\[
\begin{aligned}
\min\left(\sqrt{J/A},\frac{T}{sH\sqrt{AJ}}\right)
&\le (J/A)^{1/8}\left(\frac{T}{sH\sqrt{AJ}}\right)^{3/4}\\
&=s^{-3/4}(B/P_0)^{1/4}(J/H)^{1/4}
\le B^{-\delta/4+3\kappa/4}\le B^{-3\delta/16}.
\end{aligned}
\]
Thus $\mathcal G/B\ll_\varepsilon B^{\varepsilon-\delta/8}$;
taking $\varepsilon=\delta/12$ proves \eqref{ge:above-B-bound}.
\end{proof}

\section{Complete coprimality and the two-dimensional main term}
\label{ar:section}

The main result of this section is Corollary~\ref{ar:count}, which
expresses the fully coprime integer count as a continuous main term
plus separate arithmetic and geometric errors. The argument proceeds
in three stages. Proposition~\ref{ar:long} performs the long-variable
sieve and the lattice count. The short-variable Boolean sieve and
the character estimates then evaluate the sum of its main terms in
Theorem~\ref{ar:main}. Finally, the coarea calculation identifies
the resulting integral with the common continuous measure.
The estimates of Section~\ref{ge:section} control the geometric error
when Corollary~\ref{ar:count} is applied in Section~\ref{gl:section}.

Throughout these stages the recovered variables remain inside the
original height weight. In Section~\ref{ar:weighted-characters},
the common-prefix identity~\eqref{ar:common-prefix} separates this
joint weight from character sums over intervals independent of the
modulus. Its derivative bound is taken before the quadratic-character
average. This order is needed to apply the fixed-interval moment
estimates uniformly to the recovered weights.

\subsection{Notation and exact coprimality conditions}
\label{ar:notation}

Write
\[
a=z_{12},\quad b=z_{34},\quad c=z_{14},\quad d=z_{23},
\quad u=z_{13},\quad v=z_{24},
\]
and fix
\[
y_i=\sigma_i t_i,\qquad
\sigma_i\in\{-1,1\},\quad t_i\in\mathbb Z_{>0}.
\]
Define
\begin{equation}
\begin{aligned}
L_1(c,d)=\sigma_2ut_2c+\sigma_1vt_1d,\quad
L_2(c,d)=-\sigma_3vt_3c-\sigma_4ut_4d,
\end{aligned}
\label{ar:linear}
\end{equation}
and
\begin{equation}
\Delta
=-\sigma_2\sigma_4u^2t_2t_4
+\sigma_1\sigma_3v^2t_1t_3.
\label{ar:delta}
\end{equation}
Then
\[
\sum_{i=1}^{4}A_i(z)y_i=bL_1(c,d)-aL_2(c,d),
\]
so the torsor equation is \(bL_1=aL_2\).

By \cite[Section~2, Lemma~1]{HB2003}, the canonical coprimality
conditions involving only
\(\mathbf t=(u,v,t_1,t_2,t_3,t_4)\) are
\begin{equation}
\gcd(t_i,t_j)=1\quad(i<j),\qquad
\gcd(u,v)=\gcd(u,t_1t_3)=\gcd(v,t_2t_4)=1.
\label{ar:short-gcd}
\end{equation}
Write \(\mathbf1_{\mathcal C}(\mathbf t)\) for their indicator.
Divisor tuples use the coordinate order
\((u,v,t_1,t_2,t_3,t_4)\).

For later use, put
\[
F_1=uvt_1t_4,\qquad
F_2=uvt_2t_3,\qquad
w(n)=\prod_{p\mid n}\frac{p}{p+1},
\]
and
\begin{equation}
\rho(\mathbf t)
=\frac1{\zeta(2)}
\prod_{\substack{p\mid F_1F_2\\p\nmid uv}}\frac{p}{p+1}
\prod_{p\mid uv}\frac{p-1}{p+1}.
\label{ar:rho}
\end{equation}
Here $\rho$ is the arithmetic weight. In this section, $\eta$ denotes
the relative modulus scale; the determinant denoted by $\eta$ in
Section~\ref{ge:section} will not recur.

Under \eqref{ar:short-gcd}, we have
$\gcd(u^2t_2t_4,v^2t_1t_3)=1$ and
$\gcd(t_1t_4,t_2t_3)=1.$
Consequently,
\begin{equation}
\gcd(\Delta,uvt_1t_2t_3t_4)=1,
\qquad
\gcd(F_1,F_2)=uv.
\label{ar:units}
\end{equation}
Since $\gcd(u^2t_2t_4,v^2t_1t_3)=1$, one has
\[
\Delta=0\iff u=v=t_1=t_2=t_3=t_4=1
\quad\text{and}\quad\sigma_1\sigma_3=\sigma_2\sigma_4.
\]
Every divisor sum below is restricted to $\Delta\ne0$.
For a positive divisor \(m\mid|\Delta|\), define
\begin{equation}
\Lambda_m
=
\{(c,d)\in\mathbb Z^2:
m\mid L_1(c,d),\ m\mid L_2(c,d)\}.
\label{ar:lattice}
\end{equation}
The coefficient matrix of \(L_1,L_2\) has determinant \(\Delta\),
and all its entries are coprime to \(m\) by \eqref{ar:units}.
Hence \cite[Lemma~3]{Reuss2015} gives \(\det\Lambda_m=m\);
the case \(m=1\) is immediate.
With the original scales \(C,D>0\), set
\[
\lambda_m
=
\min_{\substack{(c,d)\in\Lambda_m\\(c,d)\ne(0,0)}}
\max\left\{\frac{|c|}{C},\frac{|d|}{D}\right\}.
\]

For positive integers \(c,d\), impose the remaining coprimality
conditions and define \(q\) by
\begin{equation}\label{ar:long-gcd}
\gcd(c,d)=\gcd(c,F_1)=\gcd(d,F_2)=1,\quad
q=\gcd(L_1(c,d),L_2(c,d))>0.
\end{equation}
Positive recovery requires \(L_1(c,d)L_2(c,d)>0\), in which case
\begin{equation}
a=\frac{\nu L_1(c,d)}q,\qquad
b=\frac{\nu L_2(c,d)}q,\qquad
\nu=\operatorname{sgn}L_1(c,d)
=\operatorname{sgn}L_2(c,d).
\label{ar:recover}
\end{equation}
These formulas give \(a,b\in\mathbb Z_{>0}\),
\(\gcd(a,b)=1\), and \(bL_1=aL_2\).

The identities
$\Delta c=-\sigma_4ut_4L_1-\sigma_1vt_1L_2$ and
$\Delta d=\sigma_3vt_3L_1+\sigma_2ut_2L_2$
show that
\[
q\mid\Delta c,\qquad q\mid\Delta d,\qquad \gcd(c,d)=1
\quad\Longrightarrow\quad q\mid|\Delta|.
\]
In particular, \eqref{ar:units} implies
$\gcd(q,uvt_1t_2t_3t_4)=1.$

To verify the remaining canonical coprimalities, reduce
\eqref{ar:linear} modulo a prime \(p\).
Under \eqref{ar:short-gcd} and \eqref{ar:long-gcd},
every displayed residue in the following table is nonzero:
\[
\begin{array}{c|cc}
p\mid & L_1\pmod p & L_2\pmod p\\ \hline
c   & \sigma_1vt_1d & -\sigma_4ut_4d\\
d   & \sigma_2ut_2c & -\sigma_3vt_3c\\
u   & \sigma_1vt_1d & -\sigma_3vt_3c\\
v   & \sigma_2ut_2c & -\sigma_4ut_4d\\
t_1 & \sigma_2ut_2c & \\
t_2 & \sigma_1vt_1d & \\
t_3 &                & -\sigma_4ut_4d\\
t_4 &                & -\sigma_3vt_3c
\end{array}
\]
Thus
\[
\begin{gathered}
\gcd(L_1L_2,cduv)=1,\qquad
\gcd(L_1,t_1t_2)=\gcd(L_2,t_3t_4)=1,\\
\gcd(ab,cduv)=1,\qquad
\gcd(a,t_1t_2)=\gcd(b,t_3t_4)=1.
\end{gathered}
\]
Together with \eqref{ar:short-gcd}, \eqref{ar:long-gcd}, and
$\gcd(a,b)=1$, these are exactly \eqref{in:coprimality}.

Conversely, the torsor equation and \eqref{in:coprimality} imply
\eqref{ar:short-gcd} and
$\gcd(c,d)=\gcd(c,F_1)=\gcd(d,F_2)=\gcd(a,b)=1$.
Since $\Delta\ne0$, the forms cannot both vanish, and
\[
bL_1=aL_2,\quad\gcd(a,b)=1
\Longrightarrow \frac{L_1}{a}=\frac{L_2}{b}\in\mathbb Z\setminus\{0\},\quad\text{and}\quad
q=\gcd(L_1,L_2)=\left|\frac{L_1}{a}\right|
=\left|\frac{L_2}{b}\right|.
\]
This recovers \eqref{ar:long-gcd} and \eqref{ar:recover}.
The overall sign convention, height restriction, and pair-sum
exclusions remain imposed separately.
\subsection{A lattice estimate away from the origin}
\label{ar:origin-subsection}

\begin{lemma}
\label{ar:origin}
Let \(L\subset\mathbb R^2\) be a full-rank lattice, and let
\(\lambda(L)\) be the length of its shortest nonzero vector,
measured in either the Euclidean or the maximum norm.
Suppose that \(f\in C_c^1(\mathbb R^2)\) is supported in a fixed
bounded set and satisfies \(f(0)=0\). Then
\[
\sum_{\mathbf x\in L}f(\mathbf x)
=
\frac1{\det L}\int_{\mathbb R^2}f(\mathbf x)\,d\mathbf x
+
O\left(\frac{\|f\|_{C^1}}{\lambda(L)}\right),
\]
with an implied constant depending only on the prescribed support set.
\end{lemma}

\begin{proof}
Use the Euclidean norm, write $\lambda=\lambda(L)$, and fix $R\ge1$
containing the prescribed support in $\{\mathbf x:|\mathbf x|_2\le R\}$.
A shortest vector extends to a basis of $L$. After rotation and reduction
of the second basis vector,
\[
\begin{gathered}
L=\{(m\lambda+n\beta,nh):m,n\in\mathbb Z\},\quad
|\beta|\le\lambda/2,\quad
h\ge\sqrt3\lambda/2,\quad\det L=\lambda h.
\end{gathered}
\]
The lower bound on $h$ follows from $\beta^2+h^2\ge\lambda^2$.
Euler's summation formula \cite[Theorem~1.2]{SS} gives
\[
\sum_{m\in\mathbb Z}g(a+m\ell)
=\ell^{-1}\int g+O\left(\int|g'|\right)
\quad(g\in C_c^1(\mathbb R)).
\]
Apply it first to each horizontal row, of which $O_R(1+h^{-1})$
meet the support, and then to $y\mapsto\int f(x,y)\,dx$. This gives
\[
\begin{aligned}
\sum_{\mathbf x\in L}f(\mathbf x)
=\frac1{\lambda h}\int f
+O_R\!\left(\|f\|_{C^1}(1+h^{-1}+\lambda^{-1})\right)
=\frac1{\det L}\int f
+O_R\!\left(\|f\|_{C^1}(1+\lambda^{-1})\right).
\end{aligned}
\]
For $\lambda\le R$, this is the required error.
If $\lambda>R$, only the origin can meet the support and $f(0)=0$, so
\[
\left|\sum_{\mathbf x\in L}f(\mathbf x)-\frac1{\det L}\int f\right|
\ll_R\frac{\|f\|_\infty}{\lambda^2}
\ll_R\frac{\|f\|_{C^1}}\lambda.
\]
Finally, $\|\mathbf x\|_\infty\le|\mathbf x|_2
\le\sqrt2\|\mathbf x\|_\infty$ gives the maximum-norm version.
\end{proof}

\subsection{The complete long-variable sieve}
\label{ar:long-subsection}

For \(W\in C_c^1(\mathbb R^2)\), write
\[
W_*(x_1,x_2)=W(Cx_1,Dx_2),\qquad
\mathcal W_1=\|W_*\|_{C^1},\qquad
I_W=\int_{\mathbb R^2}W(c,d)\,dc\,dd.
\]
Assume that, for fixed \(0<c_0<c_1\), the support of \(W_*\) lies
in
$\{(x_1,x_2):c_0\le |x_1|,|x_2|\le c_1\}.$
In particular, it avoids both coordinate axes.  Throughout the estimates
below, all integer parameters and scales that occur are bounded by
fixed powers of a parameter \(B\ge2\).  Constants may depend on these
fixed powers and on the fixed support, but not on individual divisors.

\begin{proposition}
\label{ar:long}
Assume \eqref{ar:short-gcd}, \(\Delta\ne0\), and the support and
polynomial-size conditions stated above. For \(q\mid|\Delta|\), define
\[
N_W(q)=
\sum_{c,d\in\mathbb Z}W(c,d)\,
\mathbf1_{\gcd(c,d)=1}
\mathbf1_{\gcd(c,F_1)=1}
\mathbf1_{\gcd(d,F_2)=1}
\mathbf1_{\gcd(L_1(c,d),L_2(c,d))=q}.
\]
Then, for every \(\varepsilon>0\),
\begin{equation}
N_W(q)
=
I_W\rho(\mathbf t)
\sum_{k\mid|\Delta|/q}\frac{\mu(k)w(qk)}{qk}
+
O_\varepsilon\!\left(
\frac{B^\varepsilon\mathcal W_1}{\lambda_q}
\right).
\label{ar:long-formula}
\end{equation}
The implied constant depends only on \(\varepsilon\), the prescribed
normalized support, and the fixed polynomial bounds.
\end{proposition}

\begin{proof}
The adjugate identities in Section~\ref{ar:notation} give
\[
\gcd(c,d)=1
\quad\Longrightarrow\quad
\gcd(L_1(c,d),L_2(c,d))\mid|\Delta|.
\]
Hence, on this primitive set, M\"obius inversion gives
\begin{equation}
\mathbf1_{\gcd(L_1,L_2)=q}
=
\sum_{k\mid|\Delta|/q}
\mu(k)\mathbf1_{(c,d)\in\Lambda_{qk}}.
\label{ar:gcd-filter}
\end{equation}
For \(m\mid|\Delta|\), put
\[
M_m(W)=
\sum_{(c,d)\in\Lambda_m}W(c,d)\,
\mathbf1_{\gcd(c,d)=1}
\mathbf1_{\gcd(c,F_1)=1}
\mathbf1_{\gcd(d,F_2)=1}.
\]
Thus
$N_W(q)=\sum_{k\mid|\Delta|/q}\mu(k)M_{qk}(W).$

Set
\[
L_{\mathrm{div}}=c_1\min(C,D),\quad
K_n(e,f)=\Lambda_n\cap(e\mathbb Z\times f\mathbb Z),\quad
n\mid|\Delta|,\qquad e\mid F_1,\quad f\mid F_2.
\]
By \eqref{ar:units}, \(\gcd(n,ef)=1\).
After substituting \(c=ex,\ d=fy\), the coefficient matrix has
determinant \(ef\Delta\), divisible by \(n\), and all its entries
are coprime to \(n\). Applying the same lemma to this matrix and
accounting for the scaling by \(e,f\) gives \(\det K_n(e,f)=efn\);
the case \(n=1\) is immediate.

Under the conditions \(\gcd(c,F_1)=\gcd(d,F_2)=1\),
$\gcd\bigl(\gcd(c,d),F_1F_2\bigr)=1.$
Expanding the three coprimality indicators and writing
\(r=\gcd(m,\ell)\) inside the sums gives
\begin{equation*}
M_m(W)
=
\sum_{e\mid F_1}\sum_{f\mid F_2}\mu(e)\mu(f)
\sum_{\substack{\ell\le L_{\mathrm{div}}\\\gcd(\ell,F_1F_2)=1}}
\mu(\ell)
\sum_{\mathbf x\in K_{m/r}(e,f)}W(\ell\mathbf x).
\end{equation*}
Indeed, $m\mid\ell L_i(\mathbf x)
\iff
m/r\mid L_i(\mathbf x),$
and the support condition implies
\[
W(\ell\mathbf x)\ne0
\quad\Longrightarrow\quad
x_1x_2\ne0,\qquad
\ell\le c_1C,\quad \ell\le c_1D.
\]
Thus \(\ell\le L_{\mathrm{div}}\) is an exact support restriction, including
when \(L_{\mathrm{div}}<1\).

The inclusion $rK_{m/r}(e,f)\subseteq\Lambda_m$ gives
\[
\min_{\mathbf x\in\ell K_{m/r}(e,f)\setminus\{0\}}
\max\left(\frac{|x_1|}{C},\frac{|x_2|}{D}\right)
\ge\frac{\ell\lambda_m}{r}\quad\text{and}\quad
\det\bigl(\ell K_{m/r}(e,f)\bigr)=\frac{efm\ell^2}{r}.
\]
Lemma~\ref{ar:origin}, after normalization by $C,D$, therefore yields
\begin{equation*}
\sum_{\mathbf x\in K_{m/r}(e,f)}W(\ell\mathbf x)
=\frac{I_Wr}{efm\ell^2}
+O\left(\frac{\mathcal W_1r}{\ell\lambda_m}\right).
\end{equation*}
For every $L_{\mathrm{div}}>0$,
\[
\sum_{\ell\le L_{\mathrm{div}}}\frac{\gcd(m,\ell)}\ell
\le\sum_{r\mid m}\sum_{h\le L_{\mathrm{div}}/r}\frac1h
\ll\tau(m)\log(2+L_{\mathrm{div}}).
\]
Thus the total lattice error is
$O(\mathcal W_1\lambda_m^{-1}\tau(F_1)\tau(F_2)\tau(m)\log(2+L_{\mathrm{div}}))$.
To extend the main term to all $\ell\ge1$, use
\begin{equation*}
\sum_{\ell>L_{\mathrm{div}}}\frac{\gcd(m,\ell)}{\ell^2}
\le\sum_{r\mid m}\frac1r\sum_{h>L_{\mathrm{div}}/r}\frac1{h^2}
\ll\frac{\tau(m)}{L_{\mathrm{div}}}.
\end{equation*}
This also holds for $L_{\mathrm{div}}<1$. Since
$|I_W|\ll CD\mathcal W_1$ and $(m,0),(0,m)\in\Lambda_m$,
\[
\begin{aligned}
\lambda_m\le\frac m{\max(C,D)},\qquad
\frac{|I_W|\tau(m)}{efmL_{\mathrm{div}}}
\ll\frac{\mathcal W_1\tau(m)\max(C,D)}{efm}\le\frac{\mathcal W_1\tau(m)}{ef\lambda_m}.
\end{aligned}
\]
Summing over $e,f$ gives the same total error bound.

It remains to evaluate the absolutely convergent main term.
By \eqref{ar:units}, we have
$\gcd(m,F_1F_2)=1$ and $\gcd(F_1,F_2)=uv$.
Hence
\[
\begin{aligned}
\sum_{\substack{\ell\ge1\\\gcd(\ell,F_1F_2)=1}}
\frac{\mu(\ell)\gcd(m,\ell)}{\ell^2}
=
\prod_{p\nmid mF_1F_2}(1-p^{-2})
\prod_{p\mid m}(1-p^{-1})
=
\frac{w(m)}{\zeta(2)}
\prod_{p\mid F_1F_2}(1-p^{-2})^{-1}.
\end{aligned}
\]
Therefore
\begin{equation*}
\begin{aligned}
&
\left(\sum_{e\mid F_1}\frac{\mu(e)}e\right)
\left(\sum_{f\mid F_2}\frac{\mu(f)}f\right)
\sum_{\substack{\ell\ge1\\\gcd(\ell,F_1F_2)=1}}
\frac{\mu(\ell)\gcd(m,\ell)}{\ell^2}\\
&\qquad=
\frac{w(m)}{\zeta(2)}
\prod_{\substack{p\mid F_1F_2\\p\nmid uv}}\frac{p}{p+1}
\prod_{p\mid uv}\frac{p-1}{p+1}
=\rho(\mathbf t)w(m).
\end{aligned}
\end{equation*}
Combining these estimates gives
\begin{equation}
M_m(W)
=
\frac{I_W\rho(\mathbf t)w(m)}m
+
O\!\left(
\frac{\mathcal W_1}{\lambda_m}
\tau(F_1)\tau(F_2)\tau(m)\log(2+L_{\mathrm{div}})
\right).
\label{ar:long-intermediate}
\end{equation}

Finally,
\[
\Lambda_{qk}\subseteq\Lambda_q,\qquad
\lambda_{qk}\ge\lambda_q,\qquad
\tau(qk)\le\tau(|\Delta|)
\quad(k\mid|\Delta|/q).
\]
Substituting \eqref{ar:long-intermediate} into
\eqref{ar:gcd-filter}, the total error is
$\ll
\frac{\mathcal W_1}{\lambda_q}
\tau(F_1)\tau(F_2)\tau(|\Delta|)^2\log(2+L_{\mathrm{div}}).$
The assumed polynomial bounds and the standard divisor estimate imply
\[
\tau(F_1)\tau(F_2)\tau(|\Delta|)^2\log(2+L_{\mathrm{div}})
\ll_\varepsilon B^\varepsilon,
\]
which proves \eqref{ar:long-formula}.
\end{proof}

\subsection{The short-variable Boolean sieve}
\label{ar:boolean-subsection}

Retain all the conditions in \eqref{ar:short-gcd}.
Identify the indices \(1,\ldots,6\) with
$\mathbf t=(u,v,t_1,t_2,t_3,t_4),$
in this order. For \(S\subseteq\{1,\ldots,6\}\), the local factor
of \(\mathbf1_{\mathcal C}\rho\), excluding \(\zeta(2)^{-1}\), is
\[
F_p(S)=
\begin{cases}
1,&S=\varnothing,\\[1mm]
\dfrac{p}{p+1},&S=\{j\},\quad 3\le j\le6,\\[2mm]
\dfrac{p-1}{p+1},
&S\in\bigl\{
\{1\},\{2\},\{1,4\},\{1,6\},\{2,3\},\{2,5\}
\bigr\},\\[2mm]
0,&\text{otherwise}.
\end{cases}
\]
Thus the admissible supports consist of the empty set, six
singletons, and four pairs.

Define
$c_p(D)=
\sum_{E\subseteq D}(-1)^{|D|-|E|}F_p(E).$
Finite M\"obius inversion gives
\[
F_p(S)=\sum_{D\subseteq S}c_p(D).
\]
For \(\mathbf d=(d_1,\ldots,d_6)\in\mathbb Z_{>0}^6\), put
\[
c(\mathbf d)=
\begin{cases}
\displaystyle\prod_p c_p\bigl(\{j:p\mid d_j\}\bigr),
&d_1,\ldots,d_6\text{ are squarefree},\\
0,&\text{otherwise}.
\end{cases}
\]
Multiplying the finite local identities yields
\begin{equation}
\mathbf1_{\mathcal C}(\mathbf t)\rho(\mathbf t)
=
\frac1{\zeta(2)}
\sum_{\substack{d_1\mid u,\ d_2\mid v\\
d_{i+2}\mid t_i\ (1\le i\le4)}}
c(\mathbf d).
\label{ar:boolean}
\end{equation}

The coefficients satisfy $c_p(\varnothing)=1$, $|c_p(D)|\le11$, and
\[
\begin{array}{c|c}
D & c_p(D)\\ \hline
\{j\},\ j=1,2 & -2/(p+1)\\
\{j\},\ 3\le j\le6 & -1/(p+1)\\
|D|=2,\ F_p(D)\ne0 & 1/(p+1)\\
|D|=2,\ F_p(D)=0 & -1+O(p^{-1}).
\end{array}
\]
Consequently, uniformly in $j$,
\[
\begin{aligned}
\sum_D|c_p(D)|p^{-|D|}&=1+O(p^{-2}),\\
\sum_D|c_p(D)|p^{-|D|+\frac12\mathbf1_{j\in D}}
&=1+O(p^{-3/2}),\\
\sum_D|c_p(D)|p^{-|D|+\mathbf1_{j\in D}}
&=1+5p^{-1}+O(p^{-2}).
\end{aligned}
\]
In the last line, the coefficient is $2+3$ for $j=1,2$ and $1+4$
otherwise, from the singleton and forbidden pairs containing $j$.
The first two Euler products converge, giving
\begin{equation}
\sum_{\mathbf d}\frac{|c(\mathbf d)|}{\prod_i d_i}<\infty,
\qquad
\sum_{\mathbf d}\frac{|c(\mathbf d)|d_j^{1/2}}{\prod_i d_i}<\infty
\quad(1\le j\le6).
\label{ar:boolean-convergence}
\end{equation}
For every fixed $K>0$, Mertens' product theorem \cite{Mertens} gives
\begin{equation}
\sum_{\substack{\mathbf d\in\mathbb Z_{>0}^6\\d_i\le B^K\ (1\le i\le6)}}
\frac{|c(\mathbf d)|}{\prod_{i\ne j}d_i}
\le\prod_{p\le B^K}\left(1+\frac5p+O(p^{-2})\right)
\ll_K(\log(2B))^5.
\label{ar:missing-denominator}
\end{equation}

For the six short-variable scales, write
\[
\mathbf H=(U,V,Y_1,Y_2,Y_3,Y_4),\qquad
\prod_{j=1}^6H_j=UVP_y,\qquad H_{\min}=\min_jH_j.
\]
The vector $\mathbf H$ is distinct from the scalar $H=CD$.
For a smooth function \(f\), define
\[
\|f\|_{6,*}
=
\max_{|\alpha|\le6}
\left\|
\partial_{\mathbf x}^{\alpha}
\bigl[f(H_1x_1,\ldots,H_6x_6)\bigr]
\right\|_\infty.
\]
All normalized supports lie in a fixed compact subset of
\(\mathbb R_{>0}^6\).

\begin{lemma}
\label{ar:principal}
Suppose \(H_j\ge1\), and suppose that \(m\) and all \(H_j\) are
bounded by fixed powers of \(B\ge2\).
For a smooth compactly supported weight \(\Psi\),
\begin{equation}
\begin{aligned}
\sum_{\mathbf t\in\mathbb Z_{>0}^6}
\mathbf1_{\mathcal C}(\mathbf t)\rho(\mathbf t)
\mathbf1_{\gcd(uvt_1t_2t_3t_4,m)=1}
\Psi(\mathbf t)
=
\mathfrak d(m)\int_{\mathbb R_{>0}^6}\Psi(\mathbf t)\,d\mathbf t
+
O_\varepsilon\!\left(
B^\varepsilon UVP_yH_{\min}^{-1/2}\|\Psi\|_{6,*}
\right),
\end{aligned}
\label{ar:principal-formula}
\end{equation}
where
\begin{equation}
\mathfrak d(m)=
\frac1{\zeta(2)}
\prod_{p\mid m}(1-p^{-1})^6
\prod_{p\nmid m}D_p,
\qquad
D_p=
\frac{(1-p^{-1})^5(1+6p^{-1}+p^{-2})}{1+p^{-1}}.
\label{ar:dp}
\end{equation}
Here \(\mathfrak d(m)\) is the principal-character density.
The same estimate holds uniformly when the sum and integral are
restricted to any rectangular prefix, with error measured at the full
scales $\mathbf H$.
\end{lemma}

\begin{proof}
Choose a fixed \(c\ge1\) such that the normalized support lies
in \((0,c)^6\). First consider intervals
$I_j\subseteq(0,cH_j]$ $(1\le j\le6)$, with $|I_j|$ denoting real length.
For \(\gcd(d_j,m)=1\), ordinary M\"obius inversion
\cite[Chapter~I, Theorem~2.3]{SS} gives
\[
\begin{aligned}
\#\{n\in I_j:d_j\mid n,\ \gcd(n,m)=1\}
=
\sum_{r\mid m}\mu(r)
\#\{n\in I_j:d_jr\mid n\}
=
\frac{|I_j|}{d_j}\frac{\varphi(m)}m
+O\!\left(2^{\omega(m)}\right).
\end{aligned}
\]
If \(\gcd(d_j,m)>1\), this count is zero.

Insert \eqref{ar:boolean}. Every contributing divisor satisfies
\(d_j\le cH_j\). On this fixed divisor range,
\[
\#\{n\in I_j:d_j\mid n,\ \gcd(n,m)=1\}
\ll\frac{H_j}{d_j},
\qquad
\frac{|I_j|}{d_j}\frac{\varphi(m)}m
\ll\frac{H_j}{d_j},
\]
including when the count vanishes.
Telescoping the six factors therefore gives
\[
\begin{aligned}
\sum_{\mathbf t\in
(I_1\times\cdots\times I_6)\cap\mathbb Z_{>0}^6}&
\mathbf1_{\mathcal C}(\mathbf t)\rho(\mathbf t)
\mathbf1_{\gcd(uvt_1t_2t_3t_4,m)=1}\\
=&
\frac{\prod_{j=1}^6|I_j|}{\zeta(2)}
\left(\frac{\varphi(m)}m\right)^6
\sum_{\substack{d_i\le cH_i\ (1\le i\le6)\\
\gcd(d_1\cdots d_6,m)=1}}
\frac{c(\mathbf d)}{\prod_{i=1}^6d_i}\\
&+
O\!\left(
2^{\omega(m)}UVP_y
\sum_{j=1}^6\frac1{H_j}
\sum_{\substack{\mathbf d\\d_i\le cH_i\ (1\le i\le6)}}
\frac{|c(\mathbf d)|}{\prod_{i\ne j}d_i}
\right).
\end{aligned}
\]
By \eqref{ar:missing-denominator}, this error is
$\ll
2^{\omega(m)}UVP_y(\log(2B))^5
\sum_{j=1}^6H_j^{-1}.$

The extension of the main coefficient to all \(\mathbf d\)
is justified by \eqref{ar:boolean-convergence}:
\[
\begin{aligned}
\sum_{\substack{\mathbf d\\
d_j>cH_j\text{ for some }j}}
\frac{|c(\mathbf d)|}{\prod_{i=1}^6d_i}
\le
\sum_{j=1}^6(cH_j)^{-1/2}
\sum_{\mathbf d}
\frac{|c(\mathbf d)|d_j^{1/2}}{\prod_{i=1}^6d_i}
\ll\sum_{j=1}^6H_j^{-1/2}.
\end{aligned}
\]
Thus the extension costs
$
O\!\left(UVP_y\sum_{j=1}^6H_j^{-1/2}\right).$

For \(r=p^{-1}\), finite inversion gives
\[
\sum_Dc_p(D)r^{|D|}
=\sum_SF_p(S)r^{|S|}(1-r)^{6-|S|}
=\frac{(1-r)^5(1+6r+r^2)}{1+r}=D_p.
\]
In particular, \(D_p=1+O(p^{-2})\).
For \(p\mid m\), all six coordinates are required to be units.
Apart from the common factor \(\zeta(2)^{-1}\), the contribution
at such a prime is \((1-p^{-1})^6\). Absolute convergence therefore
identifies the main coefficient with \(\mathfrak d(m)\).

The resulting rectangle error is
\[
\begin{aligned}
\ll
UVP_y\left(
2^{\omega(m)}(\log(2B))^5\sum_{j=1}^6H_j^{-1}
+\sum_{j=1}^6H_j^{-1/2}
\right)
\ll_\varepsilon
B^\varepsilon UVP_yH_{\min}^{-1/2},
\end{aligned}
\]
uniformly in the interval endpoints.

Apply partial summation \cite[Chapter~I, Theorem~1.1]{SS} in each
normalized coordinate. Since
\[
\max_{\alpha\in\{0,1\}^6}
\left\|\partial_{\mathbf x}^{\alpha}\Psi(\mathbf H\mathbf x)\right\|_\infty
\le\|\Psi\|_{6,*},
\]
and all normalized integration ranges have fixed size, the integral and
prefix boundary terms multiply the rectangle error by
$O(\|\Psi\|_{6,*})$. This proves \eqref{ar:principal-formula}.
\end{proof}

\subsection{Character estimates with independently varying interval parameters}
\label{ar:characters-subsection}

Characters are extended by zero on nonunits. Write
$S_H(\chi)=\sum_{n\in I}\chi(n)$ for any integer interval
$I\subseteq[1,cH]$, with $c$ fixed. The upper-endpoint condition is
needed for the product Dirichlet polynomial below; a length bound alone
does not suffice. Scales with a fixed positive lower bound may be replaced
by $\max(1,H)$.

\begin{lemma}
\label{ar:nonquadratic}
Retain the interval conventions above, and put
\(P_y=Y_1Y_2Y_3Y_4\).
Assume that the six scales are at least \(1\) and that, together
with \(m\), they are bounded by fixed powers of \(B\ge2\).
Then, for every \(\varepsilon>0\),
\begin{equation}\label{ar:nq}
\begin{aligned}
&\frac1{\varphi(m)}
\sum_{\substack{\chi\bmod m\\ \chi^2\ne\chi_0}}
|S_U(\chi^2)S_V(\overline\chi^{\,2})|
\prod_{i=1}^4|S_{Y_i}(\chi)|\\
&\qquad\ll_\varepsilon B^\varepsilon
\min\left\{
\sqrt{UVP_y}+\frac{\sqrt{UV}\,P_y}{\sqrt m},
UV\sqrt{P_y}\right\}.
\end{aligned}
\end{equation}
The constants are uniform in the six independently chosen
intervals. Conjugating individual character sums or multiplying
their arguments by units modulo \(m\) leaves these bounds unchanged.
\end{lemma}

\begin{proof}
For \(m=1\), the character sum is empty. Assume \(m\ge2\).
The Cochrane--Shi estimate \eqref{ge:CS} gives
\begin{equation}
\frac1{\varphi(m)}\sum_{\psi\ne\chi_0}|S_H(\psi)|^4
\ll_\varepsilon B^\varepsilon H^2.
\label{ar:cs}
\end{equation}
The squaring map on the character group has fibres of size at most
$2^{\omega(m)+1}\ll_\varepsilon B^\varepsilon$. Applying
\eqref{ar:cs} with smaller auxiliary exponents therefore gives
\[
\frac1{\varphi(m)}\sum_{\chi^2\ne\chi_0}|S_H(\chi^2)|^4
\le\frac{2^{\omega(m)+1}}{\varphi(m)}
\sum_{\psi\ne\chi_0}|S_H(\psi)|^4
\ll_\varepsilon B^\varepsilon H^2.
\]

For the four intervals $I_i\subseteq[1,cY_i]$, complete multiplicativity gives
\[
\prod_{i=1}^4S_{Y_i}(\chi)=\sum_{n\le c^4P_y}a_n\chi(n)\quad\text{and}\quad
a_n=\#\{(n_1,\ldots,n_4)\in I_1\times\cdots\times I_4:
n_1n_2n_3n_4=n\}.
\]
Here $0\le a_n\le d_4(n)\ll_\varepsilon B^\varepsilon$ and
$\sum_na_n\ll P_y$. Orthogonality therefore gives
\begin{equation*}
\begin{aligned}
\frac1{\varphi(m)}\sum_{\chi\bmod m}
\left|\prod_{i=1}^4S_{Y_i}(\chi)\right|^2
&=\sum_{\substack{n,n'\le c^4P_y\\n\equiv n'\pmod m\\\gcd(nn',m)=1}}
a_na_{n'}
\ll\max_na_n\left(1+\frac{P_y}{m}\right)\sum_na_n
\\&\ll_\varepsilon B^\varepsilon\left(P_y+\frac{P_y^2}{m}\right).
\end{aligned}
\end{equation*}

Since $|S_H(\overline\chi)|=|S_H(\chi)|$, H\"older's inequality with
exponents $4,4,2$, together with the preceding moment estimates, gives
the first bound in \eqref{ar:nq}:
\[
\begin{aligned}
(B^\varepsilon U^2)^{1/4}(B^\varepsilon V^2)^{1/4}
\left(B^\varepsilon(P_y+P_y^2/m)\right)^{1/2}
\ll_\varepsilon B^\varepsilon
\left(\sqrt{UVP_y}+\frac{\sqrt{UV}\,P_y}{\sqrt m}\right).
\end{aligned}
\]
For the second bound in \eqref{ar:nq}, use
$|S_U(\chi^2)S_V(\overline\chi^{\,2})|\ll UV$
and $\chi^2\ne\chi_0\Rightarrow\chi\ne\chi_0$. Fourfold H\"older
and \eqref{ar:cs} give
\[
UV\prod_{i=1}^4
\left(\frac1{\varphi(m)}\sum_{\chi\ne\chi_0}|S_{Y_i}(\chi)|^4\right)^{1/4}
\ll_\varepsilon B^\varepsilon UV\sqrt{P_y}.
\]
Finally, multiplying an argument by a unit multiplies its character sum
by a factor of absolute value one, and conjugation preserves absolute
values. This proves the stated invariance.
\end{proof}

\begin{lemma}
\label{ar:quadratic}
Let \(M\ge1\), \(Y_i\ge1\), and let
\(I_i\subseteq[1,cY_i]\) be four fixed, independently chosen
integer intervals, where \(c\ge1\) is fixed. Put
\[
S_{Y_i}(\chi)=\sum_{n\in I_i}\chi(n),
\qquad
P_y=Y_1Y_2Y_3Y_4.
\]
Then, for every \(\varepsilon>0\),
\begin{equation}
\sum_{M\le m<2M}
\sum_{\substack{\chi\bmod m\\
\chi^2=\chi_0,\ \chi\ne\chi_0}}
\prod_{i=1}^4|S_{Y_i}(\chi)|
\ll_{\varepsilon,c}
(MP_y)^\varepsilon M P_y^{2/3}.
\label{ar:q-average}
\end{equation}
The inner sum includes all nonprincipal quadratic characters,
whether primitive or imprimitive.
\end{lemma}

\begin{proof}
Write \(S_I(\chi)=\sum_{n\in I}\chi(n)\).
Throughout this proof,
\(\sum_{\chi^*\bmod d}^{*}\) denotes summation over
nonprincipal primitive quadratic characters modulo \(d\),
including those of even conductor.

Fix \(H\ge1\) and an integer interval \(I\subseteq[1,cH]\).
Complete multiplicativity gives
\[
S_I(\chi^*)^2
=
\sum_{n\le c^2H^2}a_n\chi^*(n),
\qquad
a_n=\#\{(r,s)\in I^2:rs=n\}\le\tau(n).
\]
For arbitrary coefficients, Heath-Brown's maximum-coefficient-norm
estimate \cite[Corollary~3]{HB1995} states
\[
\begin{aligned}
\sum_{\substack{\chi^*\ \mathrm{primitive\ real}\\
\operatorname{cond}(\chi^*)\le Q}}
\left|\sum_{n\le N}a_n\chi^*(n)\right|^2
\ll_\varepsilon (QN)^\varepsilon N(Q+N)\max_{n\le N}|a_n|^2.
\end{aligned}
\]
In particular, no squarefree support condition on $a_n$ is required.
Taking $Q=\lceil2D_0\rceil$, $N=\lceil c^2H^2\rceil$, and using
the divisor bound gives
\[
\sum_{D_0\le d<2D_0}\sum_{\chi^*\bmod d}^{*}|S_I(\chi^*)|^4
\ll_{\varepsilon,c}(D_0H)^\varepsilon(D_0H^2+H^4).
\]

There are at most two primitive quadratic characters of any conductor.
P\'olya--Vinogradov \cite[Chapter~23]{Davenport1980} gives
$|S_I(\chi^*)|\ll\sqrt d\log(2d)$, and hence
\begin{equation*}
\begin{aligned}
\sum_{D_0\le d<2D_0}\sum_{\chi^*\bmod d}^{*}|S_I(\chi^*)|^4
\ll_{\varepsilon,c}(D_0H)^\varepsilon
\min\{D_0H^2+H^4,D_0^3\}
\ll_{\varepsilon,c}(D_0H)^\varepsilon D_0H^{8/3}.
\end{aligned}
\end{equation*}
Indeed, use $D_0^3\le D_0H^{8/3}$ when $D_0\le H^{4/3}$,
and $D_0H^2+H^4\le2D_0H^{8/3}$ otherwise.
Dyadic summation, for $T_0\ge1$, gives
\begin{equation}
\begin{aligned}
\sum_{d\le T_0}\frac1d\sum_{\chi^*\bmod d}^{*}|S_I(\chi^*)|^4
&\le\sum_{\substack{j\ge0\\2^j\le T_0}}2^{-j}
\sum_{2^j\le d<2^{j+1}}\sum_{\chi^*\bmod d}^{*}|S_I(\chi^*)|^4\\
&\ll_{\varepsilon,c}H^{8/3+\varepsilon}
\sum_{\substack{j\ge0\\2^j\le T_0}}2^{j\varepsilon}
\ll_{\varepsilon,c}(T_0H)^\varepsilon H^{8/3}.
\end{aligned}
\label{ar:weighted-primitive}
\end{equation}

Every nonprincipal quadratic character modulo $m$ is induced by a unique
primitive quadratic character $\chi^*$ of conductor $d\mid m$
\cite[Chapter~5]{Davenport1980}. Write $m=d\ell$; no condition on
$\gcd(d,\ell)$ is needed, since
$\chi(n)=\chi^*(n)\mathbf1_{\gcd(n,\ell)=1}$.
For $I/a=\{n\in\mathbb Z_{>0}:an\in I\}$, M\"obius inversion gives
\[
S_I(\chi)=\sum_{\substack{a\mid\ell\\\gcd(a,d)=1}}
\mu(a)\chi^*(a)S_{I/a}(\chi^*)\quad\text{and}\quad
|S_I(\chi)|^4\le\tau(m)^3\sum_{a\mid\ell}|S_{I/a}(\chi^*)|^4.
\]
The inducing moduli and nonempty intervals satisfy
\[
\#\{m\in[M,2M):da\mid m\}\le\frac{2M}{da},\qquad
I/a\ne\varnothing\Longrightarrow
a\le cH,\qquad
 I/a\subseteq[1,cH/a]\]
and $\max(1,H/a)\asymp_c H/a.$
For fixed $a$, the interval $I/a$ is independent of $d$.
Apply \eqref{ar:weighted-primitive} with $T_0=2M/a$ and scale
$\max(1,H/a)$, together with $\tau(m)^3\ll_\varepsilon M^{\varepsilon/2}$:
\begin{equation}\label{ar:quadratic-fourth}
\begin{aligned}
&
\sum_{M\le m<2M}\sum_{\substack{\chi\bmod m\\
\chi^2=\chi_0,\ \chi\ne\chi_0}}|S_I(\chi)|^4\\
&\quad\ll_\varepsilon M^{1+\varepsilon/2}
\sum_{a\le\min(cH,2M)}\frac1a
\sum_{d\le2M/a}\frac1d\sum_{\chi^*\bmod d}^{*}|S_{I/a}(\chi^*)|^4\\
&\ll_{\varepsilon,c}M^{1+\varepsilon/2}(MH)^{\varepsilon/2}
H^{8/3}\sum_{a\ge1}a^{-11/3}
\ll_{\varepsilon,c}(MH)^\varepsilon MH^{8/3}.
\end{aligned}
\end{equation}

Applying H\"older's inequality on the common set of pairs \((m,\chi)\),
the fourth-moment bound \eqref{ar:quadratic-fourth} gives
\[
\begin{aligned}
\sum_{M\le m<2M}
\sum_{\substack{\chi\bmod m\\ \chi^2=\chi_0,\ \chi\ne\chi_0}}
\prod_{i=1}^4|S_{Y_i}(\chi)|
\ll_{\varepsilon,c}
\prod_{i=1}^4\bigl((MY_i)^\varepsilon MY_i^{8/3}\bigr)^{1/4}
\ll_{\varepsilon,c}(MP_y)^\varepsilon MP_y^{2/3}.
\end{aligned}
\]
This proves \eqref{ar:q-average}.
\end{proof}

\subsection{Insertion of the complete short sieve and joint weights}
\label{ar:weighted-characters}

Retain the six short-variable scales $\mathbf H$ and $P_y=\prod_{i=1}^4Y_i$.
All moduli and scales under consideration are bounded by fixed
powers of \(B\ge2\).

Let \(\Phi\) be smooth and compactly supported in a fixed positive
box after the substitution
$(a,b,c,d,u,v)=(Z_ax_a,Z_bx_b,Cx_c,Dx_d,Ux_u,Vx_v)$ and
$t_i=Y_ix_i$ with $1\le i\le4.$
Extend \(\Phi\) by zero. The norm \(\|\Phi\|_{j,*}\) is the maximum
of the sup norms of all derivatives of total order at most \(j\)
of this normalized composite function.
For \(\nu\in\{-1,1\}\), put
\begin{equation}
\begin{aligned}
I_\nu(q,\mathbf t)
=
\int_{\mathbb R^2}
\Phi\!\left(
\frac{\nu L_1}{q},\frac{\nu L_2}{q},c,d,\mathbf t
\right)\,dc\,dd\quad\text{and}\quad
J_\nu(q)
=
\int_{\mathbb R_{>0}^6}I_\nu(q,\mathbf t)\,d\mathbf t.
\end{aligned}
\label{ar:integrated-weight}
\end{equation}

The natural scales satisfy
$CUY_2=DVY_1=QZ_a$ and $CVY_3=DUY_4=QZ_b.$
With
\[
c=Cx_c,\quad d=Dx_d,\quad
u=Ux_u,\quad v=Vx_v,\quad t_i=Y_ix_i,
\]
we obtain
\[
\begin{aligned}
\frac{L_1}{qZ_a}
=
\frac Qq
\bigl(\sigma_2x_ux_2x_c+\sigma_1x_vx_1x_d\bigr)\quad\text{and}\quad
\frac{L_2}{qZ_b}
=
-\frac Qq
\bigl(\sigma_3x_vx_3x_c+\sigma_4x_ux_4x_d\bigr).
\end{aligned}
\]
Every positive-order derivative of these polynomials is
\(O(Q/q)\) on the fixed normalized box. Thus, for
\(q\asymp\eta Q\), \(0<\eta\le1\), the chain rule gives
\begin{equation}
\begin{aligned}
\|I_\nu(q,\cdot)\|_{6,*}
&\ll CD\,\eta^{-6}\|\Phi\|_{6,*},\\
\sup_{\mathbf t}
\left\|
\Phi\!\left(
\frac{\nu L_1(C\,\cdot,D\,\cdot)}q,
\frac{\nu L_2(C\,\cdot,D\,\cdot)}q,
C\,\cdot,D\,\cdot,\mathbf t
\right)
\right\|_{C^1}
&\ll\eta^{-1}\|\Phi\|_{1,*}.
\end{aligned}
\label{ar:derivative-loss}
\end{equation}
Indeed, an order-\(j\) chain-rule term contains at most \(j\)
derivatives of the recovery coordinates, so its contribution is
\[
\ll_j
\max(1,Q/q)^j\|\Phi\|_{j,*}
\ll_j\eta^{-j}\|\Phi\|_{j,*}.
\]
The substitution \((c,d)=(Cx_c,Dx_d)\) contributes the factor \(CD\).

Choose a fixed \(b_0\ge1\) such that all normalized short-variable
supports are compactly contained in \((0,b_0)^6\).
The Boolean expansion only uses \(d_j\le b_0H_j\), so
\[
H_j/d_j\ge b_0^{-1},\qquad
\max(1,H_j/d_j)\asymp_{b_0}H_j/d_j,\qquad
\frac{d_jn_j}{H_j}=\frac{n_j}{H_j/d_j}.
\]
Thus the interval estimates remain uniform after division by $d_j$,
without an additional derivative loss. Recall from
Section~\ref{ar:boolean-subsection} that
$|c_p(\{j\})|\ll p^{-1}$ and $|c_p(E)|\ll1$.
The divisor exponents supplied by
Lemmas~\ref{ar:nonquadratic} and~\ref{ar:quadratic} are therefore
\[
\begin{array}{c|c|c}
\text{bound}
&
\boldsymbol\alpha=(\alpha_1,\ldots,\alpha_6)
&
\displaystyle
\sum_{E\subseteq\{1,\ldots,6\}}
|c_p(E)|p^{-\sum_{j\in E}\alpha_j}
\\ \hline
\sqrt{UVP_y}
&
(\tfrac12,\tfrac12,\tfrac12,\tfrac12,\tfrac12,\tfrac12)
&
1+O(p^{-1})
\\
\sqrt{UV}\,P_y/\sqrt m
&
(\tfrac12,\tfrac12,1,1,1,1)
&
1+O(p^{-1})
\\
UV\sqrt{P_y}
&
(1,1,\tfrac12,\tfrac12,\tfrac12,\tfrac12)
&
1+O(p^{-1})
\\
UVMP_y^{2/3}
&
(1,1,\tfrac23,\tfrac23,\tfrac23,\tfrac23)
&
1+O(p^{-4/3})
\end{array}
\]
The singleton contribution is $O(p^{-1-\alpha_j})$, and each larger
support contributes $O(p^{-\sum_{j\in E}\alpha_j})$. Hence
\[
\sum_{\substack{\mathbf d\\d_j\le b_0H_j\ (1\le j\le6)}}
\frac{|c(\mathbf d)|}{\prod_jd_j^{\alpha_j}}
\ll_\varepsilon B^\varepsilon
\]
for the first three rows, while the fourth gives
$\sum_{\mathbf d}\frac{|c(\mathbf d)|}
{d_1d_2(d_3d_4d_5d_6)^{2/3}}<\infty.$

Fix \(\nu\), and define the fully sieved character sum
\[
\mathcal T_{m,\chi}(q)
=
\sum_{\mathbf t\in\mathbb Z_{>0}^6}
\mathbf1_{\mathcal C}(\mathbf t)\rho(\mathbf t)
I_\nu(q,\mathbf t)
\chi(u^2t_2t_4)\overline{\chi(v^2t_1t_3)}.
\]
Write
$(\chi_1,\ldots,\chi_6)
=(\chi^2,\overline\chi^{\,2},\overline\chi,\chi,\overline\chi,\chi)$ and
$F_q(\mathbf x)=I_\nu(q,\mathbf H\mathbf x),$
where vector products are coordinatewise.
Since characters vanish on nonunits,
\[
\prod_{j=1}^6\chi_j(d_j)=0
\quad\text{if }\gcd(d_1\cdots d_6,m)>1,
\qquad
\left|\prod_{j=1}^6\chi_j(d_j)\right|=1
\quad\text{otherwise}.
\]

The weight depends on $q$, so we first express it using prefixes common
to all moduli. Since $F_q$ vanishes near the upper faces of $[0,b_0]^6$,
repeated partial summation \cite[Chapter~I, Theorem~1.1]{SS}, or the identity
$
F_q(\mathbf z)=\int_{z_1}^{b_0}\cdots\int_{z_6}^{b_0}
\partial_1\cdots\partial_6F_q(\mathbf x)\,d\mathbf x,$
gives
\begin{equation}
\begin{aligned}
\mathcal T_{m,\chi}(q)
={}&
\frac1{\zeta(2)}
\sum_{\substack{d_j\le b_0H_j\ (1\le j\le6)\\
\gcd(d_1\cdots d_6,m)=1}}
c(\mathbf d)\prod_{j=1}^6\chi_j(d_j)\\
&\times
\int_{[0,b_0]^6}
\partial_1\cdots\partial_6F_q(\mathbf x)
\prod_{j=1}^6
\left(
\sum_{1\le n\le H_jx_j/d_j}\chi_j(n)
\right)\,d\mathbf x.
\end{aligned}
\label{ar:common-prefix}
\end{equation}

For each fixed \(\mathbf d,\mathbf x\), these six prefix intervals
are independent of \(m,\chi,q\).
Moreover, \eqref{ar:derivative-loss} gives
\[
\sup_{q\asymp\eta Q}
\|\partial_1\cdots\partial_6F_q\|_\infty
\ll CD\,\eta^{-6}\|\Phi\|_{6,*}.
\]

For fixed \(m,q\), apply Lemma~\ref{ar:nonquadratic} inside
\eqref{ar:common-prefix}, and then sum the Boolean coefficients
using the first three rows of the table. This yields both bounds
\begin{equation}\label{ar:weighted-nq}
\begin{aligned}
\frac1{\varphi(m)}
\sum_{\substack{\chi\bmod m\\\chi^2\ne\chi_0}}
|\mathcal T_{m,\chi}(q)|
\ll_\varepsilon{}&
B^\varepsilon CD\,\eta^{-6}\|\Phi\|_{6,*}\min\!\left\{
\sqrt{UVP_y}+\frac{\sqrt{UV}\,P_y}{\sqrt m},
\ UV\sqrt{P_y}
\right\}.
\end{aligned}
\end{equation}

For quadratic nonprincipal characters,
\[
\left|
\sum_{1\le n\le H_jx_j/d_j}\chi_j(n)
\right|
\ll\frac{H_j}{d_j}
\qquad(j=1,2).
\]
Lemma~\ref{ar:quadratic}, applied to the remaining four fixed
prefix intervals, consequently gives
\[
\begin{aligned}
\sum_{M\le m<2M}
\sum_{\substack{\chi\bmod m\\
\chi^2=\chi_0,\ \chi\ne\chi_0}}
\prod_{j=1}^6
\left|
\sum_{1\le n\le H_jx_j/d_j}\chi_j(n)
\right|\ll_\varepsilon
B^\varepsilon
\frac{UVMP_y^{2/3}}
{d_1d_2(d_3d_4d_5d_6)^{2/3}},
\end{aligned}
\]
uniformly for $\mathbf x\in[0,b_0]^6$.

Let \(b_{m,q}\ge0\) be a finitely supported family, with
\(b_{m,q}=0\) unless \(q\asymp\eta Q\).
Taking absolute values in \eqref{ar:common-prefix}, applying
the preceding estimate at each fixed \(\mathbf x,\mathbf d\),
and using the convergent Boolean sum gives
\begin{equation}\label{ar:weighted-q}
\begin{aligned}
\sum_{M\le m<2M}
\sum_{\substack{\chi\bmod m\\
\chi^2=\chi_0,\ \chi\ne\chi_0}}
\sum_q b_{m,q}|\mathcal T_{m,\chi}(q)|
\ll_\varepsilon
B^\varepsilon CD\,\eta^{-6}\|\Phi\|_{6,*}
\,UVMP_y^{2/3}
\sup_{M\le m<2M}\sum_q b_{m,q}.
\end{aligned}
\end{equation}
In particular,
\[
b_{m,q}=0\quad(q\nmid m)
\quad\Longrightarrow\quad
\sum_q b_{m,q}
\le\tau(m)\sup_q b_{m,q}.
\]
This coefficient bound accommodates the divisor-supported weights used
in Theorem~\ref{ar:main}, while all prefix intervals remain common to
the moduli.

All implied constants depend only on the fixed normalized
support, the fixed polynomial bounds, the comparison constants
in \(q\asymp\eta Q\), and \(\varepsilon\).

\subsection{A uniform main-term transfer}
\label{ar:transfer-subsection}

Use the common scale notation
\[
 A=Z_aZ_b,\quad H=CD,\quad J=UV,\quad P_0=AHJ,\quad
 T=(BP_0)^{1/3},\quad Q=T/A,\quad Q'=T/H,
\]
and \(Y_i=T/A_i(\mathbf Z)\).  In particular,
\begin{equation*}
 U^2Y_2Y_4=V^2Y_1Y_3=QQ',\qquad
 CD\,UVP_y=BQ,\qquad P_y=\prod_iY_i.
\end{equation*}
The normalized coefficients of \(L_1/(QZ_a)\) and \(L_2/(QZ_b)\)
and all their fixed-order derivatives are bounded on the fixed box.
Assume $H_{\min},Q,Q'\ge1$, that all scales have polynomial size in $B$,
and that every integer short tuple in the support satisfying $\mathcal C$
has $\Delta\ne0$.

Fix \(\eta_{\max}\ge1\), let \(0<\eta\le\eta_{\max}\), put
\(L=\eta Q\ge1\), and let \(\psi\) be smooth and supported in \([1,2]\).
Define the sum of the geometric main terms
\begin{align}
 \mathcal M_\Phi={}&
 \sum_{\mathbf t\in\mathbb Z_{>0}^6}
      \mathbf 1_{\mathcal C}(\mathbf t)\rho(\mathbf t)
 \sum_{q\mid|\Delta|}\psi(q/L)
       \sum_{\nu=\pm1}I_\nu(q,\mathbf t)
       \sum_{k\mid|\Delta|/q}\frac{\mu(k)w(qk)}{qk}.
 \label{ar:geometric-main-sum}
\end{align}
Proposition~\ref{ar:long} controls the difference between this sum and
the original integer count.

\begin{theorem}
\label{ar:main}
Retain the hypotheses above and the definition
\eqref{ar:geometric-main-sum} of \(\mathcal M_\Phi\). Put
\begin{equation}
\mathcal R=
\min\left\{
\frac Q{\sqrt{P_y}},
\frac Q{\sqrt{JP_y}}+\sqrt{\frac QJ}
\right\}
+P_y^{-1/3}+H_{\min}^{-1/2}+Q^{-1}+{Q'}^{-1}.
\label{ar:relative-error}
\end{equation}
For every \(\varepsilon>0\),
\begin{equation}
\begin{aligned}
\mathcal M_\Phi
={}
\mathfrak S\sum_{\nu=\pm1}
\int_0^\infty
\psi(q/L)J_\nu(q)\frac{dq}{q^2}
+
O_\varepsilon\!\left(
B^{1+\varepsilon}\eta^{-7}\mathcal R
\|\Phi\|_{6,*}\|\psi\|_{C^1}
\right),
\end{aligned}
\label{ar:main-formula}
\end{equation}
where
\[
\mathfrak S=
\prod_p(1-p^{-1})^7(1+7p^{-1}+p^{-2}).
\]
The implied constant depends only on \(\varepsilon\),
\(\eta_{\max}\), the fixed normalized support, and the fixed
polynomial bounds.
\end{theorem}

\begin{proof}
Since $\max(1,\eta^{-1})\le\eta_{\max}\eta^{-1},$
the derivative estimates \eqref{ar:derivative-loss} remain valid
throughout \(0<\eta\le\eta_{\max}\), with constants depending on
\(\eta_{\max}\).

Choose a fixed $c_3\ge2\eta_{\max}$ such that $|\Delta|\le c_3QQ'$
on the support, and set $K_0=\lfloor c_3Q'/\eta\rfloor\asymp Q'/\eta$.
Then
\[
\Delta\ne0,\quad L\le q\le2L,\quad qk\mid\Delta
\Longrightarrow k\le K_0,
\qquad qk\le2c_3QQ'.
\]
Thus the cutoff is independent of the short tuple, and
\begin{equation*}
\sum_{\substack{k\ge1\\qk\mid\Delta}}\frac{\mu(k)w(qk)}{qk}
=\sum_{k\le K_0}\frac{\mu(k)w(qk)}{qk}\mathbf1_{qk\mid\Delta}.
\end{equation*}
The outer $q$-sum may now run over all positive integers.
Since $\eta^{-1}\le Q$, all moduli and divisor ranges remain
polynomially bounded in $B$, uniformly in $\eta$.

By \eqref{ar:units} and character orthogonality,
\begin{equation*}
\begin{aligned}
\mathbf1_{\mathcal C}(\mathbf t)\rho(\mathbf t)
\mathbf1_{m\mid\Delta}
={}
\frac{\mathbf1_{\mathcal C}(\mathbf t)\rho(\mathbf t)}
{\varphi(m)}
\sum_{\chi\bmod m}
\chi(\sigma_1\sigma_2\sigma_3\sigma_4)
\chi(u^2t_2t_4)\overline{\chi(v^2t_1t_3)}.
\end{aligned}
\end{equation*}
If a short variable is not a unit modulo \(m\), both sides vanish.
We use the divisor-sum formula \eqref{ar:boolean} to expand
\(\mathbf1_{\mathcal C}(\mathbf t)\rho(\mathbf t)\)
in this identity, keeping all other weight factors
in each summand.

For the principal character, use Lemma~\ref{ar:principal} and
\eqref{ar:derivative-loss}. Since
\[
\begin{aligned}
\sum_{q\asymp L}\sum_{k\le K_0}
\frac1{qk\varphi(qk)}
&\ll_\varepsilon
B^\varepsilon
\sum_{q\asymp L}q^{-2}\sum_{k\ge1}k^{-2}
\ll_\varepsilon B^\varepsilon L^{-1},\\
CDJP_y\,\eta^{-6}L^{-1}
&=BQ\eta^{-6}(\eta Q)^{-1}
=B\eta^{-7},
\end{aligned}
\]
the total principal-character error is
\begin{equation}
O_\varepsilon\!\left(
B^{1+\varepsilon}\eta^{-7}H_{\min}^{-1/2}
\|\Phi\|_{6,*}\|\psi\|_\infty
\right).
\label{ar:principal-error}
\end{equation}

For the nonquadratic characters, apply \eqref{ar:weighted-nq}
with the full short sieve and the joint weight.
The required divisor sums satisfy
\[
\sum_{q\asymp L}\sum_{k\le K_0}(qk)^{-1}
\ll\log(2B),
\qquad
\sum_{q\asymp L}\sum_{k\le K_0}(qk)^{-3/2}
\ll L^{-1/2}.
\]
Using \(CDJP_y=BQ\), we obtain
\[
CD\eta^{-6}J\sqrt{P_y}
=B\eta^{-6}\frac{Q}{\sqrt{P_y}},\quad
CD\eta^{-6}\sqrt{JP_y}
=B\eta^{-6}\frac{Q}{\sqrt{JP_y}},\]
and $CD\eta^{-6}\sqrt J\,P_yL^{-1/2}=B\eta^{-13/2}\sqrt{\frac QJ}.$
Thus the nonquadratic contribution is
\begin{equation}
O_\varepsilon\!\left(
B^{1+\varepsilon}\eta^{-6}
\min\left\{
\frac Q{\sqrt{P_y}},
\frac Q{\sqrt{JP_y}}+\eta^{-1/2}\sqrt{\frac QJ}
\right\}
\|\Phi\|_{6,*}\|\psi\|_\infty
\right).
\label{ar:nq-error}
\end{equation}

For the quadratic nonprincipal characters, take
\[
b_{m,q}=
\begin{cases}
\displaystyle
\frac{|\mu(m/q)|w(m)}{m\varphi(m)}
|\psi(q/L)|,
&q\mid m,\quad m/q\le K_0,\\[2mm]
0,&\text{otherwise}.
\end{cases}
\]
Then
\[
\sum_qb_{m,q}
\le
\frac{\tau(m)}{m\varphi(m)}\|\psi\|_\infty
\ll_\varepsilon
B^\varepsilon m^{-2}\|\psi\|_\infty.
\]
The estimate \eqref{ar:weighted-q}
bounds the contribution of \(M\le m<2M\) by
\[
\ll_\varepsilon
B^\varepsilon CDJ\eta^{-6}
M^{-1}P_y^{2/3}
\|\Phi\|_{6,*}\|\psi\|_\infty.
\]
Every contributing modulus satisfies \(m\ge L\). Hence
\[
\sum_{\substack{M\ {\rm dyadic}\\M\ge L/2}}M^{-1}\ll L^{-1},
\qquad
\frac{CDJP_y^{2/3}\eta^{-6}}L
=B\eta^{-7}P_y^{-1/3},
\]
and the quadratic contribution is
\begin{equation}
O_\varepsilon\!\left(
B^{1+\varepsilon}\eta^{-7}P_y^{-1/3}
\|\Phi\|_{6,*}\|\psi\|_\infty
\right).
\label{ar:quadratic-error}
\end{equation}

The remaining principal expression is
\[
\sum_{\nu=\pm1}\sum_{q\ge1}
\psi(q/L)J_\nu(q)
\sum_{k\le K_0}
\frac{\mu(k)w(qk)\mathfrak d(qk)}{qk\varphi(qk)}.
\]
Direct substitution from \eqref{ar:dp} gives
\begin{equation*}
\begin{aligned}
\frac{w(m)\mathfrak d(m)}{m\varphi(m)}
={}
\frac1{m^2}
\prod_{p\nmid m}
(1-p^{-1})^6(1+6p^{-1}+p^{-2})
\prod_{p\mid m}(1-p^{-1})^6
\ll m^{-2}.
\end{aligned}
\end{equation*}
The infinite product converges absolutely.

Since
$
|J_\nu(q)|\ll CDJP_y\|\Phi\|_{0,*}
=BQ\|\Phi\|_{0,*},$
extending this principal \(k\)-series to infinity costs
\begin{equation}
O\!\left(
\frac{BQ}{LK_0}\|\Phi\|_{0,*}\|\psi\|_\infty
\right)
=
O\!\left(
\frac B{Q'}\|\Phi\|_{0,*}\|\psi\|_\infty
\right).
\label{ar:principal-k-tail}
\end{equation}

Define the multiplicative function \(f\) by
\[
f(1)=1,\qquad
f(p^e)=\frac{1-p^{-2}}{1+6p^{-1}}
\quad(e\ge1).
\]
After extracting \(q^{-2}\), the local factors of the
infinite \(k\)-series, with \(r=p^{-1}\), are
\[
\begin{array}{c|c}
p\nmid q
&
(1-r)^6\bigl((1+6r+r^2)-r^2\bigr)
=(1-r)^6(1+6r)
\\[1mm]
p\mid q
&
(1-r)^6(1-r^2).
\end{array}
\]
Therefore
\begin{equation*}
\sum_{k\ge1}
\frac{\mu(k)w(qk)\mathfrak d(qk)}{qk\varphi(qk)}
=
\frac{f(q)}{q^2}
\prod_p(1-p^{-1})^6(1+6p^{-1}).
\end{equation*}

Put $g=\mu*f$. Then
\[
g(p)=-\frac{6p^{-1}+p^{-2}}{1+6p^{-1}},\qquad
g(p^e)=0\ (e\ge2),\qquad |g(p)|\le6/p.
\]
Mertens' estimate and partial summation give
\[
\sum_{n\le X}|g(n)|\le\prod_{p\le X}(1+6/p)\ll\log^6(2X)\quad\text{and}\quad
\sum_{n>X}\frac{|g(n)|}{n}
\ll\int_X^\infty\frac{\log^6(2t)}{t^2}\,dt
\ll\frac{\log^6(2X)}X.
\]
Since $f=1*g$,
\[
\begin{aligned}
\sum_{q\le X}f(q)
=\sum_{n\le X}g(n)\left\lfloor\frac Xn\right\rfloor
=X\sum_{n\ge1}\frac{g(n)}n
+O\left(\sum_{n\le X}|g(n)|
+X\sum_{n>X}\frac{|g(n)|}n\right).
\end{aligned}
\]
Consequently,
\begin{equation}
\sum_{q\le X}f(q)=A_fX+O(\log^6(2X)),\qquad
A_f=\prod_p\left(1-p^{-1}+p^{-1}f(p)\right).
\label{ar:f-mean}
\end{equation}

Set \(G_\nu(q)=\psi(q/L)J_\nu(q)/q^2\).
With \(c,d,\mathbf t\) fixed and
\(a=\nu L_1/q\), \(b=\nu L_2/q\),
differentiation with respect to \(q\) gives
\[
q\frac{d}{dq}\Phi(\nu L_1/q,\nu L_2/q,c,d,\mathbf t)
=-(a\partial_a+b\partial_b)\Phi(a,b,c,d,\mathbf t).
\]
Consequently,
\[
|J_\nu(q)|+|qJ_\nu'(q)|\ll BQ\|\Phi\|_{1,*}.
\]
The normalized $a,b$ coordinates are bounded on the support, and
$q\asymp L$ there; hence
\[
\|G_\nu\|_\infty+\int_0^\infty|G_\nu'(q)|\,dq
\ll BQL^{-2}\|\Phi\|_{1,*}\|\psi\|_{C^1}.
\]
Partial summation \cite[Chapter~I, Theorem~1.1]{SS} applied to
\eqref{ar:f-mean} gives
\begin{equation}
\begin{aligned}
\sum_{q\ge1}f(q)G_\nu(q)
=A_f\int_0^\infty G_\nu(q)\,dq
+O\left(\frac{B\log^6(2B)}{\eta^2Q}
\|\Phi\|_{1,*}\|\psi\|_{C^1}\right).
\end{aligned}
\label{ar:q-mean-error}
\end{equation}

Finally,
\[
1-r+rf(p)
=
\frac{(1-r)(1+7r+r^2)}{1+6r},
\qquad r=p^{-1},
\]
so the main coefficient is exactly
\[
A_f\prod_p(1-p^{-1})^6(1+6p^{-1})
=
\prod_p(1-p^{-1})^7(1+7p^{-1}+p^{-2})
=
\mathfrak S.
\]

For $0<\eta\le\eta_{\max}$, all displayed powers of $\eta^{-1}$
are $O_{\eta_{\max}}(\eta^{-7})$. After extracting
$B^{1+\varepsilon}\eta^{-7}\|\Phi\|_{6,*}\|\psi\|_{C^1}$,
the five errors contribute
\[
\begin{array}{c|c}
\eqref{ar:principal-error} & H_{\min}^{-1/2}\\[1mm]
\eqref{ar:nq-error}
&\min\{Q/\sqrt{P_y},\ Q/\sqrt{JP_y}+\sqrt{Q/J}\}\\[1mm]
\eqref{ar:quadratic-error} & P_y^{-1/3}\\[1mm]
\eqref{ar:principal-k-tail} & {Q'}^{-1}\\[1mm]
\eqref{ar:q-mean-error} & Q^{-1}.
\end{array}
\]
Their sum is \eqref{ar:relative-error}, proving \eqref{ar:main-formula}.
As before, smaller auxiliary exponents absorb all divisor and logarithmic
losses into the final $B^\varepsilon$.
\end{proof}

\subsection{The continuous measure and the actual count}
\label{ar:measure-subsection}

All real integrals below are over the fixed sign sector
\(y_i=\sigma_i t_i\), \(t_i>0\), and positive variables \(z_e\).
The weight \(\Phi\) is extended by zero outside its positive support.

Put $\mathcal F=bL_1-aL_2=\sum_iA_i(z)y_i$.
Since $\partial\mathcal F/\partial y_4=A_4(z)=adu>0$, the hypersurface
is regular and $\delta_0(\mathcal F)\,d^4y\,d^6z$ is its Leray measure.
On the support of $\psi(|L_1|/(aL))$, one has $|L_1|\ge aL>0$;
therefore
\[
\mathcal F=0,\;a,b>0
\Longrightarrow L_2=bL_1/a,\;
\nu=\operatorname{sgn}L_1=\operatorname{sgn}L_2\quad\text{and}\quad
q=|L_1|/a=|L_2|/b>0.
\]
Integrating first in $b$ and then substituting for $a$ gives
\[
\delta_0(bL_1-aL_2)\,db
=\frac1{|L_1|}\delta_0(b-aL_2/L_1)\,db,\quad
a=\frac{\nu L_1}{q},\quad b=\frac{\nu L_2}{q},\quad\text{and}\quad
\frac{|da|}{|L_1|}=\frac{dq}{q^2}.
\]
The sign change $y_i=\sigma_it_i$ has absolute Jacobian one, so
\eqref{ar:integrated-weight} yields
\begin{equation}
\begin{aligned}
\sum_{\nu=\pm1}\int_0^\infty\psi(q/L)J_\nu(q)\frac{dq}{q^2}
=\int\Phi\,\psi\left(\frac{|L_1|}{aL}\right)
\delta_0\left(\sum_iA_i(z)y_i\right)\,d^4y\,d^6z.
\end{aligned}
\label{ar:coarea}
\end{equation}
On $\mathcal F=0$, $|L_1|/(aL)=|L_2|/(bL)$.
Each contributing point determines exactly one $\nu$, so no factor
$1/2$ occurs. The locus $L_1=L_2=0$ contributes nothing because $\psi(0)=0$.

\begin{corollary}
\label{ar:count}
Assume the hypotheses of Theorem~\ref{ar:main}.
Let \(\mathcal N_\Phi\) be the weighted sum over integer torsor
tuples in the fixed sign sector satisfying all canonical
coprimality conditions, with weight
$\Phi\,\psi(q/L)$ for $q=\gcd(L_1,L_2)>0.$
Then
\begin{align*}
\mathcal N_\Phi
={}&
\mathfrak S
\int
\Phi\,
\psi\!\left(\frac{|L_1|}{aL}\right)
\delta_0\!\left(\sum_{i=1}^4A_i(z)y_i\right)
\,d^4y\,d^6z+
O_\varepsilon\!\left(
B^{1+\varepsilon}\eta^{-7}\mathcal R
\|\Phi\|_{6,*}\|\psi\|_{C^1}
\right)\\
&+
O_\varepsilon\!\left(
B^\varepsilon\eta^{-1}\|\Phi\|_{1,*}\|\psi\|_\infty
\sum_{\substack{\mathbf t\text{ in the short box}\\
\mathcal C,\ \Delta\ne0}}
\ \sum_{\substack{q\mid|\Delta|\\L\le q\le2L}}
\lambda_q^{-1}
\right).
\end{align*}
The constants have the same dependencies as in Theorem~\ref{ar:main}.
The geometric-error sum may be restricted to any set containing every
$(\mathbf t,q)$ for which $\psi(q/L)W_{\nu,q,\mathbf t}\not\equiv0$
for at least one $\nu$, where
$W_{\nu,q,\mathbf t}(c,d)=\Phi(\nu L_1/q,\nu L_2/q,c,d,\mathbf t)$.
This concerns the recovered real support, including blocks containing
no admissible integer point.
\end{corollary}

\begin{proof}
Use \(W_{\nu,q,\mathbf t}\) as defined above. The bijection in
\eqref{ar:recover} gives
\[
\mathcal N_\Phi
=\sum_{\substack{\mathbf t\text{ in the short box}\\ \mathcal C,\ \Delta\ne0}}
\sum_{q\mid|\Delta|}\psi(q/L)
\sum_{\nu=\pm1}N_{W_{\nu,q,\mathbf t}}(q).
\]
By \eqref{ar:integrated-weight} and \eqref{ar:derivative-loss},
uniformly for \(L\le q\le2L\), we have
\(\int W_{\nu,q,\mathbf t}=I_\nu(q,\mathbf t)\) and
\(\|W_{\nu,q,\mathbf t}(C\,\cdot,D\,\cdot)\|_{C^1}
\ll\eta^{-1}\|\Phi\|_{1,*}\).
Summing Proposition~\ref{ar:long} therefore yields
\[
\begin{aligned}
\mathcal N_\Phi
=\mathcal M_\Phi+
O_\varepsilon\!\Bigg(
B^\varepsilon\eta^{-1}\|\Phi\|_{1,*}\|\psi\|_\infty
\sum_{\substack{\mathbf t\text{ in the short box}\\
\mathcal C,\ \Delta\ne0}}
\ \sum_{\substack{q\mid|\Delta|\\L\le q\le2L}}
\lambda_q^{-1}
\Bigg).
\end{aligned}
\]
Theorem~\ref{ar:main} and \eqref{ar:coarea} now give the stated
formula.

\end{proof}

By \eqref{ar:units}, every pair \((\mathbf t,q)\) in the
geometric-error sum of Corollary~\ref{ar:count} satisfies
\eqref{ge:q-conditions}.
Replacing \(L_2\) by \(-L_2\) leaves \(\Lambda_q\) and its
normalized shortest length \(\lambda_q\) unchanged.
Thus, whenever their scale hypotheses hold,
Propositions~\ref{ge:four-proposition} and~\ref{ge:reuss-proposition}
apply to the full geometric-error sum, including tuples whose recovered
real support contains no admissible integer point. The required scale
conditions are verified block by block in
Section~\ref{gl:geometric-core}.

\begin{remark}
\label{ar:other-lines}
Corollary~\ref{ar:count} excludes the line $12|34$, since
$A_1y_1+A_2y_2=bL_1=\nu abq\ne0.$
For any pairing $ij|kl$, an integer torsor tuple satisfying
\eqref{in:coprimality} obeys
\[
A_iy_i+A_jy_j=0
\Longrightarrow z_{jk}z_{jl}y_i=-z_{ik}z_{il}y_j\;\text{and}\;
\gcd(z_{jk}z_{jl},z_{ik}z_{il})=1
\Longrightarrow z_{ik}z_{il}\mid y_i,\quad z_{jk}z_{jl}\mid y_j.
\]
The conditions \(\gcd(y_r y_s,z_{rs})=1\) for \(r<s\) give
$\gcd(y_i,z_{ik}z_{il})=\gcd(y_j,z_{jk}z_{jl})=1$; hence
\[
A_iy_i+A_jy_j=0\Longrightarrow z_{ik}=z_{il}=z_{jk}=z_{jl}=1.
\]
On a fixed positive normalized support in the strict core,
$z_e\ge c_0Z_e\ge c_0B^\delta>1$ for sufficiently large $B$.
Thus the other two pair-sum exclusions hold automatically there.
Outside this core they must be imposed separately when applying
Corollary~\ref{ar:count} to the count on $U$.
\end{remark}

\section{The common main term in the Poisson range}
\label{po:section}

Throughout this section the six positive variables \(z_e\) remain pairwise
coprime. We first relax the conditions involving the $y$-variables and
then impose them at a fixed finite set of primes. The resulting
Poisson modulus remains fixed as $B\to\infty$.

\subsection{Normalization after relaxing all y-conditions}

Retain the torsor equation and coordinate map in
\eqref{in:torsor}, and use the unreduced height
$H_{\mathrm{raw}}$ from \eqref{in:raw-height}.

Weak tuples are as in Definition~\ref{ta:weak-definition},
except that both overall signs are allowed.
A tuple is \emph{y-relaxed} if the coprimality conditions
$\gcd(y_i,z_{ij})=1$ with $i\ne j$
are also omitted. In both cases we retain
\[
z_e\in\mathbb Z_{>0},\qquad
y_i\in\mathbb Z\setminus\{0\},\qquad
\gcd(z_e,z_f)=1\quad(e\ne f),
\]
and
\[
A_1(z)y_1+A_j(z)y_j\ne0
\qquad(j=2,3,4).
\]

In this section, $\mathcal N_{\mathrm{weak}}(X)$ counts weak
tuples with $H_{\mathrm{raw}}\le X$, and
$\mathcal T_{\mathrm{weak},>R}(X)$ counts those for which
\[
p\mid\gcd(y_i,y_j)
\quad\text{for some prime }p>R
\quad\text{and some }1\le i<j\le4.
\]
Both counts include $(z,y)$ and $(z,-y)$, so they are twice
the corresponding counts in Section~\ref{ta:section}.
Corollary~\ref{ta:weak-count} and
Theorem~\ref{ta:weak-tail} therefore give
\begin{equation}
\mathcal N_{\mathrm{weak}}(X)
\ll X\log^6(2X)\quad\text{and}\quad
\mathcal T_{\mathrm{weak},>R}(X)
\ll X\log^6(2X)R^{-5/288}.
\label{po:weak-inputs}
\end{equation}
These estimates hold uniformly for $X\ge1$ and $R\ge2$,
with the range $1\le X<2$ following by monotonicity.
The implied constants are independent of $X$ and $R$.

\begin{lemma}
\label{po:incidence-normalization}
For a y-relaxed tuple $(z,y)$, put
\[
d_{ij}=\gcd(z_{ij},|y_i|,|y_j|),\qquad
D_0=\prod_e d_e,\qquad
z'_e=\frac{z_e}{d_e},\qquad
y'_i=\frac{y_i}{C_i(d)}.
\]
Then $(z',y')$ is weak, and
\begin{equation}
A_i(z)y_i=D_0A_i(z')y'_i,\quad
x_{\mathrm{raw},i}(z,y)
=D_0^2x_{\mathrm{raw},i}(z',y').
\label{po:normalization-identities}
\end{equation}
In particular,
\[
H_{\mathrm{raw}}(z,y)
=D_0^2H_{\mathrm{raw}}(z',y').
\]
For fixed $(z',y')$ and $D_0$, there are at most $d_6(D_0)$
preimages. Moreover, for every prime $p$,
\[
p\mid z_{ij},\qquad p\mid y_i y_j
\quad\Longrightarrow\quad p\mid D_0.
\]
\end{lemma}

\begin{proof}
Since $d_e\mid z_e$ and the variables \(z_e\) are pairwise coprime,
\[
\gcd(d_e,d_f)=1\quad(e\ne f),\qquad
C_i(d)=\prod_{j\ne i}d_{ij}\mid y_i.
\]
Thus $z'_e$ are positive integers, $y'_i$ are nonzero integers,
and $\gcd(z'_e,z'_f)=1$ for $e\ne f$.

The identities $A_i(d)C_i(d)=D_0$ and $\prod_iC_i(d)=D_0^2$ give
\[
A_i(z)y_i=D_0A_i(z')y'_i,\qquad
x_{\mathrm{raw},i}(z,y)=D_0^2x_{\mathrm{raw},i}(z',y').
\]
Thus the torsor equation, the three pair-sum exclusions and the
claimed height identity are preserved.

For every index pair $ij$, the definitions imply
\[
\gcd(z'_{ij},|y'_i|,|y'_j|)
\mid
\gcd\left(
\frac{z_{ij}}{d_{ij}},
\frac{|y_i|}{d_{ij}},
\frac{|y_j|}{d_{ij}}
\right)=1.
\]
Let $p\mid z'_{ij}$ be prime. Pairwise coprimality of the new variables \(z'_e\) gives
\[
p\nmid A_i(z')A_j(z'),\qquad
p\mid A_k(z')\quad(k\notin\{i,j\}).
\]
Reducing the new torsor equation modulo $p$, we obtain
\[
A_i(z')y'_i+A_j(z')y'_j\equiv0\pmod p,
\qquad
p\mid y'_i\Longleftrightarrow p\mid y'_j.
\]
The preceding gcd identity excludes simultaneous divisibility.
Hence
\[
\gcd(z'_{ij},y'_i)=\gcd(z'_{ij},y'_j)=1,
\]
and $(z',y')$ is weak.

The same reduction applied to the original tuple gives
\[
p\mid z_{ij}
\quad\Longrightarrow\quad
\bigl(p\mid y_i\Longleftrightarrow p\mid y_j\bigr).
\]
Therefore
\[
p\mid z_{ij},\quad p\mid y_i y_j
\quad\Longrightarrow\quad
p\mid\gcd(z_{ij},|y_i|,|y_j|)
=d_{ij}\mid D_0.
\]

For fixed $(z',y')$ and $d$, the inverse formulas
$z_e=d_ez'_e$ and $y_i=C_i(d)y'_i$ determine at most one preimage;
there are $d_6(D_0)$ ordered factorizations $\prod_e d_e=D_0$.
\end{proof}

\begin{corollary}
\label{po:free-tail}
Let $\mathcal N_{\mathrm{free}}(B)$ count y-relaxed tuples with
$H_{\mathrm{raw}}\le B$, and let
$\mathcal T_{\mathrm{free},>R}(B)$ count those for which, for some
prime $p>R$ and some $i<j$,
\[
p\mid\gcd(y_i,y_j)
\qquad\text{or}\qquad
p\mid\gcd(z_{ij},y_i y_j).
\]
Uniformly for $B,R\ge2$,
\begin{equation}
\begin{aligned}
\mathcal N_{\mathrm{free}}(B)
\ll B\log^6(2B),\quad\text{and}\quad
\mathcal T_{\mathrm{free},>R}(B)
\ll B\log^6(2B)R^{-5/288}.
\end{aligned}
\label{po:free-tail-bound}
\end{equation}
\end{corollary}

\begin{proof}
Use Lemma~\ref{po:incidence-normalization}. Since the normalized
coordinates are nonzero integers,
\[
1\le H_{\mathrm{raw}}(z',y')
=D_0^{-2}H_{\mathrm{raw}}(z,y)\le B/D_0^2,
\qquad D_0\le\sqrt B.
\]
The multiplicity bound and \eqref{po:weak-inputs} give
\[
\mathcal N_{\mathrm{free}}(B)
\le\sum_{D_0\le\sqrt B}d_6(D_0)
\mathcal N_{\mathrm{weak}}(B/D_0^2)
\ll B\log^6(2B)\sum_{D_0\ge1}\frac{d_6(D_0)}{D_0^2}
=\zeta(2)^6B\log^6(2B).
\]

If $p\nmid D_0$, then $y_i=C_i(d)y'_i$ and $C_i(d)\mid D_0$
give $p\mid y_i\Longleftrightarrow p\mid y'_i$, and the normalization
lemma excludes simultaneous divisibility
$p\mid z_{ij}$ and $p\mid y_i y_j$ for any $i<j$.
Thus a forbidden prime
$p>R$ either divides $D_0$ or remains a forbidden prime for the
normalized weak tuple. Consequently,
\[
\mathcal T_{\mathrm{free},>R}(B)
\le\sum_{D_0\le\sqrt B}d_6(D_0)\Bigg\{
\mathcal T_{\mathrm{weak},>R}(B/D_0^2)
+\mathcal N_{\mathrm{weak}}(B/D_0^2)
\sum_{\substack{p>R\\p\mid D_0}}1\Bigg\}.
\]
Here and below $p$ ranges over primes. The Euler product gives
\[
\sum_{\substack{D_0\ge1\\p\mid D_0}}
\frac{d_6(D_0)}{D_0^2}
=\zeta(2)^6\bigl[1-(1-p^{-2})^6\bigr]
\le6\zeta(2)^6p^{-2}.
\]
Applying \eqref{po:weak-inputs}, using
$\log(2B/D_0^2)\le\log(2B)$, and rearranging nonnegative sums yield
\[
\begin{aligned}
\mathcal T_{\mathrm{free},>R}(B)
&\ll B\log^6(2B)\zeta(2)^6
\left(R^{-5/288}+\sum_{p>R}\bigl[1-(1-p^{-2})^6\bigr]\right)\\
&\ll B\log^6(2B)\left(R^{-5/288}+R^{-1}\right)
\ll B\log^6(2B)R^{-5/288}.
\end{aligned}
\]
All implied constants are independent of $B$ and $R$.
\end{proof}

\subsection{Size conditions for Poisson summation
and a common real measure}

Retain $P_0,T,Y_i$ from
\eqref{in:scales}--\eqref{in:short-scales}, and set
\[
(R_1,R_2,R_3)
=(Z_{12}Z_{34},Z_{14}Z_{23},Z_{13}Z_{24})=(A,H,J),
\qquad Q_j=T/R_j\quad(1\le j\le3).
\]
We write $z\asymp\mathbf Z$ for the block
$Z_e\le z_e<2Z_e$ for every index pair $e$.
For a fixed $\delta>0$, assume
\begin{equation}
\min_e Z_e\ge B^\delta,\qquad
\frac{P_0}{B}\le B^{-\delta},\qquad
\frac{P_0Z_e}{T^2}\le B^{-\delta}\quad\text{for every }e,
\label{po:strict-1}
\end{equation}
and
\begin{equation}
\min_iY_i\ge B^\delta,\qquad
\min_jQ_j\ge B^\delta.
\label{po:strict-2}
\end{equation}

Let $W(\mathbf s,\mathbf t)$ be a smooth, possibly complex-valued
function on
\[
(0,\infty)^6\times
\left\{\mathbf t\in\mathbb R^4:\sum_i t_i=0\right\},
\]
with
$
\operatorname{supp}W
\subseteq
\left\{
1\le s_e\le2,\quad
\sum_i t_i=0,\quad
c_0\le |t_i|\le c_0^{-1}
\right\},$
where $0<c_0<1$ is fixed. The weight may depend on $B$.
Fix $N=16$ and define
\[
K_N(W)=
\max_{|\alpha|+|\beta|\le N}
\left\|
\partial_{\mathbf s}^{\alpha}
\partial_{t_1}^{\beta_1}
\partial_{t_2}^{\beta_2}
\partial_{t_3}^{\beta_3}
\bigl[
W(\mathbf s,t_1,t_2,t_3,-t_1-t_2-t_3)
\bigr]
\right\|_\infty.
\]

At a tuple of integer values of \(z_e\), use
\begin{equation}
T_z=(BP(z))^{1/3},\qquad
n_i=A_i(z)y_i,\qquad
t_i=\frac{n_i}{T_z},\qquad
s_e=\frac{z_e}{Z_e}.
\label{po:adaptive}
\end{equation}
On the prescribed support,
\[
P(z)\asymp P_0,\qquad
T_z\asymp T,\qquad
A_i(z)\asymp A_i(\mathbf Z),\qquad
|y_i|=\frac{T_z|t_i|}{A_i(z)}\asymp Y_i,
\]
with constants depending only on the fixed support.

By \eqref{in:torsor-monomials},
\[
A_i(z)C_i(z)=P(z),\qquad
\prod_i A_i(z)=P(z)^2.
\]
Thus the raw coordinate formula in \eqref{in:torsor} gives
\begin{equation}
\frac{x_{\mathrm{raw},i}}{B}
=
\frac{C_i(z)A_i(z)T_z^3}{BP(z)^2}
\prod_{j\ne i}t_j
=
\prod_{j\ne i}t_j.
\label{po:true-height}
\end{equation}
In particular,
\[
\frac{H_{\mathrm{raw}}(z,y)}{B}
=
\max_i\left|\prod_{j\ne i}t_j\right|,
\qquad
\frac{|x_{\mathrm{raw},i}|}{|x_{\mathrm{raw},j}|}
=
\frac{|t_j|}{|t_i|}.
\]
We may therefore include smooth weight functions in \(W\)
to restrict these coordinate heights and ratios.
Define
\[
J_W(\mathbf s)
=
\int_{\mathbb R^3}
W(\mathbf s,t_1,t_2,t_3,-t_1-t_2-t_3)
\,dt_1\,dt_2\,dt_3
\]
and
\begin{equation}
\mathcal V_W(B)
=
B\int_{(0,\infty)^6}
J_W(\mathbf s)\prod_e\frac{ds_e}{s_e}.
\label{po:continuous}
\end{equation}
With the Leray-measure convention of \eqref{ar:coarea},
\begin{equation}
\mathcal V_W(B)
=
\int_{z>0}\int_{\mathbb R^4}
W\left(
(z_e/Z_e)_e,\,
(A_i(z)y_i/T_z)_i
\right)
\delta_0\left(\sum_iA_i(z)y_i\right)
\,d^4y\,d^6z.
\label{po:dirac-measure}
\end{equation}

For fixed $z$, the change $A_i(z)y_i=T_zt_i$, followed by
$t_4=-t_1-t_2-t_3$, gives, with $d^3t=dt_1\,dt_2\,dt_3$,
\[
\delta_0\!\left(\sum_iA_i(z)y_i\right)d^4y
=\frac{T_z^3}{\prod_iA_i(z)}\,d^3t
=\frac{B}{P(z)}\,d^3t,
\qquad
\frac{d^6z}{P(z)}=\prod_e\frac{ds_e}{s_e}.
\]
This proves \eqref{po:dirac-measure}; all integrals are absolutely
convergent on the fixed normalized support.

For fixed $R\ge2$, define
\[
\begin{aligned}
N_R(W)
={}
\sum_{\substack{(z,y)\\ \text{y-relaxed}}}
W\left((z_e/Z_e)_e,(n_i/T_z)_i\right)
\prod_{\substack{p\le R\\p\ \mathrm{prime}}}
\prod_{1\le i<j\le4}
\mathbf{1}_{p\nmid\gcd(y_i,y_j)}
\mathbf{1}_{p\nmid\gcd(z_{ij},y_i y_j)}.
\end{aligned}
\]
The y-relaxed tuples retain the torsor equation, pairwise
coprimality of the six variables \(z_e\), and all three pair-sum exclusions.
The displayed indicators impose both kinds of y-coprimality
at primes $p\le R$.

Both overall signs are included in $N_R(W)$ and
$\mathcal V_W(B)$; no projective factor $1/2$ is inserted.

\subsection{The original lattice and finite residue classes}

Fix pairwise coprime positive integers $z_e$, and write
$A_i=A_i(z)$ and $P=P(z)$.
With $n_4=-n_1-n_2-n_3$, define the full solution lattice
\[
\begin{aligned}
\Gamma_z
&=
\left\{
(A_1y_1,A_2y_2,A_3y_3):
y\in\mathbb Z^4,\ \sum_iA_i y_i=0
\right\}\\
&=
\left\{
(n_1,n_2,n_3)\in\mathbb Z^3:
A_i\mid n_i\ (i\le3),\
A_4\mid n_1+n_2+n_3
\right\}.
\end{aligned}
\]
Pairwise coprimality of the variables \(z_e\) implies
\[
p\mid z_{ij}\Longrightarrow p\nmid A_iA_j,
\qquad
\gcd(A_1,A_2,A_3,A_4)=1.
\]
Consequently, the homomorphism
\[
\mathbb Z^3\longrightarrow\mathbb Z/A_4\mathbb Z,
\qquad
(y_1,y_2,y_3)\longmapsto
A_1y_1+A_2y_2+A_3y_3
\]
is surjective. Hence
$[A_1\mathbb Z\times A_2\mathbb Z\times A_3\mathbb Z:\Gamma_z]
=A_4,$
and therefore
\begin{equation}
\det\Gamma_z=A_1A_2A_3A_4=P^2.
\label{po:det}
\end{equation}
Here and below, $\det$ denotes covolume in the three displayed
coordinates.

Identify $(h_1,h_2,h_3,0)$ with $(h_1,h_2,h_3)$. Then
\begin{equation}
\Gamma_z^*=P^{-1}L_z,\qquad
L_z=
\left\{
\mathbf h\in\mathbb Z^4:
h_4=0,\ z_{ij}\mid h_i-h_j\ (i<j)
\right\}.
\label{po:dual}
\end{equation}
Indeed, for $\mathbf n\in\Gamma_z$ and complementary pairs
$\{i,j\},\{k,l\}$,
\[
z_{ij}\mid n_k,\ n_l,\qquad
n_i+n_j\equiv0\pmod{z_{ij}}.
\]
Thus, for $\mathbf h\in L_z$,
\[
\sum_{a=1}^4h_an_a
\equiv(h_i-h_j)n_i
\equiv0\pmod{z_{ij}}.
\]
Since the six variables \(z_e\) are pairwise coprime, one has
$P\mid\sum_{a=1}^4h_an_a$ and
$P^{-1}L_z\subseteq\Gamma_z^*.$
The Chinese remainder theorem gives the surjection
\[
\mathbb Z^3\longrightarrow
\prod_{i<j}\mathbb Z/z_{ij}\mathbb Z,
\qquad
\mathbf h\longmapsto(h_i-h_j)_{i<j},
\qquad h_4=0.
\]
Indeed, each component linear form has a coefficient $1$ or $-1$,
and its modulus is coprime to the other five moduli.
Its kernel is $L_z$, so $[\mathbb Z^3:L_z]=P.$
Using the standard covolume identities
\cite[Section~1, equations~(3)--(6)]{Elkies2009},
\[
\det(P^{-1}L_z)
=P^{-3}\det L_z
=P^{-2}
=(\det\Gamma_z)^{-1}
=\det\Gamma_z^*.
\]
The inclusion is therefore equality, proving \eqref{po:dual}.

For fixed $R\ge2$, put
\[
W_R=\prod_{\substack{p\le R\\p\ \mathrm{prime}}}p.
\]
Since $\gcd(A_1,A_2,A_3,A_4)=1$, every solution modulo $W_R$
lifts to an integral solution: if $\sum_iA_i y_i=W_Rk$, choose
$\sum_iA_i b_i=k$ and replace $y$ by $y-W_Rb$.
The kernel of reduction is $W_R$ times the integral solution lattice.
Consequently,
\[
\Gamma_z/W_R\Gamma_z\simeq
\{\bar y\in(\mathbb Z/W_R\mathbb Z)^4:\textstyle\sum_iA_i\bar y_i=0\},
\qquad \#(\Gamma_z/W_R\Gamma_z)=W_R^3.
\]

The finite-prime coprimality conditions select a union of these
cosets. The pair-sum exclusions are imposed separately using
\eqref{po:line-removal}; coordinate-zero contributions already
vanish by the support condition on $W$. For each prime \(p\), let \(a_{0,p}\) and \(a_{1,p}\)
denote the proportions of solutions to
\(\sum_i A_i y_i=0\pmod p\) satisfying the coprimality
conditions at \(p\), in the cases \(p\nmid P(z)\) and
\(p\mid P(z)\), respectively.

If $p\nmid P$, the substitution $n_i=A_i y_i$ is invertible
modulo $p$. The forbidden locus in $\sum_i n_i=0$ is the union
of the six lines $n_i=n_j=0$. Distinct lines intersect only at
the origin, so
\[
\#\left\{
\mathbf n\in\mathbb F_p^4:
\sum_i n_i=0,\
n_i=n_j=0\text{ for some }i<j
\right\}
=1+6(p-1).
\]
The allowed proportion is therefore
\begin{equation}
a_{0,p}
=
\frac{p^3-1-6(p-1)}{p^3}
=
1-6p^{-2}+5p^{-3}.
\label{po:local-none}
\end{equation}

If $p\mid P$, it divides exactly one variable, say $z_{ij}$.
For the complementary indices $k,l$, we have
$p\nmid A_iA_j$ and $p\mid A_k,A_l.$
The local conditions \(p\nmid\gcd(y_r y_s,z_{rs})\) and
\(p\nmid\gcd(y_r,y_s)\) for \(r<s\) reduce to
\[
A_i y_i+A_j y_j=0,\qquad
y_i,y_j\in\mathbb F_p^\times,\qquad
(y_k,y_l)\ne(0,0).
\]
There are $p-1$ choices for $(y_i,y_j)$ and $p^2-1$ choices
for $(y_k,y_l)$. Hence
\begin{equation}
a_{1,p}
=
\frac{(p-1)(p^2-1)}{p^3}
=
(1-p^{-1})(1-p^{-2}).
\label{po:local-one}
\end{equation}
These formulas also hold for $p=2$ and depend only on whether
$p$ divides $P$, not on the positive valuations of the variables \(z_e\).

By the Chinese remainder theorem,
\[
\frac{\#\{\text{admissible cosets of }W_R\Gamma_z\}}{W_R^3}
=
\prod_{\substack{p\le R\\p\nmid P}}a_{0,p}
\prod_{\substack{p\le R\\p\mid P}}a_{1,p}.
\]
By lattice Poisson summation
\cite[Theorem~2, equation~(26)]{Elkies2009}, the zero-frequency
contribution before removing the pair-sum planes is the following,
since $T_z^3/P^2=B/P$:
\[
\frac{B}{P}J_W(z/\mathbf Z)
\prod_{\substack{p\le R\\p\nmid P}}a_{0,p}
\prod_{\substack{p\le R\\p\mid P}}a_{1,p}.
\]

\subsection{Removing the pairing planes without an additive constant}

Fix a pairing $ij|kl$ of $\{1,2,3,4\}$.
By \eqref{in:torsor-monomials} and pairwise coprimality of the variables \(z_e\),
$\gcd(A_i(z),A_j(z))=z_{kl},$
and hence
\[
\operatorname{lcm}(A_i(z),A_j(z))
=\frac{P(z)}{z_{ij}},
\qquad
\operatorname{lcm}(A_k(z),A_l(z))
=\frac{P(z)}{z_{kl}}.
\]
If $\mathbf n\in\Gamma_z$ lies on the plane $n_i+n_j=0$, then
$n_j=-n_i$ and $n_l=-n_k$. Hence
\[
n_i\in\frac{P(z)}{z_{ij}}\mathbb Z,\qquad
n_k\in\frac{P(z)}{z_{kl}}\mathbb Z.
\]
The pair $(n_i,n_k)$ determines the entire tuple.
On the support of $W$, \eqref{po:adaptive} gives
\[
0<|n_i|,\ |n_k|\le c_0^{-1}T_z.
\]
Therefore
\[
\begin{aligned}
\#\left\{
\mathbf n\in\Gamma_z:
n_i+n_j=0,\
W(z/\mathbf Z,n/T_z)\ne0
\right\}
&\le
4\left\lfloor\frac{T_z z_{ij}}{c_0P(z)}\right\rfloor
 \left\lfloor\frac{T_z z_{kl}}{c_0P(z)}\right\rfloor\\
&\le
\frac{4T_z^2z_{ij}z_{kl}}{c_0^2P(z)^2}
\ll
\frac{T^2Z_{ij}Z_{kl}}{P_0^2}.
\end{aligned}
\]
This bound has no additive constant: if either floor is zero,
there are no contributing points.

On the prescribed block in the variables \(z_e\),
\[
P(z)\asymp P_0,\qquad T_z\asymp T,\qquad
\#\{z\in\mathbb Z_{>0}^6:Z_e\le z_e<2Z_e\}
\le\prod_e(Z_e+1)\ll P_0.
\]
Summing over the three pairings and using
$T^3=BP_0$ and $Q_r=T/R_r$, we obtain
\begin{equation}
\begin{aligned}
\sum_{\substack{Z_e\le z_e<2Z_e\\
                 \gcd(z_e,z_f)=1\ (e\ne f)}}
 \ \sum_{\substack{\mathbf n\in\Gamma_z\\
                   \prod_{j=2}^4(n_1+n_j)=0}}
 \left|W(z/\mathbf Z,n/T_z)\right|
&\ll
\|W\|_\infty\,\frac{T^2}{P_0}\sum_{r=1}^3R_r\\
&=
\|W\|_\infty\,B\sum_{r=1}^3Q_r^{-1}
\ll
\|W\|_\infty\,B^{1-\delta},
\end{aligned}
\label{po:line-removal}
\end{equation}
where the last inequality is \eqref{po:strict-2}.
The implied constant depends only on the fixed support.

All y-sieve conditions were omitted in this upper bound, so
\eqref{po:line-removal} is uniform in the finite sieve threshold $R$.
The absolute values also cover signed or complex weights.

Each plane $t_i+t_j=0$ has zero three-dimensional Leray measure
inside $\sum_i t_i=0$. Its removal therefore leaves
$\mathcal V_W(B)$ in \eqref{po:continuous} unchanged.
The coordinate planes contribute nothing because
$|t_i|\ge c_0$ on the support of $W$.

\subsection{Summing the nonzero frequencies}

Fix $0<\theta<1$. For $L\ge1$, consider
\[
h_4=0,\qquad
0<\|\mathbf h\|_\infty\le L,\qquad
z_{ij}\mid h_i-h_j.
\]
If $h_i\ne h_j$, the standard divisor bound gives
\[
\#\left\{
z_{ij}\in[Z_{ij},2Z_{ij})\cap\mathbb Z:
z_{ij}\mid h_i-h_j
\right\}
\ll_\theta(1+L)^{\theta/6}.
\]
If $h_i=h_j$, there are $O(Z_{ij})$ choices, since $Z_{ij}\ge1$.

For an equality partition with $b$ classes, the class containing
$4$ has value zero; the remaining classes have $O(L^{b-1})$
choices. The possible nonzero patterns are
\[
\begin{array}{c|c|c}
\text{class sizes}
&\text{number of frequency choices}
&\text{product of unrestricted size parameters}\\ \hline
1+1+1+1 & O(L^3) & 1\\
2+1+1   & O(L^2) & Z_{ij}\\
3+1     & O(L)   & A_i(\mathbf Z)\\
2+2     & O(L)   & R_j .
\end{array}
\]
The pattern $4$ gives only $\mathbf h=0$.
Multiplying at most six divisor bounds and summing the finitely
many equality partitions yields
\begin{equation}
\begin{aligned}
&\sum_{z\asymp\mathbf Z}
\#\left\{
\mathbf h\in L_z:
0<\|\mathbf h\|_\infty\le L
\right\}\\
&\qquad\ll_\theta
(1+L)^\theta
\left\{
L^3+
L^2\sum_eZ_e+
L\left(\sum_iA_i(\mathbf Z)+\sum_jR_j\right)
\right\}.
\end{aligned}
\label{po:frequency-count}
\end{equation}
Pairwise coprimality of the variables \(z_e\) may be omitted in this upper bound.
The estimate is uniform in the independent scales $Z_e$ and
holds for every $L\ge1$.

Use the Fourier transform in the three independent coordinates:
\[
\begin{aligned}
\widehat W(\mathbf s,\boldsymbol\xi)
=
\int_{\mathbb R^3}
W(\mathbf s,t_1,t_2,t_3,-t_1-t_2-t_3)
e^{-2\pi i(t_1\xi_1+t_2\xi_2+t_3\xi_3)}
\,dt_1\,dt_2\,dt_3.
\end{aligned}
\]
The fixed support and integration by parts give
\[
|\widehat W(\mathbf s,\boldsymbol\xi)|
\ll_N K_N(W)(1+\|\boldsymbol\xi\|_\infty)^{-N},
\qquad N=16.
\]

With the negative-exponent Fourier convention fixed above,
the translated form of Poisson summation
\cite[Theorem~2, equation~(32)]{Elkies2009} gives, for each admissible
coset representative $\mathbf r$,
\[
\begin{aligned}
\sum_{\mathbf n\in\mathbf r+W_R\Gamma_z}
W(z/\mathbf Z,\mathbf n/T_z)=
\frac{T_z^3}{W_R^3P(z)^2}
\sum_{\boldsymbol\xi\in(W_R\Gamma_z)^*}
e^{2\pi i\mathbf r\cdot\boldsymbol\xi}
\widehat W(z/\mathbf Z,T_z\boldsymbol\xi).
\end{aligned}
\]
Here the fourth coordinate is always the negative sum of the
first three.

By \eqref{po:dual}, one has
$(W_R\Gamma_z)^*
=
\frac{1}{W_RP(z)}L_z.$
The number of admissible cosets is at most $W_R^3$, and hence
\[
\left|
\sum_{\mathbf r\ \mathrm{admissible}}
e^{2\pi i\mathbf r\cdot\boldsymbol\xi}
\right|
\le W_R^3,
\qquad
\frac{T_z^3}{W_R^3P(z)^2}W_R^3
=\frac{B}{P(z)}.
\]

Let $E_{\mathrm{Pois}}(W)$ be the nonzero-frequency remainder,
summed over the block in the variables \(z_e\), before the pairing-plane subtraction.
Put $L_0=\frac{W_RP_0}{T}.$
Since $P(z)\asymp P_0$, we have
\[
\frac{T_z}{W_RP(z)}
=
\frac1{L_0}
\left(\frac{P(z)}{P_0}\right)^{-2/3}
\asymp L_0^{-1}.
\]
It follows that
\[
\begin{aligned}
|E_{\mathrm{Pois}}(W)|
\ll
K_N(W)\frac{B}{P_0}
\sum_{z\asymp\mathbf Z}
\sum_{\substack{\mathbf h\in L_z\\\mathbf h\ne0}}
\left(1+\frac{\|\mathbf h\|_\infty}{L_0}\right)^{-N}.
\end{aligned}
\]

Apply \eqref{po:frequency-count} on dyadic shells.
For $L_0\ge1$, the contribution of the $L^k$ term is
$O_{\theta,N}(L_0^{k+\theta})$ for $k=1,2,3$.
For $L_0<1$, nonzero integer frequencies have norm at least one,
and the same contribution is
\[
O_{\theta,N}\!\left(
L_0^N\sum_{m\ge0}2^{m(k+\theta-N)}
\right)=O_{\theta,N}(L_0^N)=O_{\theta,N}(L_0^k).
\]
Both estimates hold because $N=16>3+\theta$.
Consequently,
\[
\begin{aligned}
|E_{\mathrm{Pois}}(W)|
\ll_\theta
K_N(W)\frac{B}{P_0}(1+L_0)^\theta
\Bigg\{
L_0^3+L_0^2\sum_eZ_e
+L_0\left(\sum_iA_i(\mathbf Z)+\sum_jR_j\right)
\Bigg\}.
\end{aligned}
\]

Since $P_0\le B^{1-\delta}$, we have
$L_0=W_RB^{-1/3}P_0^{2/3}
\le W_RB^{(1-2\delta)/3}.$
For fixed $R$, choosing $0<\theta<\min(1,3\varepsilon)$ gives
$(1+L_0)^\theta\ll_{R,\varepsilon}B^\varepsilon$.
Thus the divisor loss is absorbed after summing the unbounded
frequency ranges.
Using $L_0=W_RP_0/T$ and $T^3=BP_0$, we obtain
\begin{equation}
\begin{aligned}
|E_{\mathrm{Pois}}(W)|
&\ll_{R,\varepsilon}
K_N(W)B^{1+\varepsilon}
\left\{
\frac{P_0}{B}
+\frac{P_0}{T^2}\sum_eZ_e
+\frac{\sum_iA_i(\mathbf Z)+\sum_jR_j}{T}
\right\}\\
&\ll_{R,\varepsilon}
K_N(W)B^{1-\delta+\varepsilon}.
\end{aligned}
\label{po:frequency}
\end{equation}
Indeed, by \eqref{po:strict-1}--\eqref{po:strict-2},
\[
\begin{aligned}
\frac{P_0}{B}
+\frac{P_0}{T^2}\sum_eZ_e
+\frac{\sum_iA_i(\mathbf Z)+\sum_jR_j}{T}
=
\frac{P_0}{B}
+\sum_e\frac{P_0Z_e}{T^2}
+\sum_iY_i^{-1}
+\sum_jQ_j^{-1}
\le14B^{-\delta}.
\end{aligned}
\]

\paragraph{Unsieved form used in Section~\ref{ta:section}.}
The first bound in \eqref{po:frequency} also holds under the weaker
assumptions of fixed normalized support and $1\le Z_e\le B^C$ for a
fixed $C>0$: then $L_0\le W_RB^{4C-1/3}$, and the same argument
absorbs $(1+L_0)^\theta$ by taking $\theta$ sufficiently small in terms
of $C$ and $\varepsilon$. In this generality the implied constant
may depend on $C$; the strict inequalities are used for the second bound.
No derivatives in the variables \(z_e\) enter this frequency estimate. Thus its $K_N(W)$
may be replaced by the uniform $C^N$ norm in the three independent
real coordinates, allowing arbitrary dependence on the parameters \(z_e\).
Taking $W_R=1$ gives the unsieved version. This argument does not use
the near-line estimates of Section~\ref{ta:section}.

\subsection{Weighted averaging over the six auxiliary variables}

Define
\[
f_R(z)=
\mathbf1_{\{\gcd(z_e,z_f)=1\ (e\ne f)\}}
\prod_{\substack{p\le R\\p\nmid P(z)}}a_{0,p}
\prod_{\substack{p\le R\\p\mid P(z)}}a_{1,p}.
\]
At primes $p>R$, only pairwise coprimality of the variables \(z_e\) is imposed.

\begin{lemma}
\label{po:boolean}
Fix $R\ge2$ and $0<\theta<1$.
Let $F\in C_c^1((0,\infty)^6)$ be supported in $[1,2]^6$.
Uniformly for $Z_e>0$,
\begin{equation}
\begin{aligned}
\sum_{z\in\mathbb Z_{>0}^6}
\frac{f_R(z)}{P(z)}F(z/\mathbf Z)
={}
\mathfrak S_R
\int_{(0,\infty)^6}F(\mathbf s)\prod_e\frac{ds_e}{s_e}
+
O_{R,\theta}\left(
\|F\|_{C^1}\sum_e Z_e^{-\theta}
\right),
\end{aligned}
\label{po:boolean-average}
\end{equation}
where
\begin{equation}
\mathfrak S_R
=
\prod_{p\le R}(1-p^{-1})^7(1+7p^{-1}+p^{-2})
\prod_{p>R}(1-p^{-1})^5(1+5p^{-1}).
\label{po:SR}
\end{equation}
The norm is taken in the normalized coordinates $\mathbf s$.
\end{lemma}

\begin{proof}
For a subset $S$ of the six index pairs, put
\begin{equation}
h_p(S)=
\begin{cases}
(-1)^{|S|}\bigl(a_{0,p}-|S|a_{1,p}\bigr),&p\le R,\\
(-1)^{|S|}(1-|S|),&p>R.
\end{cases}
\label{po:boolean-local}
\end{equation}
Set $h_R(d)=0$ unless every coordinate $d_e$ is squarefree;
otherwise set
$h_R(d)=\prod_p h_p\bigl(\{e:p\mid d_e\}\bigr).$
Only finitely many factors differ from $1$.
Coordinatewise M\"obius inversion
\cite[Chapter~I, Theorem~2.3]{SS} gives
\[
f_R(z)=\sum_{d_e\mid z_e\ \text{for all }e}h_R(d).
\]
Indeed, summing \eqref{po:boolean-local} over subsets of a
divisibility pattern gives the respective local weights
$a_{0,p},a_{1,p},0$ for patterns of sizes $0,1,\ge2$ when
$p\le R$, and $1,1,0$ when $p>R$.

For $p>R$,
\[
h_p(\varnothing)=1,\qquad
h_p(\{e\})=0,\qquad
|h_p(S)|\le5\quad(|S|\ge2).
\]
Thus, for each fixed index pair $e$,
\[
\sum_S
\frac{|h_p(S)|\,p^{\theta\mathbf1_{e\in S}}}{p^{|S|}}
=
1+O(p^{-2+\theta}).
\]
Since $\theta<1$, the corresponding Euler products converge.
The finitely many factors with $p\le R$ are finite, so
\begin{equation}
\sum_d\frac{|h_R(d)|}{P(d)}
+
\sum_d\frac{|h_R(d)|d_e^\theta}{P(d)}
\ll_{R,\theta}1
\qquad\text{for every }e.
\label{po:fractional-moments}
\end{equation}

The binomial theorem and \eqref{po:local-none}--\eqref{po:local-one}
give
\[
\sum_S\frac{h_p(S)}{p^{|S|}}
=
\begin{cases}
(1-p^{-1})^7(1+7p^{-1}+p^{-2}),&p\le R,\\
(1-p^{-1})^5(1+5p^{-1}),&p>R.
\end{cases}
\]
Thus \eqref{po:fractional-moments} implies
$\sum_d h_R(d)/P(d)=\mathfrak S_R$.

Let $g\in C_c^1((0,\infty))$ be supported in $[1,2]$.
Euler summation \cite[Chapter~I, Theorem~1.2]{SS} for $0<h\le1$,
and the bound $\sum_{1\le hn\le2}n^{-1}\le2$ for all $h>0$, give
\[
\left|\sum_{n\ge1}\frac{g(hn)}n-
\int_0^\infty g(s)\frac{ds}{s}\right|
\ll\|g\|_{C^1}\min(1,h)\qquad(h>0).
\]
Apply this estimate in each coordinate with mesh $d_e/Z_e$ and
telescope the six sum and integral operators. Their uniformly
bounded masses control the other coordinates, so only first partial
derivatives of $F$ are needed. Since $\min(1,h)\le h^\theta$,
\begin{equation}
\begin{aligned}
\sum_{n\in\mathbb Z_{>0}^6}
\frac{F((d_en_e/Z_e)_e)}{\prod_e n_e}
={}
\int_{(0,\infty)^6}F(\mathbf s)\prod_e\frac{ds_e}{s_e}
+
O\left(
\|F\|_{C^1}\sum_e(d_e/Z_e)^\theta
\right).
\end{aligned}
\label{po:riemann}
\end{equation}

Insert the divisor expansion of $f_R$:
\[
\sum_z\frac{f_R(z)}{P(z)}F(z/\mathbf Z)
=
\sum_d\frac{h_R(d)}{P(d)}
\sum_n\frac{F((d_en_e/Z_e)_e)}{\prod_e n_e}.
\]
Using \eqref{po:riemann}, the main term is
\[
\left(\sum_d\frac{h_R(d)}{P(d)}\right)
\int F(\mathbf s)\prod_e\frac{ds_e}{s_e}
=
\mathfrak S_R\int F(\mathbf s)\prod_e\frac{ds_e}{s_e}.
\]
The error is at most
\[
\ll
\|F\|_{C^1}
\sum_e Z_e^{-\theta}
\sum_d\frac{|h_R(d)|d_e^\theta}{P(d)}
\ll_{R,\theta}
\|F\|_{C^1}\sum_e Z_e^{-\theta}.
\]
The uniform sum bounds and \eqref{po:fractional-moments}
justify the rearrangements and termwise estimates.
This proves \eqref{po:boolean-average}.
\end{proof}

\subsection{The Poisson main term and removal of the finite sieve}

\begin{theorem}
\label{po:main}
Fix $R\ge2$, $\delta>0$, and $\varepsilon>0$.
Assume \eqref{po:strict-1}--\eqref{po:strict-2}, and let $W$
satisfy the stated support conditions, with $N=16$.
Then, uniformly for $B\ge2$,
\begin{equation}
N_R(W)
=
\mathfrak S_R\mathcal V_W(B)
+
O_{R,\delta,\varepsilon}
\left(K_N(W)B^{1-\delta/2+\varepsilon}\right).
\label{po:main-formula}
\end{equation}
The implied constant may also depend on the fixed support,
but is independent of $B$, the individual scales, and $W$.
Both overall signs are counted; no projective factor $1/2$
is included.
\end{theorem}

\begin{proof}
By the fixed support,
$J_W\in C_c^1((0,\infty)^6)$, $\operatorname{supp}J_W\subseteq[1,2]^6$, and
$\|J_W\|_{C^1}\ll_{c_0}K_N(W)$.
Poisson summation and \eqref{po:line-removal} give
\[
N_R(W)=B\sum_{z\in\mathbb Z_{>0}^6}\frac{f_R(z)}{P(z)}J_W(z/\mathbf Z)
+E_{\mathrm{Pois}}(W)
+O\!\left(\|W\|_\infty B\sum_jQ_j^{-1}\right).
\]
Apply Lemma~\ref{po:boolean} with $F=J_W$ and $\theta=1/2$,
and use \eqref{po:continuous} and \eqref{po:frequency}. The error is
\[
\ll_{R,\delta,\varepsilon}K_N(W)
\left(B\sum_e Z_e^{-1/2}+B^{1-\delta+\varepsilon}
+B\sum_jQ_j^{-1}\right)
\ll_{R,\delta,\varepsilon}K_N(W)B^{1-\delta/2+\varepsilon}.
\]
This proves \eqref{po:main-formula}.
\end{proof}

\paragraph{Weights with growing derivatives.}
Theorem~\ref{po:main} also permits weights with growing derivatives:
if $K_{16}(W)\ll B^{\delta/8}$ with a uniform implied constant,
then, with
$\varepsilon=\delta/8$, its error is $O_{R,\delta}(B^{1-\delta/4})$.
For a cutoff of normalized width $B^{-\kappa}$, this applies
provided the full weight satisfies $K_{16}(W)\ll B^{16\kappa}$
with a uniform implied constant and $16\kappa\le\delta/8$.
The Poisson weights in Section~\ref{gl:section} have uniformly
bounded derivatives and do not require this additional restriction.

Let $\mathfrak S$ be as in Theorem~\ref{ar:main}.
Both local factors in \eqref{po:SR} are positive and equal to
$1+O(p^{-2})$. Hence the products converge absolutely and
\[
\left|\log\frac{\mathfrak S_R}{\mathfrak S}\right|
\ll\sum_{p>R}p^{-2}\ll R^{-1}.
\]
Consequently,
\begin{equation}\label{po:S-limit}
\mathfrak S_R=\mathfrak S+O(R^{-1})
\qquad(R\ge2).
\end{equation}

\begin{corollary}
\label{po:core-main}
For each \(B\ge2\), consider a family \(\{W_\beta\}\), independent of
\(R\), on \(O((\log B)^6)\) blocks satisfying
\eqref{po:strict-1}--\eqref{po:strict-2}.
Here \(\delta>0\), the implied constant in the block count, the common
normalized support constants, and \(C_0,C_1>0\) below are fixed
independently of \(B\) and \(R\). Assume
$K_N(W_\beta)\le B^{\delta/8}$ for every $\beta,$
and, on the real torsor,
\[
\sum_\beta
\left|
W_\beta\left(
(z_e/Z_{\beta,e})_e,\,
(A_i(z)y_i/T_z)_i
\right)
\right|
\le C_0.
\]
Suppose that the union of these supports satisfies
$H_{\mathrm{raw}}(z,y)\le C_1B$, where $C_0,C_1>0$
are fixed.

Let $N_{\mathrm{full}}(W_\beta)$ be defined as
$N_R(W_\beta)$ with
$\gcd(y_i,y_j)=\gcd(z_{ij},y_i y_j)=1$ for
$1\le i<j\le4$.
Both overall signs are counted. Then
\begin{equation}
\sum_\beta N_{\mathrm{full}}(W_\beta)
=
\mathfrak S\sum_\beta\mathcal V_{W_\beta}(B)
+
o\!\left(B(\log B)^6\right).
\label{po:core-main-formula}
\end{equation}
\end{corollary}

\begin{proof}
Fix $R\ge2$. Taking $\varepsilon=\delta/8$ in
Theorem~\ref{po:main} gives
$K_N(W_\beta)B^{1-\delta/2+\varepsilon}
\le B^{1-\delta/4}.$
Summing over the blocks,
\[
\sum_\beta N_R(W_\beta)
=
\mathfrak S_R\sum_\beta\mathcal V_{W_\beta}(B)
+
O_{R,\delta}\!\left(
B^{1-\delta/4}(\log B)^6
\right).
\]

A tuple satisfying the finite sieve but failing the full
coprimality conditions has a forbidden prime divisor $p>R$.
Grouping these tuples across all blocks and applying the
absolute-overlap bound, Corollary~\ref{po:free-tail} yields
\[
\begin{aligned}
\left|
\sum_\beta
\bigl(N_{\mathrm{full}}(W_\beta)-N_R(W_\beta)\bigr)
\right|
\le
C_0\mathcal T_{\mathrm{free},>R}(C_1B)
\ll
B\log^6(2B)R^{-5/288}.
\end{aligned}
\]
The implied constant is independent of $B$ and $R$.

For any normalized point $(\mathbf s,\mathbf t)$ with
$\sum_i t_i=0$, the substitution
\[
z_e=Z_{\beta,e}s_e,
\qquad
y_i=\frac{(BP(z))^{1/3}}{A_i(z)}t_i
\]
lies on the real torsor. The overlap hypothesis therefore
implies $\|W_\beta\|_\infty\le C_0$.
Using the fixed support and \eqref{po:continuous},
\[
\left|\mathcal V_{W_\beta}(B)\right|
\le
C_0B(\log2)^6(2c_0^{-1})^3,
\qquad
\sum_\beta
\left|\mathcal V_{W_\beta}(B)\right|
\ll B(\log B)^6.
\]
Consequently, \eqref{po:S-limit} gives
\[
\left|
(\mathfrak S_R-\mathfrak S)
\sum_\beta\mathcal V_{W_\beta}(B)
\right|
\ll
B(\log B)^6R^{-1},
\]
again with a constant independent of $B$ and $R$.

Combining the three estimates,
\[
\begin{aligned}
\sum_\beta N_{\mathrm{full}}(W_\beta)
-\mathfrak S\sum_\beta\mathcal V_{W_\beta}(B)
=
O_{R,\delta}\!\left(
B^{1-\delta/4}(\log B)^6
\right)
+
O\!\left(
B\log^6(2B)
\bigl(R^{-5/288}+R^{-1}\bigr)
\right).
\end{aligned}
\]
The second implied constant is independent of $R$.
Hence, for every fixed $R\ge2$,
\[
\limsup_{B\to\infty}
\frac{
\left|
\sum_\beta N_{\mathrm{full}}(W_\beta)
-\mathfrak S\sum_\beta\mathcal V_{W_\beta}(B)
\right|
}{
B(\log B)^6
}
\ll R^{-5/288}+R^{-1}.
\]
Letting $R\to\infty$ proves
\eqref{po:core-main-formula}.
\end{proof}

\section{The common volume and Peyre's constant}\label{co:section}

We compute Peyre's constant and compare it with the coefficient of
the common continuous main term. To make the normalizations explicit,
the calculations below give, with
\(I=\int_0^2\log(1+s)\,ds/s\),
\[
 \alpha(Y)=\frac1{1920},\qquad
 \operatorname{Vol}_6(\Omega)=3\alpha(Y)=\frac1{640},\qquad
 J_\infty=16I=\frac23\omega_\infty.
\]
Consequently,
\[
 \frac12\operatorname{Vol}_6(\Omega)J_\infty
 =\alpha(Y)\omega_\infty=\frac I{80}.
\]
Here \(\Omega\) is the logarithmic polytope introduced below,
\(J_\infty\) is the integral in the torsor coordinates, and
\(\omega_\infty\) is the real Peyre density.
The Picard-slice factor \(1/3\) and the projective factor \(1/2\)
are kept explicit throughout. The finite-place calculation then
supplies the Euler product in \eqref{in:constant}.

\subsection{The resolution and the logarithmic polytope}

Let
\[
\Pi=
\left\{[t_1:t_2:t_3:t_4]\in\mathbb P^3:
\sum_i t_i=0\right\}
\simeq\mathbb P^2.
\]
By \cite[Example~14]{DT2007}, the minimal resolution $Y\to S$
is the blowup of the six intersections of the four lines
$t_i=0$ in $\Pi$. Its anticanonical morphism resolves
\[
[t]\longmapsto[\Phi_1(t):\cdots:\Phi_4(t)],
\qquad
\Phi_i(t)=\prod_{j\ne i}t_j.
\]
On the open subset with nonzero coordinates, the inverse is
$
[t_1:t_2:t_3:t_4]
=
[x_1^{-1}:x_2^{-1}:x_3^{-1}:x_4^{-1}].$

Let $L$ be the pullback of a line, and let $E_{ij}$ lie over
$t_k=t_l=0$, where $\{i,j,k,l\}=\{1,2,3,4\}$. Then
\[
\operatorname{Pic}(Y)
=
\mathbb ZL\oplus\bigoplus_e\mathbb ZE_e,
\qquad
-K_Y=3L-\sum_eE_e.
\]
By \cite[Example~14]{DT2007}, the thirteen negative curves
have the following classes:
\begin{equation}
\begin{aligned}
&E_e
&&\text{for the six index pairs},\\
&D_i=L-\sum_{e:\,i\notin e}E_e
&&(1\le i\le4),\\
&G_{ij|kl}=L-E_{ij}-E_{kl}
&&(ij|kl=12|34,\ 13|24,\ 14|23).
\end{aligned}
\label{co:negative}
\end{equation}
Here $D_i$ and $G_{ij|kl}$ are the strict transforms of
$t_i=0$ and $t_i+t_j=0$, respectively. The standard blowup
intersection form gives four disjoint curves $D_i$ and
\[
D_i^2=-2,\qquad G_{ij|kl}^2=-1,\qquad (-K_Y)^2=3.
\]
Thus $Y$ is a split generalized del Pezzo surface of degree
$3$ and type $4A_1$. By \cite[Theorem~3.10]{DJT2008},
\[
\overline{\operatorname{Eff}}(Y)
=
\operatorname{cone}\{E_e,D_i,G_{ij|kl}\}.
\]

For $v\in\operatorname{Pic}(Y)_{\mathbb R}^{\vee}$, write
$\ell=\langle v,L\rangle$ and
$u_e=\langle v,E_e\rangle.$
The anticanonical slice is
\[
\langle v,-K_Y\rangle=1
\quad\Longleftrightarrow\quad
3\ell-\sum_eu_e=1
\quad\Longleftrightarrow\quad
\ell=\ell(u):=\frac{1+\sum_eu_e}{3}.
\]
Its intersection with the dual effective cone is identified with
\begin{equation}
\Omega=
\left\{
u\in\mathbb R^6:
\begin{array}{ll}
u_e\ge0
&\text{for every index pair},\\[2pt]
\displaystyle
\ell(u)-\sum_{e:\,i\notin e}u_e\ge0
&(1\le i\le4),\\[4pt]
\ell(u)-u_{ij}-u_{kl}\ge0
&\text{for the three pairings}.
\end{array}
\right\}.
\label{co:polytope}
\end{equation}
Summing the four inequalities indexed by \(i\) gives
$0\le
4\ell(u)-2\sum_eu_e
=
\frac{4-2\sum_eu_e}{3},$
hence
\[
\Omega\subseteq
\left\{u\in\mathbb R_{\ge0}^6:\sum_eu_e\le2\right\}.
\]
For $B>1$ and $Z_e=B^{u_e}$,
\eqref{in:scales}--\eqref{in:short-scales} give
\[
T=(BP_0)^{1/3}=B^{\ell(u)},\quad
\frac{\log Y_i}{\log B}
=\ell(u)-\sum_{e:\,i\notin e}u_e,\quad
\frac{\log\bigl(T/(Z_{ij}Z_{kl})\bigr)}{\log B}
=\ell(u)-u_{ij}-u_{kl}.
\]
These identities identify the inequalities indexed by \(i\) and by \(ij|kl\)
with the corresponding logarithmic scale constraints.

The dual integral basis gives measure $d\ell\,d^6u$, and
$d\ell\,d^6u
=
\frac13\,d\!\left(3\ell-\sum_eu_e\right)d^6u.$
Hence the slice normalization in
\cite[Section~2]{DJT2008} gives
\begin{equation}
\alpha(Y)
=
\int_\Omega\frac{d^6u}{3}
=
\frac13\operatorname{Vol}_6(\Omega).
\label{co:slice}
\end{equation}
The root system is $4A_1$, whose Weyl group has order $2^4$.
By \cite[Theorem~1.3 and Table~1]{DJT2008},
\begin{equation}
\alpha(Y)
=
\frac{1}{120\cdot2^4}
=
\frac1{1920},
\qquad
\operatorname{Vol}_6(\Omega)
=
3\alpha(Y)
=
\frac1{640}.
\label{co:volume}
\end{equation}

Since \(Y\) is split,
\[
H^1\!\left(\mathbb Q,\operatorname{Pic}(\overline Y)\right)=0,
\qquad \beta(Y)=1.
\]
Weak approximation on $\mathbb P^2$ transfers to $Y$ through
a common birational open subset, which is dense in
$Y(\mathbb Q_v)$ at every place $v$ because $Y$ is smooth.
Since $Y$ is proper,
$\overline{Y(\mathbb Q)}=Y(\mathbb A_{\mathbb Q}).$
Thus Peyre's adelic closure is the full adelic space;
see \cite[Section~2.1]{Peyre1995}.

\subsection{Every finite place, including two}

Work over
\[
\Pi_{\mathbb Z}
=
\left\{\textstyle\sum_i t_i=0\right\}
\subset\mathbb P^3_{\mathbb Z},
\qquad
\Pi_{\mathbb Z}\simeq\mathbb P^2_{\mathbb Z}.
\]
The four coordinate lines have no triple intersection in any
fiber. Their six pairwise intersections are therefore disjoint
$\mathbb Z$-sections. Blowing them up gives a smooth proper
model $\mathcal Y/\mathbb Z$ of $Y$.

Near $t_3=t_4=0$, write
$(t_1,t_2,t_3,t_4)=(1,-1-a-b,a,b).$
Then
\[
(\Phi_1,\Phi_2,\Phi_3,\Phi_4)
=
\bigl(-ab(1+a+b),\ ab,\ -b(1+a+b),\ -a(1+a+b)\bigr).
\]
Since $1+a+b$ is a unit near the centre, the base ideal is
$(\Phi_1,\Phi_2,\Phi_3,\Phi_4)
\mathcal O_{\Pi_{\mathbb Z}}
=
(a,b).$
Permutation symmetry gives the same calculation at every
centre. Away from these centres the sections have no common
zero. After removing the exceptional factors, they generate
\[
\mathcal O_{\mathcal Y}\left(3L-\sum_eE_e\right)
\simeq
\omega_{\mathcal Y/\mathbb Z}^{-1},
\]
where $L,E_e$ also denote the corresponding relative divisor
classes. Thus the pullback of the coordinate maximum metric
is the integral anticanonical model metric at every prime.

For every prime $p$, each blown-up point is replaced by
$\mathbb P^1_{\mathbb F_p}$, so
\[
\#\mathcal Y(\mathbb F_p)
=
(p^2+p+1)-6+6(p+1)
=
p^2+7p+1.
\]
Properness gives
$Y(\mathbb Q_p)=\mathcal Y(\mathbb Z_p).$
The proof of \cite[Lemma~2.1.2]{Peyre1995}, applied to the
smooth integral model and its anticanonical metric, gives
\[
\omega_p=\frac{\#\mathcal Y(\mathbb F_p)}{p^2}
        =1+\frac7p+\frac1{p^2}
\]
for every prime \(p\).
The removed curves have $p$-adic measure zero.

The geometric Picard group has rank $7$ and trivial Galois
action. Consequently, the convergence factor of
\cite[Section~2.1]{Peyre1995} is
\[
L_p(s,\operatorname{Pic}(\overline Y))
=
(1-p^{-s})^{-7},
\qquad
L_p(1,\operatorname{Pic}(\overline Y))^{-1}
=
(1-p^{-1})^7.
\]
Thus the normalized local factor is
\begin{equation}
\mathfrak S_p
=
L_p(1,\operatorname{Pic}(\overline Y))^{-1}\omega_p
=
(1-p^{-1})^7(1+7p^{-1}+p^{-2}).
\label{co:finite}
\end{equation}

This agrees with the torsor density. The probabilities that
none of the variables \(z_e\) is divisible by $p$, or that exactly
one specified variable is divisible by $p$, are
$(1-p^{-1})^6$ and $p^{-1}(1-p^{-1})^5$, respectively.
Using \eqref{po:local-none}--\eqref{po:local-one},
\[
\begin{aligned}
&(1-p^{-1})^6a_{0,p}
+6p^{-1}(1-p^{-1})^5a_{1,p}\\
&\qquad=
(1-p^{-1})^6
\bigl(1-6p^{-2}+5p^{-3}
      +6p^{-1}(1-p^{-2})\bigr)\\
&\qquad=
(1-p^{-1})^7(1+7p^{-1}+p^{-2})
=
\mathfrak S_p.
\end{aligned}
\]
For either allowed divisibility pattern of the variables \(z_e\), some $A_i(z)$ is a $p$-adic
unit. For fixed values of the variables \(z_e\), each admissible solution modulo $p$
therefore has $p^{3(m-1)}$ lifts modulo $p^m$ for every
$m\ge1$; the coprimality conditions depend only on reduction
modulo $p$.

In particular,
\[
\#\mathcal Y(\mathbb F_2)=19,
\qquad
\omega_2=\frac{19}{4},
\qquad
\mathfrak S_2=\frac{19}{512}.
\]
The positivity and absolute convergence of
\(\mathfrak S=\prod_p\mathfrak S_p\) were established in
\eqref{po:S-limit}.

\subsection{The real density}

On $t_4=-t_1-t_2-t_3$, put
\begin{equation}
\mathcal D
=
\{(t_1,t_2,t_3)\in\mathbb R^3:\max_i|\Phi_i(t)|\le1\},
\qquad
J_\infty=\int_{\mathcal D}dt_1\,dt_2\,dt_3.
\label{co:real-domain}
\end{equation}
Coordinate permutations and simultaneous negation preserve
this measure: their coordinate changes have absolute determinant $1$.
After ordering the coordinates by decreasing absolute value,
for \(|t_1|>R\ge1\) we have
\[
|t_1|\ge|t_2|\ge|t_3|\ge|t_4|,
\qquad
|t_1|\le3|t_2|,
\qquad
|t_3|,|t_4|
\le\frac1{|t_1t_2|}
\le\frac3{|t_1|^2}.
\]
Using $(t_1,t_3,t_4)$ as coordinates and summing over the
finitely many orders gives
\begin{equation}
\operatorname{Vol}\{t\in\mathcal D:\max_i|t_i|>R\}
\ll
\int_R^\infty r^{-4}\,dr
\ll R^{-3}
\qquad(R\ge1).
\label{co:real-tail}
\end{equation}
Thus $J_\infty<\infty$. Coordinate and pair-sum hyperplanes
have measure zero.

In the projective charts $t_4=x_4=1$,
\[
t_3=-1-t_1-t_2,
\qquad
x_i=t_i^{-1}\quad(1\le i\le3).
\]
For $F$ in \eqref{in:surface},
\[
F_{x_3}(x_1,x_2,x_3,1)
=
x_1x_2+x_1+x_2
=
-\frac{t_3}{t_1t_2},
\qquad
|dx_1\,dx_2|
=
\frac{|dt_1\,dt_2|}{|t_1t_2|^2}.
\]
Peyre's local-measure construction
\cite[Section~2.1]{Peyre1995}, with the coordinate maximum
metric, therefore gives
\[
\begin{aligned}
d\omega_\infty
&=
\frac{|dx_1\,dx_2|}
{|F_{x_3}|\max(1,|x_1|,|x_2|,|x_3|)}
=
\frac{|dt_1\,dt_2|}
{|t_1t_2t_3|
 \max(1,|t_1|^{-1},|t_2|^{-1},|t_3|^{-1})}\\
&=
\frac{|dt_1\,dt_2|}
{\max_i|\Phi_i(t_1,t_2,t_3,1)|}.
\end{aligned}
\]
The omitted curves have measure zero, so
\[
\omega_\infty
=
\int_{\mathbb R^2}
\frac{du\,dv}{\max_i|\Phi_i(u,v,-1-u-v,1)|}.
\]

Set
$t=s(u,v,-1-u-v,1).$
This parametrizes \(t_4\ne0\), and
\[
\left|\det\frac{\partial(t_1,t_2,t_3)}
{\partial(s,u,v)}\right|=|s|^2,\qquad
\max_i|\Phi_i(t)|
=|s|^3\max_i|\Phi_i(u,v,-1-u-v,1)|.
\]
The maximum on the right vanishes only at
\((u,v)\in\{(0,0),(0,-1),(-1,0)\}\), a set of measure zero.
Outside this set,
\[
\int_{\{|s|^3\max_i|\Phi_i(u,v,-1-u-v,1)|\le1\}}
|s|^2\,ds
=\frac{2}{3\max_i|\Phi_i(u,v,-1-u-v,1)|}.
\]
Tonelli's theorem yields
\begin{equation}
J_\infty
=\int_{\mathbb R^2}
\frac{2\,du\,dv}{3\max_i|\Phi_i(u,v,-1-u-v,1)|}
=\frac23\omega_\infty.
\label{co:real-normalization}
\end{equation}

There are eight regions with fixed coordinate signs in which
one coordinate has the opposite sign to the other three.
By permutation and simultaneous negation, take
\[
t=s(1,-a,-b,-(1-a-b)),\qquad
s>0,\quad a,b>0,\quad a+b<1.
\]
The Jacobian is $s^2$. The six orders of $a,b,1-a-b$ contribute
equally; on $a\ge b\ge1-a-b>0$,
\[
\max_i|\Phi_i(t)|=s^3ab,\qquad
\int_0^{(ab)^{-1/3}}s^2\,ds=\frac1{3ab}.
\]
With
\[
(a,b,1-a-b)=\frac{(r,1,v)}{r+1+v},\qquad
r\ge1,\quad0<v\le1,
\]
we have $|da\,db|/(ab)=dr\,dv/[r(r+1+v)]$. The total contribution from these eight regions is
\[
\begin{aligned}
\frac{8\cdot6}{3}\int_0^1\int_1^\infty
\frac{dr\,dv}{r(r+1+v)}
=16\int_0^1\frac{\log(2+v)}{1+v}\,dv
=16\int_1^2\frac{\log(1+s)}s\,ds.
\end{aligned}
\]

In each of the remaining six regions, two coordinates are
positive and two are negative. Take
\[
t=s(u,1-u,-v,-(1-v)),\qquad s>0,\quad0<u,v<1.
\]
The Jacobian is again $s^2$. The symmetries $u\mapsto1-u$,
$v\mapsto1-v$ and $(u,v)\mapsto(v,u)$ reduce the square to
$1/2\le v\le u\le1$ with multiplicity $8$. On its interior,
\[
\max_i|\Phi_i(t)|=s^3uv(1-v),\qquad
\int_0^{[uv(1-v)]^{-1/3}}s^2\,ds=\frac1{3uv(1-v)}.
\]
The total contribution from these six regions is
\[
\frac{6\cdot8}{3}\int_{1/2}^1\int_v^1
\frac{du\,dv}{uv(1-v)}
=16\int_{1/2}^1\frac{-\log v}{v(1-v)}\,dv
=16\int_0^1\frac{\log(1+s)}s\,ds,
\]
where the last equality uses $v=(1+s)^{-1}$.
The boundaries of all these regions have measure zero. Thus, with
$I=\int_0^2\frac{\log(1+s)}s\,ds$, the two contributions and
\eqref{co:real-normalization} give
\begin{equation}
J_\infty=16I,\qquad\omega_\infty=24I.
\label{co:real}
\end{equation}

\subsection{The full continuous domain}
Let $z_e\ge1$ and let the $y_i$ be real and signed.  Define
$q_{ij|kl}=\frac{A_i(z)y_i+A_j(z)y_j}{z_{ij}z_{kl}}.$
Let $\mathcal R_B$ be the domain on the torsor hyperplane with
\begin{equation}\label{co:full-domain}
 z_e\ge1,\quad |y_i|\ge1,\quad
 |q_{ij|kl}|\ge1\ \text{for all three pairings},\quad
 \max_i\left|C_i(z)\prod_{j\ne i}y_j\right|\le B.
\end{equation}
Set
\begin{equation}\label{co:V}
 \mathcal V(B)=\int\mathbf1_{\mathcal R_B}
          \delta_0\left(\sum_i A_i(z)y_i\right)d^4y\,d^6z.
\end{equation}
Here \(\delta_0\) is the Dirac distribution at zero; on the
torsor hyperplane,
\[
\delta_0\!\left(\sum_iA_i(z)y_i\right)d^4y
=\frac{dy_1\,dy_2\,dy_3}{A_4(z)},
\qquad
y_4=-\frac{\sum_{i=1}^3A_i(z)y_i}{A_4(z)}.
\]

\begin{proposition}
\label{co:continuous}
For the domain and measure in
\eqref{co:full-domain}--\eqref{co:V},
\[
\mathcal V(B)
\sim
\frac{J_\infty}{640}B(\log B)^6
\qquad(B\to\infty).
\]
\end{proposition}

\begin{proof}
For $B>1$, use
\[
z_e=B^{u_e},
\qquad
T_z=(BP(z))^{1/3}=B^{\ell(u)},
\qquad
t_i=\frac{A_i(z)y_i}{T_z},
\qquad
t_4=-t_1-t_2-t_3.
\]
By \eqref{po:dirac-measure}, the Leray factor is $B\,d^3t/P(z)$.
Since $dz_e=(\log B)z_e\,du_e$, we obtain
\begin{equation}
\delta_0\!\left(\sum_iA_i(z)y_i\right)d^4y\,d^6z
=
B(\log B)^6\,d^6u\,dt_1\,dt_2\,dt_3.
\label{co:Jacobian}
\end{equation}
Tonelli's theorem permits these successive changes of variables.

By \eqref{po:true-height}, the height condition becomes
$t\in\mathcal D$. The remaining conditions are
\[
u_e\ge0,
\qquad
|y_i|
=
B^{\ell(u)-\sum_{e:\,i\notin e}u_e}|t_i|
\ge1,
\qquad
|q_{ij|kl}|
=
B^{\ell(u)-u_{ij}-u_{kl}}|t_i+t_j|
\ge1.
\]
Here $ij|kl$ runs over $12|34,13|24,14|23$.
Thus, writing $d^3t=dt_1\,dt_2\,dt_3$,
\[
\begin{aligned}
\frac{\mathcal V(B)}{B(\log B)^6}
={}&
\int_{\mathbb R_{\ge0}^6}\int_{\mathcal D}
\prod_{i=1}^4
\mathbf1_{\{B^{\ell(u)-\sum_{e:\,i\notin e}u_e}|t_i|\ge1\}}
\prod_{ij|kl}
\mathbf1_{\{B^{\ell(u)-u_{ij}-u_{kl}}|t_i+t_j|\ge1\}}
\,d^3t\,d^6u.
\end{aligned}
\]
Outside the seven affine hyperplanes where an exponent vanishes
and the seven hyperplanes $t_i=0$ or $t_i+t_j=0$, the product
of indicators converges to $\mathbf1_\Omega(u)$ by
\eqref{co:polytope}. These exceptional sets have measure zero.

For every point of $\mathcal R_B$, the monomial identities in
\eqref{in:torsor-monomials} give
\[
P(z)^2
\le
P(z)^2\prod_i|y_i|^3
=
\prod_i\left|C_i(z)\prod_{j\ne i}y_j\right|
\le B^4.
\]
Hence every contributing point satisfies
$u_e\ge0$ and $\sum_eu_e\le2.$
The integrand, extended by zero to
$\mathbb R^6\times\mathbb R^3$, is therefore bounded by
$\mathbf1_{\{u_e\ge0,\ \sum_eu_e\le2\}}
\mathbf1_{\mathcal D}(t).$
This majorant is integrable because its \(u\)-support is
bounded and \(J_\infty<\infty\) by \eqref{co:real-tail}.

Dominated convergence and \eqref{co:volume} now yield
\[
\begin{aligned}
\lim_{B\to\infty}
\frac{\mathcal V(B)}{B(\log B)^6}
=
\int_\Omega d^6u\int_{\mathcal D}d^3t
=
\operatorname{Vol}_6(\Omega)J_\infty
=
\frac{J_\infty}{640}.
\end{aligned}
\]
Since $J_\infty=16I>0$ by \eqref{co:real}, this proves the
asserted asymptotic.
\end{proof}

\begin{corollary}
\label{co:constant}
For the maximum-coordinate height, Peyre's constant is
\[
c_{\mathrm{Peyre}}
=
\frac{I}{80}
\prod_p(1-p^{-1})^7(1+7p^{-1}+p^{-2}).
\]
Moreover,
\[
\frac12\mathfrak S\mathcal V(B)
\sim
c_{\mathrm{Peyre}}B(\log B)^6
\qquad(B\to\infty).
\]
\end{corollary}

\begin{proof}
\[
L\!\left(s,\operatorname{Pic}(\overline Y)\right)=\zeta(s)^7,
\qquad
\lim_{s\to1}(s-1)^7L\!\left(s,\operatorname{Pic}(\overline Y)\right)=1.
\]
Together with \(\beta(Y)=1\) and the full adelic closure established
above, Peyre's construction \cite[Section~2]{Peyre1995} yields
\[
c_{\mathrm{Peyre}}
=\alpha(Y)\omega_\infty\mathfrak S
=\frac{24I}{1920}\mathfrak S
=\frac I{80}\mathfrak S.
\]
Moreover, \eqref{co:slice} and \eqref{co:real-normalization} give
\[
\frac12\operatorname{Vol}_6(\Omega)J_\infty
=\frac12\bigl(3\alpha(Y)\bigr)
\left(\frac23\omega_\infty\right)
=\alpha(Y)\omega_\infty.
\]
The second assertion now follows from
Proposition~\ref{co:continuous}; the factor \(1/2\) identifies
opposite primitive representatives.
\end{proof}

\section{Global assembly}\label{gl:section}

We implement the division of the parameter ranges in the
\hyperref[ta:proof-strategy]{preliminary roadmap}, using a common
partition of unity and the same convention for signed torsor
representatives. The argument first controls the discarded integer
ranges, then verifies the geometric and arithmetic errors on every
retained block, and finally compares the assembled main terms with
the continuous volume. The two-dimensional method uses the complete
sieve; the Poisson method first uses a fixed finite sieve, whose tail
is removed only after the height limit. The final error inequalities
\eqref{gl:limit-error} record the dependence on the remaining
truncations and justify their successive removal.

\subsection{Truncations and the order of choices}

Fix
\[
M>1,\qquad
0<h<\min\{1,\sqrt M-1\},\qquad
0<\delta<\frac1{100},\qquad
0<\kappa<\frac{\delta}{100}.
\]
By Proposition~\ref{ta:height-tail},
\[
\#\{x\in U(\mathbb Q):H(x)\le B,\ \min_i|x_i|<B/M\}
\ll M^{-1/2}B\log^6(2B).
\]
We therefore first restrict to
\begin{equation*}
B/M\le |x_i|\le B\qquad(1\le i\le4).
\end{equation*}

Choose a smooth function \(0\le\chi_h\le1\) satisfying
\[
\chi_h(r)=1\quad(0\le r\le1),\qquad
\chi_h(r)=0\quad(r\ge1+h),\qquad
\|\chi_h^{(j)}\|_\infty\ll_j h^{-j}.
\]
Then
\[
\begin{aligned}
\prod_{i=1}^4
\chi_h\!\left((1+h)\frac{|x_i|}{B}\right)
\chi_h\!\left((1+h)\frac{B}{M|x_i|}\right)
&\le
\mathbf1_{\{B/M\le|x_i|\le B\ (1\le i\le4)\}}\\
&\le
\prod_{i=1}^4
\chi_h\!\left(\frac{|x_i|}{B}\right)
\chi_h\!\left(\frac{B}{M|x_i|}\right).
\end{aligned}
\]
Both products extend smoothly by zero at zero coordinates, and
their supports satisfy
\[
\frac{B}{M(1+h)}\le |x_i|\le(1+h)B.
\]

Choose a fixed nonnegative \(\psi\in C_c^\infty((1,2))\) such that
$\sum_{k\in\mathbb Z}\psi(r/2^{k/2})=1$ for
$r>0.$
Use the product partition
\[
Z_e\in\{2^{k/2}:k\in\mathbb Z\},\qquad
\sum_{\mathbf Z}\prod_e\psi(z_e/Z_e)=1.
\]
These are two interlaced dyadic grids in each coordinate, with
bounded overlap. Proposition~\ref{ta:boundary} applies after
enlarging the logarithmic intervals by \(O(1)\); the other local
estimates use only fixed multiplicative support bounds.

Retain the scales in
\eqref{in:scales}--\eqref{in:short-scales}, ordered so that
\[
A\ge H\ge J,\qquad
(Q_1,Q_2,Q_3)
=\left(\frac TA,\frac TH,\frac TJ\right),
\qquad Q=Q_1.
\]
Resolve ties in a fixed order. On the height supports,
\[
|y_i|\asymp_{M,h}Y_i,\qquad
C_i(\mathbf Z)\prod_{j\ne i}Y_j=B,
\quad\text{and}\quad
\frac{x_i}{B}
=
\frac{C_i(z)}{C_i(\mathbf Z)}
\prod_{j\ne i}\frac{y_j}{Y_j}.
\]
Write the height envelopes in the normalized \(t\)-coordinates as
\[
\begin{aligned}
w_{M,h}^{+}(t)
&=\prod_{i=1}^4
\chi_h\!\left(|\Phi_i(t)|\right)
\chi_h\!\left(\frac1{M|\Phi_i(t)|}\right),\\
w_{M,h}^{-}(t)
&=\prod_{i=1}^4
\chi_h\!\left((1+h)|\Phi_i(t)|\right)
\chi_h\!\left(\frac{1+h}{M|\Phi_i(t)|}\right).
\end{aligned}
\]
Both are extended by zero where some \(\Phi_i(t)=0\).

Take
\[
W^{\pm}(\mathbf s,\mathbf t)
=w_{M,h}^{\pm}(\mathbf t)\prod_e\psi(s_e),
\qquad
\Phi_{\mathbf Z}^{\pm}(z,y)
=W^{\pm}\!\left((z_e/Z_e)_e,(A_i(z)y_i/T_z)_i\right).
\]
Restrict these weights to the chosen sign sector and write the resulting
weights as \(\Phi\) and \(W\). Their normalized supports are fixed for
fixed \(M,h\), and
$\|\Phi\|_{6,*}+K_{16}(W)\ll_{M,h}1.$
For each envelope, the sign sectors are finite and the nonnegative
local weights have total mass at most one.
Retain only blocks satisfying
\begin{equation}\label{gl:core}
\min\{\min_e Z_e,\min_iY_i,\min_jQ_j\}\ge B^\delta\quad\text{and}\quad
P_0\notin(B^{1-\delta},B^{1+\delta}),
\end{equation}
and delete the transition band
\begin{equation}\label{gl:arithmetic-band}
B^{-\delta}<\frac{HJ}{T}<B^\delta.
\end{equation}

Assign a retained block to the Poisson method precisely when
\[
P_0\le B^{1-\delta}
\quad\text{and}\quad A\le F_0B^{3\delta/4}.
\]
Assign every other retained block to the two-dimensional method.
We verify below that each assigned method has a power-saving error.
To estimate the discarded ranges, observe that every block containing a canonical point satisfies
\[
Z_e\gg1,\qquad Y_i\gg_{M,h}1,\qquad
1\le |v_{ij}|
=
\frac{|A_i(z)y_i+A_j(z)y_j|}{z_{ij}z_{kl}}
\ll_{M,h}\frac{T}{Z_{ij}Z_{kl}},
\]
by \eqref{in:q-identity}. Moreover,
\[
\begin{aligned}
\log Y_i
&=\frac13\log B+\frac13\sum_e\log Z_e
  -\log A_i(\mathbf Z),\\
\log\frac{T}{Z_{ij}Z_{kl}}
&=\frac13\log B+\frac13\sum_e\log Z_e
  -\log Z_{ij}-\log Z_{kl},\\
\log\frac{HJ}{T}
&=-\frac13\log B-\frac13\log A
  +\frac23\log H+\frac23\log J.
\end{aligned}
\]
Thus the deleted ranges lie in finitely many slabs of width
\(O_{M,h}(\delta\log B+1)\) for fixed nonconstant affine functions
of \((\log Z_e)_e\). Proposition~\ref{ta:boundary} bounds their
total weighted count by
\[
O_{M,h}\!\left(
\delta B\log^6(2B)+B\log^5(2B)
\right).
\]

In a block assigned to the two-dimensional method, use the
pairing \(ij|kl\) for which \(Z_{ij}Z_{kl}\) is largest
and put \(\rho_0=B^{-\kappa}\).
By \eqref{ar:linear} and \eqref{ar:recover},
\[
q=\gcd(L_1,L_2)
=\frac{|L_1|}{a}
=\frac{|L_2|}{b}
=\frac{|A_1(z)y_1+A_2(z)y_2|}{ab}
=|v_{12}|,
\qquad ab\asymp A.
\]
Since \(AQ=T\), the range
$0<q<\rho_0Q$
and any fixed enlargement of it satisfy
\[
0<|A_1(z)y_1+A_2(z)y_2|
\ll B^{-\kappa}T.
\]
Apply Theorem~\ref{ta:near-line} with height bound \((1+h)B\),
height ratio \(M(1+h)^2\), and core exponent \(\delta/2\).
Under this change of height, the scales \(T,Y_i,Q_j\) are all multiplied
by \((1+h)^{1/3}\), and the cutoff differs from that of the theorem
only by a fixed positive factor. Thus \eqref{gl:core} implies the
hypotheses for sufficiently large \(B\), and the fixed-multiple version
of the theorem applies, since
$\kappa<\frac{\delta}{100}<\frac{\delta/2}{32}.$
There are \(O(\log^6(2B))\) relevant blocks in the variables \(z_e\). Hence the total
near-line contribution is
\[
O_{M,h,\delta,\kappa}\!\left(
B^{1-\kappa/12}\log^6(2B)
\right)
=o\!\left(B(\log B)^6\right).
\]

Implement the truncation through the separate smooth weights
\[
\psi\!\left(\frac{q}{sQ}\right),
\qquad
s=\rho_0\,2^{m/2},\quad m\ge0.
\]
They satisfy
\[
0\le
\sum_{m\ge0}
\psi\!\left(\frac{q}{\rho_0\,2^{m/2}Q}\right)
\le1,
\]
with value zero for \(q\le\rho_0Q\) and value one for
\(q\ge\sqrt2\,\rho_0Q\). The omitted weight is therefore covered
by the same near-line estimate. Since \(q\ll_{M,h}Q\) on the
support,
\[
\rho_0\le s\ll_{M,h}1,\qquad
L=sQ\ge B^{\delta-\kappa}\ge1,\qquad
\#\{\text{active bands}\}\ll_{M,h,\kappa}\log(2B).
\]
Apply Theorem~\ref{ar:main} and Corollary~\ref{ar:count} with
\(\eta=s\) and fixed sufficiently large \(\eta_{\max}\).
We keep the factor \(\psi(q/L)\) separate from \(\Phi\),
so the normalized derivative bounds for \(\Phi\) remain unchanged.

For each fixed \(M,h,\delta\), choose \(\kappa\) as above and fix
the sieve threshold \(R\ge2\) before letting \(B\) grow.
The limits are taken in the order
\[
B\longrightarrow\infty,\qquad
R\longrightarrow\infty,\qquad
\delta\longrightarrow0,\qquad
h\longrightarrow0,\qquad
M\longrightarrow\infty.
\]
The corresponding continuous boundary estimates and removal of
the height envelopes are established in the final subsection.

\subsection{The geometric error on the whole core}\label{gl:geometric-core}

Retain \eqref{gl:core} and fix
\[
0<\kappa<\frac{\delta}{100},\qquad
0<\varepsilon<\frac{\delta}{100},\qquad
B^{-\kappa}\le s\ll_{M,h}1.
\]
For one block in the variables \(z_e\), sign sector, and band, put
\[
\mathcal G=
\sum_{\substack{\mathbf t\text{ in the short box}\\
\mathcal C,\ \Delta\ne0}}
\ \sum_{\substack{q\mid|\Delta|\\sQ\le q\le2sQ}}
\lambda_q^{-1}.
\]
This is the full majorant in Corollary~\ref{ar:count}, including
tuples whose recovered supports contain no admissible integer point.
Section~\ref{ar:notation} gives
\[
\mathcal C,\ \Delta=0
\quad\Longrightarrow\quad
u=v=t_1=t_2=t_3=t_4=1,
\]
which is excluded by the growing short-variable scales for
sufficiently large \(B\). By \eqref{ar:units},
\[
q\mid|\Delta|
\quad\Longrightarrow\quad
\gcd(q,uvt_1t_2t_3t_4)=1.
\]
The sign change \(L_2\mapsto-L_2\) preserves \(\Lambda_q\) and
\(\lambda_q\). Thus Propositions~\ref{ge:four-proposition}
and~\ref{ge:reuss-proposition}, applied with exponent
\(\varepsilon/2\), give
\begin{align}
\frac{\mathcal G}{B}
\ll_{M,h,\varepsilon} B^{\varepsilon/2}
\biggl(&
s^{-1/2}\sqrt{\frac QH}
+s\frac{\sqrt{HJ}}A
+\frac Q{\sqrt{HJ}}
+\sqrt{\frac JA}
+s^{-1}(C^{-1}+D^{-1})
\biggr),
\label{gl:four}\\
\frac{\mathcal G}{B}
\ll_{M,h,\varepsilon} B^{\varepsilon/2}
\biggl(&
s^{-1/2}\sqrt{\frac QH}
+\frac Q{\sqrt{HJ}}
+s^{-1}\frac{T}{H\sqrt{AJ}}
+s^{-1}(C^{-1}+D^{-1})\notag\\
&+Z_a^{-1}+Z_b^{-1}
+\frac{\sqrt{J/A}}{\sqrt{P_y}}
\biggr).
\label{gl:reuss}
\end{align}
Here \eqref{ge:identities} gives
\[
\frac{Q^2}{\sqrt{HJP_y}}=\frac{\sqrt{HJ}}A,
\qquad
\frac{Q}{\sqrt H\,P_y^{1/4}}=\sqrt{\frac JA}.
\]
The exceptional contributions are already included in
\eqref{gl:reuss}. For every retained block, we have
\[
C^{-1}+D^{-1}+Z_a^{-1}+Z_b^{-1}\ll B^{-\delta},
\qquad
\frac{\sqrt{J/A}}{\sqrt{P_y}}
\le P_y^{-1/2}\le B^{-2\delta}.
\]

Suppose first that \(P_0\ge B^{1+\delta}\).
Since \(A\ge H\ge J\),
\begin{equation*}
\sqrt{\frac QH}
=\left(\frac B{P_0}\right)^{1/6}
  \left(\frac{J^2}{AH}\right)^{1/6}
\le B^{-\delta/6}\quad\text{and}\quad
\frac Q{\sqrt{HJ}}
=\left(\frac B{P_0}\right)^{1/3}
  \left(\frac{HJ}{A^2}\right)^{1/6}
\le B^{-\delta/3}.
\end{equation*}
If \(A/J\ge B^\delta\), then
\[
\frac{\sqrt{HJ}}A
\le\sqrt{\frac JA}
\le B^{-\delta/2},
\]
and we use \eqref{gl:four}. Otherwise,
\[
P_0=AHJ<B^\delta HJ^2\le B^\delta H^3,
\]
so
\[
\frac{T}{H\sqrt{AJ}}
=\left(\frac B{P_0}\right)^{1/3}
 \left(\frac{P_0}{H^3}\right)^{1/6}
\le B^{-\delta/6},
\]
and we use \eqref{gl:reuss}.

Now suppose that \(P_0\le B^{1-\delta}\).
On a block assigned to the Poisson method, we have
$A\le F_0B^{3\delta/4}.$
For complementary pairs \(ij|kl\), using
\(Y_kY_l=F_0Z_{kl}/Z_{ij}\), we obtain
\begin{equation*}
\frac{Z_{ij}}{F_0}
=
\frac{\sqrt{Z_{ij}Z_{kl}/F_0}}{\sqrt{Y_kY_l}}
\le
\frac{\sqrt{A/F_0}}{B^\delta}
\le B^{-5\delta/8}.
\end{equation*}
Consequently,
$
\min_e Z_e,\ \min_iY_i,\ \min_jQ_j
\ge B^{5\delta/8},\quad
\frac{P_0}{B}\le B^{-5\delta/8}
\quad\text{and}\quad
\frac{P_0Z_e}{T^2}
=\frac{Z_e}{F_0}
\le B^{-5\delta/8}
\;\text{for every }e.$
Thus \eqref{po:strict-1}--\eqref{po:strict-2} hold with
\(5\delta/8\) in place of \(\delta\). For fixed \(R\),
Theorem~\ref{po:main} and \(K_{16}(W)\ll_{M,h}1\) yield
\[
N_R(W)
=
\mathfrak S_R\mathcal V_W(B)
+
O_{M,h,R,\delta,\varepsilon}
\!\left(B^{1-5\delta/16+\varepsilon/2}\right).
\]

On the remaining blocks, \(A>F_0B^{3\delta/4}\).
Since \(F_0^3=B^2/P_0\), we have
\[
\frac Q{\sqrt{HJ}}
=\sqrt{\frac{F_0}{A}}
<B^{-3\delta/8},\quad
\sqrt{\frac QH}
=\left(\frac{F_0}{A}\right)^{1/4}
  \left(\frac JH\right)^{1/4}
<B^{-3\delta/16},\]
and
$\frac{P_0}{A^3}
<
\left(\frac{P_0}{B}\right)^2B^{-9\delta/4}
\le B^{-17\delta/4}.$
Hence
\[
\frac{\sqrt{HJ}}A
=
\left(\frac{P_0}{A^3}\right)^{1/2}
<B^{-17\delta/8},
\qquad
\sqrt{\frac JA}
\le
\left(\frac{P_0}{A^3}\right)^{1/4}
<B^{-17\delta/16},
\]
and \eqref{gl:four} applies.

These cases exhaust \eqref{gl:core}. On every block assigned
to the two-dimensional method,
\[
\frac{\mathcal G}{B}
\ll_{M,h,\varepsilon}
B^{\varepsilon/2+\kappa-\delta/6}.
\]
Corollary~\ref{ar:count}, with exponent \(\varepsilon/2\), and
\(\|\Phi\|_{1,*}\|\psi\|_\infty\ll_{M,h}1\) give the geometric error
per block and band as
\[
B^{\varepsilon/2}s^{-1}\mathcal G
\ll_{M,h,\varepsilon}B^{1+\varepsilon+2\kappa-\delta/6}.
\]
The local Poisson error is
\(O_{M,h,R,\delta,\varepsilon}(B^{1-5\delta/16+\varepsilon/2})\).
Both are summed with the arithmetic error below.

\subsection{The arithmetic error on those same blocks}

Consider the blocks assigned to
the two-dimensional method, with
$0<\kappa,\varepsilon<\frac{\delta}{100}$ and
$B^{-\kappa}\le s\ll_{M,h}1.$
For \(\mathcal R\) defined in \eqref{ar:relative-error},
the core conditions give
\[
P_y\ge B^{4\delta},\qquad
\min(U,V,Y_1,Y_2,Y_3,Y_4,Q,Q')\ge B^\delta.
\]

Outside the deleted band \eqref{gl:arithmetic-band}, we have
$\frac{HJ}{T}\le B^{-\delta}$ or
$\frac{HJ}{T}\ge B^\delta.$
In the first range,
\[
\frac Q{\sqrt{P_y}}
=\frac{T/A}{T^2/P_0}
=\frac{HJ}{T}
\le B^{-\delta}.
\]
In the second range, \(A\ge H\) gives
\[
\frac Q{\sqrt{JP_y}}
=\frac{H\sqrt J}{T}
\le\frac{\sqrt{AHJ}}T
=P_y^{-1/4}
\le B^{-\delta}\quad\text{and}\quad
\sqrt{\frac QJ}
=\sqrt{\frac T{AJ}}
\le\sqrt{\frac T{HJ}}
\le B^{-\delta/2}.
\]
Furthermore,
\[
\begin{aligned}
P_y^{-1/3}
+\min(U,V,Y_1,Y_2,Y_3,Y_4)^{-1/2}
+Q^{-1}+{Q'}^{-1}
\le
B^{-4\delta/3}+B^{-\delta/2}+2B^{-\delta}.
\end{aligned}
\]
Consequently, on all retained two-dimensional blocks,
we have
$\mathcal R\ll B^{-\delta/2}.$
The weights constructed above satisfy
\[
\eta=s,\qquad
L=sQ\ge B^{\delta-\kappa}\ge1,\qquad
\|\Phi\|_{6,*}\|\psi\|_{C^1}\ll_{M,h}1.
\]
Apply the fixed-\(\eta_{\max}\) version of Theorem~\ref{ar:main},
with \(\eta_{\max}\) chosen to cover \(s\ll_{M,h}1\).
Its arithmetic error per block and band is
\[
B^{1+\varepsilon}s^{-7}\mathcal R
\|\Phi\|_{6,*}\|\psi\|_{C^1}
\ll_{M,h,\varepsilon}B^{1+\varepsilon+7\kappa-\delta/2}.
\]
The arithmetic and geometric errors are added. Since
\(0<\varepsilon,\kappa<\delta/100\), we have
\[
\varepsilon+2\kappa-\delta/6<-\delta/10,\qquad
\varepsilon+7\kappa-\delta/2<-\delta/10,\qquad
\varepsilon/2-5\delta/16<-\delta/10.
\]
There are \(O(\log^6(2B))\) blocks in the variables \(z_e\), finitely many sign sectors,
and \(O_{M,h,\kappa}(\log(2B))\) active bands per two-dimensional block. Thus
\[
\begin{array}{c|c|c}
\text{error}&\text{local bound}&\text{total bound}\\
\hline
\text{geometric}&B^{1+\varepsilon+2\kappa-\delta/6}
 &B^{1-\delta/10}\log^7(2B)\\
\text{arithmetic}&B^{1+\varepsilon+7\kappa-\delta/2}
 &B^{1-\delta/10}\log^7(2B)\\
\text{Poisson}&B^{1-5\delta/16+\varepsilon/2}
 &B^{1-\delta/10}\log^6(2B)
\end{array}
\]
All totals are \(o(B(\log B)^6)\). Their constants may depend on
\(M,h,\delta,\kappa,\varepsilon\), and, in the Poisson row, on fixed
\(R\). Poisson blocks require no additional band summation.

The power saving above concerns the counting errors on the retained
blocks, with all auxiliary parameters fixed: their sum is
\[
 O_{M,h,\delta,\kappa,\varepsilon,R}
 \!\left(B^{1-\delta/10}\log^7(2B)\right).
\]
After the height limit, the limsup inequalities
\eqref{gl:limit-error} still contain the normalized truncation costs
\(C_{M,h}\delta+CR^{-5/288}+CM^{-1/2}\).
The continuous-volume and height-envelope limits are justified by
dominated convergence, without extracting a quantitative rate.
Moreover, the saving exponent \(\delta/10\) tends to zero when
\(\delta\to0\), and the implied constants may depend on the other
auxiliary parameters. Thus the local power saving does not establish
a global remainder \(O(B^{1-c})\) for any fixed \(c>0\).
The ordered limits below yield the asserted
\(o(B(\log B)^6)\) remainder.

\subsection{A single main measure}

By \eqref{ar:coarea},
\[
\sum_{\nu=\pm1}\int_0^\infty
\psi(q/L)J_\nu(q)\frac{dq}{q^2}
=
\int
\Phi\,\psi\!\left(\frac{|L_1|}{aL}\right)
\delta_0\!\left(\sum_iA_i(z)y_i\right)
\,d^4y\,d^6z.
\]
For the Poisson weights, \eqref{po:dirac-measure} gives
\[
\mathcal V_W(B)
=
\int
W\!\left((z_e/Z_e)_e,(A_i(z)y_i/T_z)_i\right)
\delta_0\!\left(\sum_iA_i(z)y_i\right)
\,d^4y\,d^6z.
\]
Thus both expressions use the same Leray measure, with their
original partition weights. The additional factor
\(\psi(|L_1|/(aL))\) remains on the two-dimensional blocks.

Write \(\sum_{\mathrm{Pois}}\) for summation over the retained
Poisson weights, with all sign sectors included. Fix \(R\ge2\).
By Theorem~\ref{po:main} and the preceding error estimates,
\[
\sum_{\mathrm{Pois}}N_R(W)
=
\mathfrak S_R\sum_{\mathrm{Pois}}\mathcal V_W(B)
+
o_{M,h,\delta,R}\!\left(B(\log B)^6\right).
\]
Equation~\eqref{po:S-limit} gives
\begin{equation}\label{gl:finite-product}
\mathfrak S_R=\mathfrak S+O(R^{-1}).
\end{equation}

The nonnegative partition weights have total mass at most one,
and their union satisfies \(H_{\mathrm{raw}}\le(1+h)B\le2B\).
By Corollary~\ref{po:free-tail},
\begin{equation}\label{gl:sieve-error}
\begin{aligned}
0
\le
\sum_{\mathrm{Pois}}
\bigl(N_R(W)-N_{\mathrm{full}}(W)\bigr)
\le
\mathcal T_{\mathrm{free},>R}(2B)
\ll B\log^6(2B)R^{-5/288}.
\end{aligned}
\end{equation}
This estimate applies once to the union of the supports and is
uniform in \(B,R,M,h\). All heights are unreduced.

For sufficiently large \(B\), the real Poisson supports satisfy
\(z_e,|y_i|\ge1\). Hence
\[
P(z)^2
\le
P(z)^2\left(\prod_i|y_i|\right)^3
=
\prod_i|x_{\mathrm{raw},i}|
\le(2B)^4.
\]
In the variables of \eqref{co:Jacobian}, this implies, for \(B\ge4\),
\[
u_e\ge0,\qquad
\sum_eu_e\le2+\frac{2\log2}{\log B}\le3,
\qquad
\max_i|\Phi_i(t)|\le2.
\]
By \eqref{co:real-domain}--\eqref{co:real-tail} and cubic scaling,
the last real domain has volume \(2J_\infty\). Consequently,
\[
0\le
\sum_{\mathrm{Pois}}\mathcal V_W(B)
\le
2J_\infty B(\log B)^6
\int_{\substack{u_e\ge0\\\sum_eu_e\le3}}d^6u
\ll B(\log B)^6.
\]
Together with \eqref{gl:finite-product}, this gives
\[
\left|
(\mathfrak S_R-\mathfrak S)
\sum_{\mathrm{Pois}}\mathcal V_W(B)
\right|
\ll B(\log B)^6R^{-1}.
\]
For sufficiently large $B$, the hypotheses of
Corollary~\ref{po:core-main} hold, with core
exponent $5\delta/8$ and $K_{16}(W)\ll_{M,h}1$.
Consequently,
\[
\sum_{\mathrm{Pois}}N_{\mathrm{full}}(W)
=\mathfrak S\sum_{\mathrm{Pois}}\mathcal V_W(B)
+o\!\left(B(\log B)^6\right).
\]
The preceding bounds retain the uniform quantitative costs
$O(B\log^6(2B)R^{-5/288})$ and $O(B(\log B)^6R^{-1})$.
Both methods therefore contribute the common measure
\begin{equation}\label{gl:common-measure}
\mathfrak S\,
\delta_0\!\left(\sum_iA_i(z)y_i\right)
\,d^4y\,d^6z,
\end{equation}
multiplied by their respective partition weights. Comparison with
the untruncated continuous domain is deferred to the next subsection.

By Remark~\ref{ar:other-lines}, the canonical tuples in the core
satisfy all three pair-sum exclusions. The pairing-plane correction
is already included in Theorem~\ref{po:main}; by
\eqref{po:line-removal}, its sum is
$O_{M,h}(B^{1-\delta}\log^6(2B))=o(B(\log B)^6)$.
These planes have zero Leray measure.

Finally, the sign \(\nu=\operatorname{sgn}L_1\) in
\eqref{ar:coarea} is uniquely determined; it introduces no
factor \(1/2\). Both overall signs are included throughout, and
\[
(z,y)\longmapsto(z,-y)
\quad\Longrightarrow\quad
x_{\mathrm{raw}}\longmapsto-x_{\mathrm{raw}}.
\]
The factor \(1/2\) is applied once, in the final passage from
canonical signed representatives to projective points.

\subsection{Continuous boundaries and the final passage to the limit}

In the coordinates of \eqref{co:Jacobian},
\[
\frac{x_{\mathrm{raw},i}}B=\Phi_i(t),\qquad
\delta_0\!\left(\sum_iA_i(z)y_i\right)d^4y\,d^6z
=
B(\log B)^6\,d^6u\,d^3t.
\]
Recall the height envelopes \(w_{M,h}^{\pm}\) defined above.
Let \(\mathcal V_{M,h}^{\pm}(B)\) be obtained from \eqref{co:V}
by replacing only the height indicator by \(w_{M,h}^{\pm}(t)\),
retaining
\[
z_e\ge1,\qquad |y_i|\ge1,\qquad
|q_{ij|kl}|\ge1
\quad\text{for all three pairings}.
\]
In particular, the upper envelope is not intersected with
\(H_{\mathrm{raw}}\le B\).

On either height support,
\[
\frac1{M(1+h)}\le|\Phi_i(t)|\le1+h,
\qquad
|t_i|
=
\frac{\bigl(\prod_j|\Phi_j(t)|\bigr)^{1/3}}
{|\Phi_i(t)|}
\asymp_{M,h}1.
\]
For the full smoothed domain, and for the retained local supports
when \(B\) is sufficiently large,
\[
P(z)^2\le\prod_i|x_{\mathrm{raw},i}|\le(2B)^4,
\qquad
u_e\ge0,\qquad
\sum_eu_e
\le2+\frac{2\log2}{\log B}
\le3
\quad(B\ge4).
\]
Moreover,
\[
|q_{ij|kl}|
=
\frac{T_z}{z_{ij}z_{kl}}|t_i+t_j|
\asymp
\frac{T}{Z_{ij}Z_{kl}}|t_i+t_j|.
\]
Thus every block meeting the full smoothed domain satisfies
\[
Z_e\gg1,\qquad
Y_i\gg_{M,h}1,\qquad
Q_j\gg_{M,h}1.
\]
Since
\[
\frac{\log Z_e}{\log B}
=
u_e+O\!\left(\frac1{\log B}\right),
\]
the parameter ranges excluded in
\eqref{gl:core}--\eqref{gl:arithmetic-band}
are contained in finitely many regions between pairs of
parallel hyperplanes with fixed nonzero normal vectors.
The distance between the hyperplanes in each pair is
\(O_{M,h}(\delta+1/\log B)\), with fixed
nonzero normal vectors. Fubini's theorem on the bounded
\(u\)-simplex gives
\[
\operatorname{Vol}_6(\text{deleted logarithmic slabs})
\ll_{M,h}\delta+\frac1{\log B}.
\]

On retained local supports, \(z_e,|y_i|\ge1\) for sufficiently
large \(B\). For each of the three pairings,
\[
|q_{ij|kl}|<1
\quad\Longrightarrow\quad
|t_i+t_j|\ll_{M,h}B^{-\delta}.
\]
For the selected pairing, the omitted near-line weight,
including the smooth transition, satisfies
\[
|q_{ij|kl}|\ll B^{-\kappa}Q
\quad\Longrightarrow\quad
|t_i+t_j|\ll_{M,h}B^{-\kappa}.
\]
Fubini's theorem in coordinates containing \(t_i+t_j\) gives
\[
\operatorname{Vol}_3
\{t\in\operatorname{supp}w_{M,h}^{+}:|t_i+t_j|\le r\}
\ll_{M,h}r
\qquad(0<r\le1).
\]
Here the assembled main term is the sum of the
two-dimensional and Poisson continuous integrals, without the
common arithmetic factor \(\mathfrak S\), and with their
assigned partition weights and, on the two-dimensional blocks,
all retained \(q\)-bands.
These nonnegative weights have total mass at most one.
Thus replacing this sum by \(\mathcal V_{M,h}^{\pm}(B)\) has total cost
\begin{equation}\label{gl:continuous-boundary}
O_{M,h}\!\left(
B(\log B)^6
\left\{
\delta+\frac1{\log B}+B^{-\kappa}+B^{-\delta}
\right\}
\right).
\end{equation}
This includes insertion of all missing auxiliary lower bounds;
no factor counting the dyadic blocks is introduced.

Apply the preceding counting estimates to the two envelopes.
Proposition~\ref{ta:height-tail} controls the integer height tail;
Proposition~\ref{ta:boundary} and Theorem~\ref{ta:near-line}
control the remaining integer deletions. Together with
\eqref{gl:sieve-error} and \eqref{gl:continuous-boundary}, they give
\begin{equation}\label{gl:limit-error}
\begin{aligned}
\limsup_{B\to\infty}
\frac{2N_U(B)-\mathfrak S\mathcal V_{M,h}^{+}(B)}
{B(\log B)^6}
&\le C_{M,h}\delta+CR^{-5/288}+CM^{-1/2},\\
\limsup_{B\to\infty}
\frac{\mathfrak S\mathcal V_{M,h}^{-}(B)-2N_U(B)}
{B(\log B)^6}
&\le C_{M,h}\delta+CR^{-5/288}.
\end{aligned}
\end{equation}
Here \(C\) is absolute and \(C_{M,h}\) is independent of
\(\delta,R,B\). The term \(R^{-1}\) from
\eqref{gl:finite-product} is absorbed into \(R^{-5/288}\).

The proof of Proposition~\ref{co:continuous} applies with the
bounded weight $w_{M,h}^{\pm}$.
On $\{u_e\ge0\}$, the product of the four indicators
$|y_i|\ge1$ and the three indicators $|q_{ij|kl}|\ge1$
converges almost everywhere to $\mathbf1_\Omega(u)$.
An integrable majorant, independent of \(M\) and \(0<h<1\), is
\[
\mathbf1_{\{u_e\ge0,\ \sum_eu_e\le3\}}
\mathbf1_{\{\max_i|\Phi_i(t)|\le2\}}.
\]
Indeed, \eqref{co:real-tail} and cubic scaling give
\[
\int_{\substack{u_e\ge0\\\sum_eu_e\le3}}d^6u
\int_{\{\max_i|\Phi_i(t)|\le2\}}d^3t
=
\frac{3^6}{6!}\,2J_\infty
<\infty.
\]
Dominated convergence and \eqref{co:volume} yield
\[
\lim_{B\to\infty}
\frac{\mathcal V_{M,h}^{\pm}(B)}{B(\log B)^6}
=
\operatorname{Vol}_6(\Omega)
\int_{\mathbb R^3}w_{M,h}^{\pm}(t)\,d^3t
=
\frac1{640}
\int_{\mathbb R^3}w_{M,h}^{\pm}(t)\,d^3t.
\]
For fixed \(M\), outside polynomial zero sets,
\[
w_{M,h}^{\pm}(t)
\longrightarrow
\mathbf1_{\{1/M\le|\Phi_i(t)|\le1\ \text{for all }i\}}
\qquad(h\to0).
\]
These indicators increase almost everywhere to
\(\mathbf1_{\mathcal D}(t)\) as \(M\to\infty\). Hence
\[
\lim_{M\to\infty}\lim_{h\to0}\lim_{B\to\infty}
\frac{\mathfrak S\mathcal V_{M,h}^{\pm}(B)}
{B(\log B)^6}
=
\frac{\mathfrak S J_\infty}{640}.
\]

Take limits in \eqref{gl:limit-error} in the order
\[
B\to\infty,\qquad
R\to\infty,\qquad
\delta\to0,\qquad
h\to0,\qquad
M\to\infty,
\]
choosing \(0<\kappa,\varepsilon<\delta/100\) before the first limit.
For fixed \(M,h\), the first three limits give
\[
\begin{aligned}
\frac{\mathfrak S}{640}
\int_{\mathbb R^3}w_{M,h}^{-}(t)\,d^3t
&\le
\liminf_{B\to\infty}
\frac{2N_U(B)}{B(\log B)^6}
\le
\limsup_{B\to\infty}
\frac{2N_U(B)}{B(\log B)^6}\\
&\le
\frac{\mathfrak S}{640}
\int_{\mathbb R^3}w_{M,h}^{+}(t)\,d^3t
+CM^{-1/2}.
\end{aligned}
\]
Letting \(h\to0\) and then \(M\to\infty\) in these inequalities,
and using Corollary~\ref{co:constant}, gives
\[
\lim_{B\to\infty}\frac{N_U(B)}{B(\log B)^6}
=\frac{\mathfrak S J_\infty}{1280}
=\frac{\mathfrak S I}{80}
=c_{\mathrm{Peyre}}.
\]
Equivalently, by Proposition~\ref{co:continuous},
\begin{equation}\label{gl:comparison}
2N_U(B)=\mathfrak S\mathcal V(B)+o\bigl(B(\log B)^6\bigr).
\end{equation}

\subsection*{Acknowledgements}

This work is supported by the National Natural Science Foundation of China
(Grant Nos.\ 12571002, 12501007 and 12526407), Young Elite Scientists
Sponsorship Program by Henan Association for Science and Technology
(Grant No.\ 2026HYTP001), High-level Talent Research Start-up Project
Funding of Henan Academy of Sciences (Grant Nos.\ 241819033, 252019083 and
242019104) and High-level Achievement Award and Cultivation Project Funding
of Henan Academy of Sciences (Grant Nos.\ 20262319007, 20262319001 and
20262319003).

\subsection*{AI Declaration}
The author used ChatGPT to assist with drafting, language editing,
and the examination of mathematical arguments. The author is
responsible for verifying all statements, proofs and references
in this manuscript.

\bigskip
\noindent
\begin{minipage}{\textwidth}
\footnotesize
\textsc{Institute of Mathematics, Henan Academy of Sciences,
Zhengzhou, Henan 450046, P.~R.~China}\par
\noindent\textit{E-mail address:}
{\hypersetup{urlcolor=black}%
\href{mailto:zhaoxiaodong@hnas.ac.cn}{\texttt{zhaoxiaodong@hnas.ac.cn}}}
\end{minipage}

\end{document}